\documentclass[12pt,amssymb]{amsart}
\usepackage{latexsym}
\begin{document}
\newtheorem{propo}{Proposition}[section]
\newtheorem{defi}[propo]{Definition}
\newtheorem{conje}[propo]{Conjecture}
\newtheorem{lemma}[propo]{Lemma}
\newtheorem{corol}[propo]{Corollary}
\newtheorem{theor}[propo]{Theorem}
\newtheorem{examp}[propo]{Example}
\newtheorem{remar}[propo]{Remark}
\newtheorem{ques}[propo]{Question}
\newtheorem{prob}[propo]{Problem}

\newcommand{\Gal}{\mathop{\rm {Gal}}\nolimits}
\newcommand{\Aut}{\mathop{\rm {Aut}}\nolimits}
\newcommand{\Out}{\mathop{\rm {Out}}\nolimits}
\newcommand{\Mult}{{\rm {Mult}}}
\newcommand{\Inn}{{\rm {Inn}}}
\newcommand{\Irr}{\mathop{\rm {Irr}}\nolimits}
\newcommand{\IBR}{{\rm {IBr}}}
\newcommand{\IBRL}{{\rm {IBr}}_{\ell}}
\newcommand{\Id}{{\rm Id}}
\def\Ker{\mathop{\mathrm{ Ker}}\nolimits}
\def\id{\mathop{\mathrm{ id}}\nolimits}
\renewcommand{\Im}{{\rm {Im}}}
\newcommand{\Ind}{{\rm {Ind}}}
\newcommand{\diag}{{\rm {diag}}}
 \newcommand{\mar}{\marginpar}
\newcommand{\soc}{{\rm {soc}}}
\newcommand{\End}{\mathop{\rm {End}}\nolimits}
\newcommand{\sol}{{\rm {sol}}}
\newcommand{\Hom}{{\rm {Hom}}}
\def\rank{\mathop{\mathrm{ rank}}\nolimits}
\def\ch{\mathop{\mathrm{ char}}\nolimits}
\newcommand{\Syl}{\mathop{\rm {Syl}}\nolimits}
\newcommand{\Tr}{{\rm {Tr}}\nolimits}
\newcommand{\tr}{{\rm {tr}}\nolimits}

\newcommand{\sss}{\subset}

\newcommand{\Spec}{\mathop{\rm Spec}\nolimits}
\newcommand{\spec}{\mathop{\rm Spec}\nolimits}
\newcommand{\ad}{\mathop{\rm ad}\nolimits}
\newcommand{\Sym}{\mathop{\rm Sym}\nolimits}

\newcommand{\Char}{\mathop{\rm char}\nolimits}
\newcommand{\vtri}{\rhd}
\newcommand{\pr}{\mathop{\rm pr}\nolimits}
\newcommand{\Rad}{\mathop{\rm Rad}\nolimits}
\newcommand{\codim}{\mathop{\rm codim}\nolimits}
\newcommand{\ind}{\mathop{\rm ind}\nolimits}
\newcommand{\CC}{{\mathbb C}}
\newcommand{\RR}{{\mathbb R}}
\newcommand{\QQ}{{\mathbb Q}}
\newcommand{\ZZ}{{\mathbb Z}}
\newcommand{\GG}{{\mathbb G}}
\newcommand{\SSS}{{\mathbb S}}
\newcommand{\AAA}{{\mathbb A}}
\newcommand{\FF}{{\mathbb F}}
\newcommand{\FQ}{{\mathbb F}_{q}}
\newcommand{\EC}{{\mathcal E}}
\newcommand{\GC}{{\mathcal G}}
\newcommand{\GCF}{\GC^{F}}
\newcommand{\GCD}{\GC^{*}}
\newcommand{\GCDF}{(\GCD)^{\FD}}
\newcommand{\FD}{F^{*}}
\newcommand{\HC}{{\mathcal H}}
\newcommand{\HCF}{\HC^{F}}
\newcommand{\HCD}{\HC^{*}}
\newcommand{\HCDF}{(\HCD)^{\FD}}
\newcommand{\OC}{{\mathcal O}}
\newcommand{\LC}{{\mathcal L}}
\newcommand{\LCF}{\LC^{F}}
\newcommand{\TC}{{\mathcal T}}
\newcommand{\TCF}{\TC^{F}}
\newcommand{\DC}{{\mathcal D}}
\newcommand{\DCF}{\DC^{F}}
\newcommand{\CL}{{\mathcal C}}
\newcommand{\CLQ}{\CL(q)}
\newcommand{\CLL}{{\rm {Cl}}_{\ell'}}
\newcommand{\HAA}{\hat{\alpha}}
\newcommand{\HAB}{\hat{\beta}}
\newcommand{\HAC}{\hat{\gamma}}
\newcommand{\HAZ}{\hat{\zeta}}
\newcommand{\ZC}{{\mathcal Z}}
\newcommand{\ZCF}{\ZC^{F}}
\newcommand{\XC}{{\mathcal X}}
\newcommand{\ECL}{{\mathcal E}_{\ell}}
\newcommand{\ECS}{\EC(\GCF,(s))}
\newcommand{\ECST}{\EC(\GCF,(st))}
\newcommand{\SC}{{\mathcal S}}
\newcommand{\MC}{{\mathcal M}}
\newcommand{\dvg}{d_{V}(g)}
\newcommand{\bg}{\bar{g}}
\newcommand{\tG}{\tilde{G}}
\newcommand{\bga}{\bar{g}_{A}}
\newcommand{\bgb}{\bar{g}_{B}}
\newcommand{\ga}{g_{A}}
\newcommand{\gb}{g_{B}}
\newcommand{\bfa}{\bar{f}_{A}}
\newcommand{\bfb}{\bar{f}_{B}}
\newcommand{\fa}{f_{A}}
\newcommand{\fb}{f_{B}}
\newcommand{\bj}{{\bf j}}
\newcommand{\RGT}{R_{\TC,\,\theta}}
\newcommand{\la}{\langle}
\newcommand{\ra}{\rangle}
\newcommand{\eps}{\epsilon}
\newcommand{\varep}{\varepsilon}
\newcommand{\lam}{\lambda}
\newcommand{\om}{\omega}
\newcommand{\se}{\subseteq}
\renewcommand{\mod}{\bmod \,}
\newcommand{\ep}{\varepsilon}
\newcommand{\epf}{\hfill $\Box$}
\newcommand{\al}{\alpha} 
\newcommand{\ta}{\hspace{0.5mm}^{2}\hspace*{-0.8mm}}
\newcommand{\tb}{\hspace{0.5mm}^{3}\hspace*{-0.8mm}}
\newcommand{\tn}{\hspace{0.5mm}^{t}\hspace*{-0.8mm}}

\marginparsep-0.5cm
\renewcommand{\thefootnote}{\fnsymbol{footnote}}
\footnotesep6.5pt

\title[Hall-Higman type theorems]
{Hall-Higman type theorems for semisimple elements of finite classical groups}
\author{Pham Huu Tiep}
\address{Department of Mathematics,
University of Florida, Gainesville, FL 32611-8105, U.S.A.}
\email{tiep@math.ufl.edu}
\author[A. E. Zalesski\u\i]{A. E. Zalesski\u\i}
\address{School of Mathematics, University of
East Anglia, Norwich NR4 7TJ, United Kingdom}
\email{a.zalesskii@uea.ac.uk}
\thanks{The first author gratefully acknowledges the support of the NSF 
(grant DMS-0600967), the NSA (grant H98230-04-0066), and the EPSRC (grant GR/S56511/01). 
The second author acknowledges a partial support from
the Levehulme Trust (grant EM/2006/0030).
The authors also thank F. L\"ubeck for helping them with some computer calculation.}


\maketitle

\section{Introduction}

Let $G$ be a finite group. For any element $g \in G$ and any irreducible 
representation $\Theta$ of $G$ over an algebraically closed 
field $\FF$ of characteristic $\ell$, let $\deg (\Theta(g))$ denote the degree 
of the minimal polynomial of the matrix $\Theta(g)$. For $g \in G$, let $o(g)$ denote the order of 
$g$ modulo $Z(G)$. Clearly, $\deg(\Theta(g)) \leq o(g)$. 
Furthermore, if the characteristic of $\FF$ is coprime to the order $|g|$ of 
$g$, then $\deg(\Theta(g))$ is just the number of distinct eigenvalues of $\Theta(g)$. 
On the other hand, if $\Char(\FF) = \ell > 0$ and $g$ is an $\ell$-element, then 
$\deg(\Theta(g))$ is the largest size of Jordan blocks of $\Theta(g)$. 

There are many applications demonstrating the importance of knowing the eigenvalues and minimal 
polynomials of elements of linear groups, in other words, of group elements in finite 
dimensional representations. Investigations on this topic go back at least to the classical results of 
Blichfeldt \cite{B}. In 1956, a fundamental result was obtained by Hall and Higman \cite{HH} 
which describes the minimal polynomial of any $p$-element of a finite $p$-solvable group in a 
$p$-modular representation. The Hall-Higman Theorem led to various reductions for the Restricted
Burnside Problem, cf. \cite{HH}. It also proved to be invaluable for the development of linear 
methods in group theory, cf. \cite{HB}. Later on, many important results were 
established by Thompson \cite{Th}, Shult \cite{Sh}, Berger \cite{Be}, Robinson \cite{R}, and others. 
To a large extent, these results can be viewed as contributions to the following general problem:

\begin{prob}\label{min} 
{\sl Determine all possible values for $\deg(\Theta(g))$, and if possible, all triples $(G,\Theta,g)$ 
with $\deg(\Theta(g)) < o(g)$, in the first instance under the condition that $o(g)$ is a $p$-power.}
\end{prob} 

The celebrated theorem of Hall-Higman \cite{HH} is concerned with the case where $G$ is 
$p$-solvable, $p = \ell$ and $O_{p}(G) = 1$. Work of Thompson and later of Ho on quadratic modules, 
and also some more recent work along the lines of Meierfrankenfeld-Stellmacher-Stroth's
program on a third-generation proof of the classification theorem for finite simple groups,
are devoted to the case $o(g) = p = \ell$. 

The priority in studying Problem \ref{min} should be given to 
the groups $G$ that are ``close'' to be simple. Already for 
the alternating groups, Problem \ref{min} turns out to be very challenging: let us mention that 
this problem, under the assumptions that $G = \AAA_{p}$ and $g$ is a $p$-cycle, has been solved only in 
the case where $p > \ell > (p-1)/2$ by work of Thompson (unpublished manuscript), and 
$\ell \in \{0,p\}$ by work of
Kleshchev and the second author \cite{KZ}. We consider the case $G \in Lie(r)$, a finite group 
of Lie type in characteristic $r$. The {\it defining} characteristic representations (that is, 
$r = \ell$) can be treated rather efficiently using  
representation theory of algebraic groups. The main tool to handle semisimple elements $g$ is
the theory of weights together with Premet's theorem \cite{P}; however, it is impossible to 
provide an explicit solution to the problem. This is also true for the behaviour of  unipotent elements in defining characteristic representations, cf. \cite{Su}
for results of asymptotic nature. We  focus on {\it cross 
characteristic} representations, i.e. $r \neq \ell$. If $r = p$ (that is, $g$ is {\it unipotent}),
then the problem is basically solved in \cite{Z1} and \cite{DZ2}. The case where $G$ is classical and $g$ is contained in a proper parabolic subgroup is settled in \cite{DZ1}. For $\ell =0$ and $g$ of order $p$, the cases with minimum polynomial of $g$ of degree less than $p$ have been 
determined in \cite{Z4} and \cite{Z5}. Also, substantial 
results in the case $\ell = p$ have been obtained in \cite{Z3}. However, the methods developed in these papers do not work in some key cases.
Consider for example an element $g$ of order $2^{13}-1=8191$ in $G=SL_{13}(2)$ and any
nontrivial irreducible representation $\Theta$ of $G$ in characteristic $\ell $ dividing $|G|$ for $2<\ell <8191$. All the 
currently available general results yield only the bound $\deg(\Theta(g)) \geq 13$. However, various heuristic 
arguments seem to indicate that a much stronger bound, $\deg(\Theta(g)) \geq 8189$, should hold true. 
More generally, if $g \in G$ has order $p^{a}$, then heuristic arguments
lead to the belief that one should be able to prove that $\deg(\Theta(g)) \geq p^{a-1}(p-1)-1$, and 
in fact this bound should be the optimal bound. 

\smallskip
The goal of the present paper is to prove this optimal bound for semisimple elements of 
finite classical groups in cross characteristic representations.    
The exceptional groups of Lie type will be 
addressed in a sequel of the paper. Let $\MC$ denote the collection of 
finite classical groups with exceptional Schur multiplier, that is, the collection
of finite groups $G$ with $\soc(G/Z(G))$, the socle of $G/Z(G)$, equal to 
$PSL_{2}(4)$, $PSL_{3}(2)$, $PSL_{4}(2)$, $PSU_{4}(2)$, $Sp_{6}(2)$, $PSL_{2}(9)$, 
$PSL_{3}(4)$, $PSU_{4}(3)$, $SU_{6}(2)$, $\Omega_{7}(3)$, $\Omega^{+}_{8}(2)$.      

\smallskip
The first main result of the paper is the following theorem:

\begin{theor}\label{main1}
{\sl Let $G$ be a finite quasi-simple classical group, not belonging to $\MC$, and let $\Theta$ 
be a nontrivial irreducible representation of $G$ over a field of cross characteristic
$\ell$. Let $g \in G \setminus Z(G)$ be a semisimple element of prime order $p$. 
Then one of the following holds.

{\rm (i)} $\deg(\Theta(g)) = p$. 

{\rm (ii)} $p > 2$, $\deg(\Theta(g)) = p-1$ and Sylow $p$-subgroups of $G/Z(G)$ are cyclic.
Moreover, the conclusions of Proposition \ref{p-cyclic} hold with $a = 1$.

{\rm (iii)} $p = q+1$ is a Fermat prime, $G/Z(G) = PSU_{n}(q)$, 
$g^{q+1} = 1$, $\Theta$ is a Weil
representation, and $\deg (\Theta(g)) = p-1$. Furthermore, either $n \leq 3$, or $g$ 
corresponds to a pseudoreflection in $GU_{n}(q)$ and $(n,q+1) = 1$. 

{\rm (iv)} $o(g) = p = (q^{n}-1)/(q-1)$, $G/Z(G) = PSL_{n}(q)$, $n$ a prime, $\Theta$ is a Weil 
representation of degree $p-1$ or $p-2$, and $\dim(\Theta) = \deg(\Theta(g))$. Furthermore, 
either $n > 2$ and $\ell = p$, or $n = 2$ and $q$ is even.}
\end{theor}

For applications, as well as for the treatment of exceptional groups of Lie type, it is 
helpful to expand the class of groups of Lie type considered in Theorem \ref{main1}.  
For a prime power $q$, let $\CLQ$ be the list of insoluble groups $G$ of the 
following form: $SL^{\eps}_{n}(q) \leq G \leq GL^{\eps}_{n}(q)$ with $n \geq 2$,
where $\eps = +$ for $SL$ and $\eps = -$ for $SU$; 
$Sp_{2n}(q) \leq G \leq CSp_{2n}(q)$ with $n \geq 2$;
$Spin(V) \leq G \leq \Gamma^{+}(V)$ with $V = \FF^{2n+1}_{q}$, $n \geq 3$ and $q$ odd;
$Spin(V) \leq G \leq \Gamma(V)$ with $V = \FF^{2n}_{q}$ and $n \geq 4$; cf. \cite{TZ5} for the 
definition of the Clifford group $\Gamma(V)$ and $\Gamma^{+}(V)$. Setting in addition 
$V = \FF_{q}^{n}$ in the case of $GL_{n}(q)$, $V = \FF_{q^{2}}^{n}$ in the case of 
$GU_{n}(q)$, and $V = \FF_{q}^{2n}$ in the case of $CSp_{2n}(q)$, we will refer to $V$ as
the {\it natural module} for $G$. We say an element $g \in G$ is {\it irreducible} (on $V$),
if the induced action of $g$ on $V$ is irreducible. Furthermore, $g$ is called a 
{\it pseudoreflection} if the action of $g$ on $V$ can be represented by a diagonal matrix 
$\diag(\al,\beta, \beta, \ldots ,\beta)$ for some $\al \neq \beta$. The $\ell$-modular {\it Weil 
representations} are defined to be the composition factors of degree $> 1$ of the $\ell$-modular
reduction of complex Weil representations \cite{S}. 

\smallskip
The second main result of the paper is the following theorem:

\begin{theor}\label{main2}
{\sl Let $G \in \CLQ$ for a prime power $q$, $p$ a prime not dividing $q$, and $S := \soc(G/Z(G))$. 
Let $\Theta$ be an irreducible representation of $G$ over a field of characteristic $\ell$ coprime 
to $q$, of degree $> 1$. Let $g \in G$ be a semisimple element of prime power order $p^{a} > 1$ modulo 
$Z(G)$. Assume in addition that $\dim(\Theta) > 2$ if $G/Z(G) = PGO^{+}_{4n}(q)$ with $n \geq 2$ 
and $q$ odd. Then one of the following holds.

{\rm (i)} $p^{a} \geq \deg(\Theta(g)) > p^{a-1}(p-1)$. 

{\rm (ii)} $p > 2$, $\deg(\Theta(g)) = p^{a-1}(p-1)$ and Sylow $p$-subgroups of $G/Z(G)$ are cyclic.
Furthermore, either $a = 1$, or $\ell \neq p$. Moreover, 
the conclusions of Proposition \ref{p-cyclic} hold.

{\rm (iii)} $S = PSU_{n}(q)$, $o(g) = p = q+1$ is a Fermat prime, $g$ belongs to  
$GU_{1}(q)^{n}$, $\Theta$ is a Weil representation, and $\deg(\Theta(g)) = p-1$. Furthermore, 
either $n \leq 3$, or $g$ is a pseudoreflection in $GU_{n}(q)$.

{\rm (iv)} $a \geq 2$, $S = PSU_{n}(q)$ with $n \equiv 1 (\mod p^{a-1})$, $q+1 = p$ is a 
Fermat prime, $g^{p^{a-1}}$ is a pseudoreflection in $GU_{n}(q)$, and 
$\Theta$ is a Weil representation. Furthermore, either $\deg(\Theta(g)) = p^{a-1}(p-1)$, or 
$(n,p^{a},q) = (4,9,2)$ and $\deg(\Theta(g)) = 5$. 

{\rm (v)} $o(g) = p^{a} = (q^{n}-1)/(q-1)$, $S = PSL_{n}(q)$, $n$ a prime, $\Theta$ is a Weil 
representation of degree $o(g)-1$ or $o(g)-2$, and $\dim(\Theta) = \deg(\Theta(g))$. 
Furthermore, either $n > 2$ and $\ell = o(g) = p$, or $n = 2$ and $q$ is even.}
\end{theor}

Theorem \ref{main2} implies the following consequence, where $\varphi(\cdot)$ denotes the Euler 
function as usual.  

\begin{corol}\label{main3}
{\sl Let $G \in \CLQ$ for a prime power $q$, and let $\Theta$ be an irreducible representation of 
$G$ over a field of characteristic $\ell$ coprime to $q$, of degree $> 1$. Let $g \in G$ be a 
semisimple element of prime power order modulo $Z(G)$. Assume in addition that $\dim(\Theta) > 2$ 
if $G/Z(G) = PGO^{+}_{4n}(q)$ with $n \geq 2$ and $q$ odd. 
Then $\deg(\Theta(g)) \geq \varphi(o(g))-1$, and this bound is best possible.
\hfill $\Box$}
\end{corol}

The optimality of the bound $\deg(\Theta(g)) \geq \varphi(o(g))-1$ follows from Theorem 
\ref{main2}(v).

\smallskip
The information concerning the classical groups belonging to $\MC$ is collected
in Table I, the proof of whose correctness is somewhat ad hoc and omitted for
the sake of brevity. In this table, 
$G/Z(G)$ is a simple finite group of Lie type in characteristic $r$, $G$ is a 
universal  
central extension of $G/Z(G)$ such that $|Z(G)|$ is a multiple of $r$
(that is, $G/Z(G)$ has an exceptional Schur multiplier). Let $Z_{0}$ be 
a Sylow $r$-subgroup of $Z(G)$. In the table, we list all the irreducible 
$\FF G$-representations $\Theta$ such that $\deg (\Theta (g))< o(g)$ and 
$\Theta (Z_{0})\neq 1$. In fact, we also include in the table a few cases 
with $\Theta (Z_0) = 1$ when we can compute the precise value of $\deg (\Theta (g))$
here but not in our general approach described in the paper. 

\vskip10pt
\begin{figure}
 \begin{center}{{\sc Table I: Exceptional Schur multiplier cases } }

 \vspace{10pt}
\begin{tabular}{|l|c|c|c|c|c|c| }

\hline
$G/Z(G)$& $|\Theta (Z_0)|$ & $\ell\neq r $   & $|g|$ & Name &$\dim (\Theta)$  &
$\deg(\Theta (g))$   \\
\hline
 
$PSL_2(4)$&2&&5&5A,5B &$1<d<5$ &$d$ \\

 & &&3& 3A& $2$ &2 \\
 \hline
 
$PSL_2(9)$&3&5&4,5& 4A,5A,5B& $3$ &$3$ \\
 


\hline
$PSL_3(2)$&2&7&7& 7A,7B&$1<d<7$ &$d$ \\
&&7&3& 3A&$2$ &$2$ \\
&&$\neq 7$&7& 7A,7B& $d=3,4,6$ &$d$ \\
&&$\neq 7$&3&3A& $ 4$ &$2$ \\
 \hline
$PSL_3(4)$&16&&7& 7A,7B&$6 $ &$6 $ \\
\hline
$PSU_4(2)$&2&&3& 3C& $4$ &$2$ \\
&&&5,9& 5A, 9A,9B& $4$ &$4$ \\
&1&&9& 9A,9B&$5$ &$5$ \\
&&$\ell\neq 3$&9& 9A,9B&$6$ &$6$ \\
&1&&9& 9A,9B&$10$ &$7$ \\
&&$\ell=3$&9& 9A,9B&$16$ &$8$ \\
&&$\ell\neq 3$&9&9A,9B& $20$ &$8$ \\
\hline
$PSL_4(2)$&2&&3,5&3A,5A& $8$ &$|g|-1$ \\
\hline
$PSp_6(2)$&2&&3,5,9&3A,5A,9A& $8$ &$|g|-1$ \\
&1&&9& 9A&$7$ &$7$ \\
\hline
$PSU_4(3)$&3&&4,7&4A, 7A,7B& $6$ &$|g|-1$ \\
&&&8& 8A&$6$ &$6$ \\
\hline
$\Omega^+_8(2)$&2&&3,5,9&3A,3B,5A,5B,9B,9C& $8$ &$|g|-1$ \\
\hline
$\Omega_7(3)$ &3&$p=\ell =2$ &4,8 &4A,8A,8B & &$\geq 3$, $\geq 5$\\
\hline
$PSU_6(2)$&2&$p=\ell =3$ &9&9A& $9$ &$ $  $\geq 7$ \\
\hline
\end{tabular}
\end{center}
\end{figure}

\smallskip	
Briefly, our approach can be outlined as follows. The action of semisimple elements $g \in G$ on 
the natural module $V$ of $G$ distinguishes two cases: $g$ is irreducible, resp. reducible on $V$.
The critical case of irreducible elements that are minimal in a certain sense is handled in \S4. 
Deligne-Lusztig theory (cf. \cite{C}, \cite{L}), as well as fundamental results of Brou\'e and Michel
\cite{BM} play a key role in the treatment of this critical case. In fact results proved in \S4
also apply to exceptional groups of Lie type. Results of \S4 and \S5, in particular Corollary 
\ref{slprime} and Theorem \ref{slsup}, form an induction base to handle arbitrary 
irreducible elements, which are treated in Theorem \ref{clas1}. We also use 
results of \cite{DZ1}, \cite{TZ2}, and \cite{GMST} concerning semisimple elements 
that lie in a proper parabolic subgroup.
In turn, the case of reducible elements builds on the irreducible case and 
it is completed in \S6.  

We expect our main results to be useful in a number of applications. One particular application
we have in mind is the recognition of finite linear groups. Our results have also
been used in joint work of Kleshchev and the first author \cite{KT} to classify representations of 
finite general linear groups that are irreducible over proper subgroups.
 
Throughout the paper, $\FF$ is an algebraically closed field of characteristic $\ell$ coprime 
to $q$. If $\chi$ is a class function of $G$ then $\hat{\chi}$ denotes the restriction of $\chi$ 
to the $\ell'$-elements in $G$. If $V$ is an $\FF G$-module then $V^{g} := \{v \in V \mid gv = v\}$ 
is the fixed point subspace for $g \in G$. $X * Y$ denotes a central product of subgroups 
$X, Y \leq G$. If $V$ is any $\FF G$-module, $\dvg$ stands for the
degree of the  minimal polynomial of $g$ acting on $V$.  In what follows, by a 
{\it nondegenerate}, resp., {\it totally singular}, subspace of an orthogonal space we mean any 
subspace $U$ such that the bilinear form associated to the quadratic form is nondegenerate on $U$, 
resp., the quadratic form is zero on $U$. When the (Hermitian, symplectic, or quadratic) form
on $V$ is specified, $I(V)$ denotes the isometry group of the form. 

\section{Some preliminaries}

The following simple observation allows us to focus on elements of prime power order {\it in 
$G$} (not just modulo $Z(G)$):

\begin{lemma}\label{easy}
{\sl Let $G$ be a finite group and $g \in G$ with $o(g) = p^{a}$. Then there is a 
$p$-element $h \in G$ with the following properties:

{\rm (i)} $o(h) = p^{a}$;

{\rm (ii)} Let $V$ be any irreducible $\FF G$-module. Then $d_{V}(g) = d_{V}(h)$. 
Moreover, $g$ is irreducible on $V$ if and only if $h$ is.\\ 
Furthermore, if $G \lhd H$ and $C_{H}(g)G = H$ then for any $\Phi \in \IBRL(G)$ there is
$\Theta \in \IBRL(H)$ such that $\Phi$ is a constituent of $\Theta|_{G}$ and 
$d_{\Phi}(g) = d_{\Theta}(g)$. Moreover, if $Z(G) \leq Z(H)$ then the values of $o(g)$ in $G$ and 
in $H$ are the same.}
\end{lemma}

\begin{proof}
Write $g = ht$ with $h$ the $p$-part of $g$ and $t$ the $p'$-part of $g$. Since $o(g) = p^{a}$,
$t \in Z(G)$. Now it is straighforward to check (i) and (ii). Next assume $G \lhd H$, 
$C_{H}(g)G = H$, and $\Phi \in \IBRL(G)$. Then we can find $\Theta \in \IBRL(H)$ such that $\Phi$ 
is a constituent of $\Theta|_{G}$. By Clifford's theorem, 
$\Theta|_{G} = \oplus^{s}_{i=1}\Phi_{i}$, with $\Phi_{i}(g) = \Phi(x_{i}gx_{i}^{-1})$ for some
$x_{i} \in H$. Since $C_{H}(g)G = H$, we can choose $x_{i} \in C_{H}(g)$. It follows that 
$\Phi_{i}(g) = \Phi(g)$ and so $d_{\Phi}(g) = d_{\Theta}(g)$. Finally, assume $Z(G) \leq Z(H)$.
Then $Z(G) = Z(H) \cap G$. It follows that $g^{m} \in Z(G)$ if and only if $g^{m} \in Z(H)$, whence
the values of $o(g)$ in $G$ and in $H$ are the same. 
\end{proof}

\begin{lemma}\label{zgm} {\rm \cite[Lemma IX.2.7]{HB}} 
{\sl Let $p,r$ be primes and $a,b$ positive integers such that $p^{a}=r^{b}+1$. Then either

{\rm (i)} $p = 2$, $b = 1$, and $r$ is a Mersenne prime, or 

{\rm (ii)} $r = 2, a = 1$, and $p$ is a Fermat prime, or

{\rm (iii)} $p^{a} = 9$.
\hfill $\Box$}
\end{lemma}

We will frequently use the following well-known statement:

\begin{lemma}\label{div}
{\sl Let $p$ be a prime and let $q$ be an integer such that $p|(q-1)$. 
If $p = 2$, assume in addition that $4|(q-1)$. Then $(q^{p}-1)/(q-1) \equiv p (\mod p^{2})$.}  
\end{lemma}

\begin{proof} 
The lemma is obvious if $p = 2$. For $p > 2$, write $q = Ap^{c}+1$ for some integers 
$A,c \geq 1$ with $(p,A) = 1$. Then 
$q^{p}-1 = ABp^{3c} + A^{2}p^{2c+1}(p-1)/2 + Ap^{c+1} = ACp^{2c+1} + Ap^{c+1}$ for some 
integers $B,C$. Hence $(q^{p}-1)/(q-1) = Cp^{c+1} + p$.
\end{proof}

\begin{lemma}\label{fix} {\rm \cite{GT3}}
{\sl Let $V$ be an $\FF G$-module with a submodule $U$ and let $g \in G$. Then for the 
dimensions of the $g$-fixed point subspaces the following inequality holds:
$\max\{\dim(U^{g}),\dim((V/U)^{g})\} \leq \dim(V^{g})$.
\hfill $\Box$}
\end{lemma}

\begin{lemma}\label{free}
{\sl Let $A = \langle g \rangle$ be a cyclic group of order $p$ and let $V$ be a free 
$\FF_{p}A$-module. Let $U$ be any subquotient of the $A$-module $V$. Then 
$d_{U}(g) \geq p\dim(U)/\dim(V)$. In particular, if $\dim(U) > (p-1)\dim(V)/p$ then 
$d_{U}(g) = p$.}
\end{lemma}

\begin{proof}
By the assumption, $\dim(V^{g}) = \dim(V)/p$. By Lemma \ref{fix}, $\dim(U^{g}) \leq \dim(V^{g})$.
Hence $d_{U}(g) \geq \dim(U)/\dim(U^{g}) \geq p\dim(U)/\dim(V)$.   
\end{proof}

\begin{lemma}\label{perm} {\rm \cite[Proposition 2.15]{DZ1}}
{\sl Let $V = V_{1} \oplus \ldots \oplus V_{s}$ be a direct sum of $\FF$-spaces and let 
$g \in GL(V)$ be an element of prime-power order that permutes the $V_{i}$'s transitively. Then 
$d_{V}(g) = s \cdot d_{V_{1}}(g^{s})$.
\hfill $\Box$}
\end{lemma}

\begin{remar}\label{unram} {\rm \cite[Corollary 3.2]{TZ3}}
{\em If $\Phi \in \IBR_{p}(G)$ lifts to a faithful complex representation which is unramified above 
$p$, and $g \in G$ is a $p$-element with $o(g) = p$, then $d_{\Phi}(g) \geq p-1$; in fact,
all Jordan blocks of $\Phi(g)$ are of size $1$, $p-1$, or $p$.}
\end{remar}

\begin{lemma}\label{equal}
{\sl Let $V$ be a nondegenerate Hermitian, symplectic, or orthogonal space over $\FF_{q}$ that 
contains the orthogonal sum of $s \geq 2$ isometric, nondegenerate subspaces 
$\oplus^{s}_{i=1}V_{i}$. Assume $g \in I(V)$ stabilizes each $V_{i}$, and, after an isometric 
identification of $V_{i}$ with each other, induces the same action on all $V_{i}$. Assume in 
addition that $s \geq 3$ if $q \equiv 3 (\mod 4)$ and $V$ is not Hermitian. Then $g$ fixes  
a nonzero totally singular subspace of $V$.}
\end{lemma}

\begin{proof}
By the assumption, we can find bases $(e_{ij}, 1 \leq j \leq d)$ in $V_{i}$, $1 \leq i \leq s$, 
such that the Gram matrices and the actions of $g$ with respect to these $s$ bases are the same. 
It follows that $g$ preserves the subspace 
$W_{\al} := \langle \sum^{s}_{i=1}\al_{i}e_{ij} \mid 1 \leq j \leq d \rangle$ for any 
$\al := (\al_{1}, \ldots ,\al_{s}) \in \FF_{q}^{s}$. If $q$ is even, choose 
$\al = (1,1,0, \ldots ,0)$. If $q$ is odd and $V$ is Hermitian, choose $\al = (1,\gamma,0, \ldots ,0)$
where $\gamma \in \FF_{q}$ has order $2(q^{1/2}+1)$. If $q \equiv 1 (\mod 4)$ and $V$ is not 
Hermitian, choose $\al = (1,\gamma,0, \ldots ,0)$ where $\gamma \in \FF_{q}$ has order $4$. Assume 
$q \equiv 3 (\mod 4)$ and $V$ is not Hermitian. Then $s \geq 3$, and we can find 
$\beta, \gamma \in \FF_{q}$ such that $1 + \beta^{2} + \gamma^{2} = 0$. Choose 
$\al = (1,\beta,\gamma,0, \ldots ,0)$ in this case. These choices of $\al$ ensure that $W_{\al}$ 
is totally singular.   
\end{proof}

In what follows, we will need the following observation to handle reducible elements of classical 
groups. 

\begin{lemma}\label{decom}
{\sl Let $V$ be a nondegenerate Hermitian, symplectic, or orthogonal space over $\FF_{q}$. Assume 
$g \in G := I(V)$ is an element that does not fix any nonzero totally singular subspace of 
$V$. Then one of the following holds.

{\rm (i)} $V$ can be decomposed into an orthogonal sum $\oplus^{s}_{i=1}V_{i}$ of nondegenerate 
subspaces, and $g$ acts irreducibly on each of them.

{\rm (ii)} $q$ is even, $V$ is an orthogonal space of even dimension $2n$, 
$G = GO^{\eps}_{2n}(q)$, and 
$g \in H$ for some subgroup $H \simeq Sp_{2n-2}(q) \times 2$ of $G$.}
\end{lemma}

\begin{proof}
We induct on $\dim(V)$. If $g$ is irreducible on $V$ then we are done. Otherwise consider a
nonzero minimal $g$-invariant subspace $U$ of $V$, its orthogonal complement $U^{\perp}$ (relative 
to the Hermitian or the bilinear form $B$ on $V$), and let $W := U \cap U^{\perp}$. If $W = 0$,
then $V = U \oplus U^{\perp}$, both $U$ and $U^{\perp}$ are nondegenerate and $g$-invariant, and
so we are done by induction hypothesis. Assume $W \neq 0$. By minimality, $U = W$, that is, 
the form $B$ is zero on $U$. By assumption, $U$ cannot be totally singular. The only case when 
this can happen is where $q$ is even and $V$ is an orthogonal space of (even) dimension $2n$. Observe 
that the set of all singular vectors in $U$ is a $g$-invariant subspace of codimension $1$ in $U$. 
So by minimality $\dim(U) = 1$, and $U$ is generated by some nonsingular vector $u \in V$. Since 
$g$ fixes $U$ and fixes the quadratic form, we see that $g(u) = u$. Thus 
$g \in H:= Stab_{G}(u) \simeq Sp_{2n-2}(q) \times 2$.    
\end{proof}

To link $\deg(g)$ to $\deg(g^{p^{i}})$ for a $p$-element $g$, we will use the following statement.
 
\begin{lemma}\label{root1}
{\sl {\rm (i)} Let $a,b \in GL_{n}(\FF)$ be matrices such that $[a,b]$ is scalar and 
$[a^{p^{m}},b] = \Id \neq [a^{p^{m-1}},b]$, for a prime $p$ and an integer $m \geq 1$. 
Then $\Spec(b)$ consists of all $p^{m}$-roots of the elements of $\Spec(b^{p^{m}})$.

{\rm (ii)} Assume $\ell = p$ and $b \in GL_{n}(\FF)$ is a $p$-element. Then 
$p^{i} \cdot \deg(b^{p^{i}}) \geq \deg(b) \geq p^{i}(\deg(b^{p^{i}}) - 1) +1$ for every
integer $i \geq 1$.}
\end{lemma}

\begin{proof}  
(i) Since $[a,b]$ centralizes $a$ and $b$, $[a,b^{p^{m}}] = [a,b]^{p^{m}} = [a^{p^{m}},b] = \Id$. 
In particular, $[a,b] = \om \cdot \Id$ with $\om^{p^{m}} = 1$. Since 
$\Id \neq [a^{p^{m-1}},b] = \om^{p^{m-1}} \cdot \Id$, $\om$ is a primitive $p^{m}$-root of unity 
in $\FF$. Let $W = \FF^{n}$ be the natural module for $GL_{n}(\FF)$, $\lam \in \Spec(b^{p^{m}})$, 
and let $W_{\lam}$ be the $\lam$-eigenspace of $b^{p^{m}}$ in $W$. Obviously, $W_{\lam}$ is 
fixed by $a$ and $b$. Let $\beta$ be an eigenvalue of $b$ on $W_{\lam}$, and $v$ a 
corresponding eigenvector for $b$. Then $ba^{i}v = [b,a^{i}]a^{i}bv = \om^{-i}\beta a^{i}v$. 
Thus $a^{i}v$ belongs to the $\om^{-i}\beta$-eigenspace of $b$ in $W_{\lam}$, 
for $i=0, \ldots ,p^{m}-1$, whence $\beta, \om\beta, \ldots ,\om^{p^{m}-1}\beta$ are 
eigenvalues of $b$ on $W_{\lam}$.

(ii) Assume $\deg(b^{p^{i}}) = k+1$. Then $b$ is annihilated by 
$(t^{p^{i}}-1)^{k+1} = (t-1)^{p^{i}(k+1)}$ and so $\deg(b) \leq p^{i}(k+1)$. If 
$\deg(b) \leq p^{i}k$, then $b$ is annihilated by $(t-1)^{p^{i}k} = (t^{p^{i}}-1)^{k}$, and 
so $b^{p^{i}}$ is annihilated by $(t-1)^{k}$, a contradiction.
\end{proof}

\begin{corol}\label{two}
{\sl Assume that $S \lhd G/Z(G) \leq \Aut(S)$ for a simple non-abelian group $S$, and 
that $(G/Z(G))/S$ is elementary abelian of order $4$, resp. cyclic. Consider any irreducible 
$\FF G$-representation $\Theta$ of degree $> 2$, resp. $> 1$. 

{\rm (i)} If $h \in G \setminus Z(G)$, then $\Theta(h)$ is not scalar. 

{\rm (ii)} If $\ell = 2$ and $g \in G$ with $o(g) = 2^{a}$, then $d_{\Theta}(g) \geq 2^{a-1}+1$.}
\end{corol}

\begin{proof} 
(i) Assume $\Theta(h)$ is scalar. Let $Z := Z(G)$, 
$K := \{x \in G \mid \Theta(x) \mbox{ is scalar}\}$ and let $L := K/Z$. Then $1 \neq L \lhd G/Z$. 
Notice that $L \cap S \neq 1$. (Otherwise $[L,S] = 1$ and so $L \leq C_{G/Z}(S) = 1$, a 
contradiction.) By simplicity, $L \geq S$ and so $G/K \simeq (G/Z)/L$ is either cyclic or 
$\simeq \ZZ_{2}^{2}$. In the former case, the quotient of $\Theta(G)$ by its central subgroup 
$\Theta(K)$ is cyclic, whence $\Theta(G)$ is abelian and so $\dim(\Theta) = 1$ by 
irreducibility, a contradiction. In the latter case, the quotient of $\Theta(G)$ by its central 
subgroup $\Theta(K)$ is isomorphic to $\ZZ_{2}^{2}$, and so $\dim(\Theta) \leq 2$, 
a contradiction. 

(ii) Setting $h = g^{2^{a-1}}$, we see that $o(h) = 2$ and so $h \notin Z(G)$. By
(i), $d_{\Theta}(h) = 2$, whence we are done by Lemma \ref{root1}(ii).
\end{proof}

Notice that members of the family $\CLQ$ satisfy the assumption made on $G$ in Corollary \ref{two}. 

\begin{lemma}\label{prod1}
{\sl Let $g \in G$ be a $p$-element of order $p^{k}$, and let $U$ and $V$ be $\FF G$-modules.
Then $d_{U \otimes V}(g) \geq d_{U}(g)$. 

{\rm (i)} Assume that $\ell \neq p$ and $d_{U}(g) + d_{V}(g) > p^{k}$. Then 
$d_{U \otimes V}(g) = p^{k}$. 

{\rm (ii)} Assume that $\ell = p$ and $k = 1$. Then 
$d_{U \otimes V}(g) = \min\{p,d_{U}(g) + d_{V}(g) -1\}$. In general, if 
$V_{i}$ are $\FF G$-modules for $1 \leq i \leq s$, then 
$d_{V_{1} \otimes \ldots \otimes V_{s}}(g) = \min\{p,1-s+\sum^{s}_{i=1}d_{V_{i}}(g)\}$.

{\rm (iii)} Assume that $\ell = p \not{|} d_{U}(g)$ and $g|_{V} \neq 1_{V}$. Then 
$d_{U \otimes V}(g) > d_{U}(g)$.}
\end{lemma}

\begin{proof} 
The inequality $d_{U \otimes V}(g) \geq d_{U}(g)$ is obvious.
 
(i) Let $S$ denote the set of all $p^{k}$-roots of unity in $\FF$. Assume that 
$d_{U \otimes V}(g) < p^{k}$, say $c \in S$ is not an eigenvalue for $g$ on $U \otimes V$. Let 
$\{a_{1}, \ldots ,a_{m}\}$, resp. $S \setminus \{b_{1}, \ldots ,b_{n}\}$ be the set of all
distinct eigenvalues of $g$ on $U$, resp. on $V$, for some integers $m, n \geq 1$.  
For each $i$ with $1 \leq i \leq m$ we have $a_{i}S = S$, hence 
$c \in \{a_{i}b_{1}, \ldots ,a_{i}b_{n}\}$ and so $a_{i}b_{j} = c$ for exactly one index $j$. 
Thus there are $m$ ordered pairs $(i,j)$ such that $a_{i}b_{j} = c$. On the other hand,
for a fixed $j$, there is at most one index $i$ such that $a_{i}b_{j} = c$. It follows that 
there are at most $n$ ordered pairs $(i,j)$ such that $a_{i}b_{j} = c$, i.e.
$m \leq n$ and so $d_{U}(g) + d_{V}(g) \leq p^{k}$.

(ii) Let $J_{a}$ denote the Jordan block of size $a$ and with eigenvalue $1$ and let 
$1 \leq l \leq m \leq p$. According to \cite[Theorem VIII.2.7]{Fe}, the Jordan canonical form of 
$J_{l} \otimes J_{m}$ equals $\diag (J_{m+l-1},J_{m+l-3},\ldots ,J_{m-l+1})$ if $l + m \leq p$, and
$$\diag (\underbrace{J_{p},\ldots ,J_{p}}_{(m+l-p) \mbox{ \tiny{times}}},
  J_{2p-m-l-1},J_{2p-l-m-3}, \ldots ,J_{m-l+1})$$
if $l + m > p$. Now the first claim follows by taking $l = d_{U}(g)$ and $m = d_{V}(g)$.
The second claim follows from the first by induction on $s$.

(iii) Setting $h := g-1$, we see that 
$h(u \otimes v) = hu \otimes hv + u \otimes hv + hu \otimes v$ for $u \in U$ and $v \in V$.
We can choose $u$ and $v$ such that $h^{n}u = 0 \neq h^{n-1}u$, $h^{2}v = 0 \neq hv$, 
where $n := d_{U}(g)$. Induction on $m$ shows that 
$h^{m}(u \otimes v) = m(h^{m}u \otimes hv + h^{m-1}u \otimes hv) + h^{m}u \otimes v$. Thus 
$h^{n}(u \otimes v) = nh^{n-1}u \otimes hv \neq 0$.    
\end{proof}

\begin{lemma}\label{filtr}
{\sl Let $V$ be a finite dimensional $\FF_{p}$-space, $G \leq GL(V)$, and let $h \in G$ be a 
$p$-element. For any positive integer $k$, let $V_{k} := \Ker((h-1)^{k})$.

{\rm (i)} If $n \geq k$ then $V_{n+1}/V_{n}$ embeds in $V_{k+1}/V_{k}$ as a $C_{G}(h)$-module.

{\rm (ii)} Assume $h = g^{p^{b}}$ for some $p$-element $g \in G$ and $n := d_{V}(h)$.
If $d_{V_{n}/V_{n-1}}(g) = m$, then $\dvg = (n-1)p^{b} + m$.}
\end{lemma}

\begin{proof}
(i) Clearly, the map $f~:~V_{n+1}/V_{n} \to V_{k+1}/V_{k}$ defined by 
$f(v+V_{n}) = (h-1)^{n-k}(v)+V_{k}$ is an injective $C_{G}(h)$-map.

(ii) By the choice of $n$, $V_{n} = V$. For any $v \in V$, $v' := (g-1)^{m}v \in V_{n-1}$, whence 
$$(g-1)^{(n-1)p^{b}+m}v = ((g-1)^{p^{b}})^{n-1}v' = (h-1)^{n-1}v' = 0~.$$
On the other hand, by the choice of $m$, there is some $u \in V$ such that
$(g-1)^{m-1}u \notin V_{n-1}$ and so 
$0 \neq (h-1)^{n-1}(g-1)^{m-1}u = (g-1)^{(n-1)p^{b}+m-1}u$. Thus $\dvg = (n-1)p^{b} + m$.  
\end{proof}

Let $\CLL(G)$ denote the set of conjugacy classes of $\ell'$-elements of a finite group $G$.
The following statement is basically due to Brauer:

\begin{lemma}\label{brauer1}
{\sl Let $G$ be a finite group. Assume a finite group $A$ acts on $\CLL(G)$ and on $\IBRL(G)$
in such a way that, for any $a \in A$, any $\chi \in \IBRL(G)$, and for any $\ell'$-element
$g \in G$, $\chi(g) = \chi^{a}(g^{a})$, where $g^{a}$ is an element in $G$ such that 
$(g^{a})^{G} = (g^{G})^{a}$. Then the following statements hold.

{\rm (i)} For any $a \in A$, the number of $a$-fixed conjugacy classes of $\ell'$-elements in $G$ is 
equal to the number of $a$-fixed irreducible $\ell$-modular Brauer characters of $G$.

{\rm (ii)} The numbers of $A$-orbits on $\CLL(G)$ and on $\IBRL(G)$ are the same.}
\end{lemma}

\begin{proof}
(i) Let $\CLL(G) = \{g_{1}^{G}, \ldots ,g_{m}^{G}\}$ and 
$\IBRL(G) = \{\varphi_{1}, \ldots ,\varphi_{m}\}$. Then the matrix 
$X := \left(\varphi_{i}(g_{j})\right)_{1 \leq i,j \leq m}$ is nondegenerate by 
\cite[Theorem 60.4]{Do}. Now one can repeat the proof of \cite[Theorem (6.32)]{Is} verbatim.
(ii) is a consequence of (i).  
\end{proof}

\begin{corol}\label{brauer2}
{\sl Let $G$ be a normal subgroup of a finite group $H$. Define the action of $h \in H$ on 
$\CLL(G)$ and $\IBRL(G)$ as follows: $(g^{G})^{h} = (h^{-1}gh)^{G}$, 
$\varphi^{h}(g) = \varphi(hgh^{-1})$. Then the following statements hold. 

{\rm (i)} The numbers of $H$-orbits on $\CLL(G)$ and on $\IBRL(G)$ are the same.

{\rm (ii)} Assume $H/(C_{H}(G)G)$ is cyclic. Then the number of $H$-stable conjugacy classes of 
$\ell'$-elements in $G$ is equal to the number of $H$-stable irreducible $\ell$-modular Brauer 
characters of $G$.

{\rm (iii)} Assume $H/G$ is cyclic. Then the number of $H$-stable conjugacy classes of 
$\ell'$-elements in $G$ is equal to the number of irreducible $\ell$-modular Brauer characters of 
$G$ that are extendible to $H$.}
\end{corol}

\begin{proof}
(i) follows from Lemma \ref{brauer1}(ii).

(ii) Clearly, $C_{H}(G)G$ acts trivially on both $\CLL(G)$ and $\IBRL(G)$. Hence we can apply 
Lemma \ref{brauer1}(i) to $a$, a generator of $H/(C_{H}(G)G)$.

(iii) follows from (ii) by observing that $\varphi \in \IBRL(G)$ is $H$-stable precisely
when it is extendible to $H$, cf. \cite[Theorem III.2.14]{Fe}.   
\end{proof}

\section{Some groups of low rank}

\begin{examp}\label{su30}
{\em We describe some properties of irreducible Brauer characters of $G := SU_{3}(q)$ and 
$H := GU_{3}(q)$ with $q = p^{f}$. 

(i) {\it The number of $\varphi \in \IBRL(G)$ that does not extend to $H$ equals} 
$$\left\{ \begin{array}{ll}0, & \mbox{if }(3,q+1) = 1 \mbox{ or }\ell = p,\\
 9, & \mbox{if }3|(q+1) \mbox{ and }\ell \neq 3,p,\\ 
 3, & \mbox{if }\ell = 3|(q+1).\end{array}\right.$$
Indeed, the statement is obvious if $(3,q+1) = 1$, as in this case $H = Z(H) \times G$. 
Assume $3|(q+1)$. It is shown in \cite{Ge} that $G$ has exactly $9$ conjugacy classes that are not
$H$-stable, where $3$ of them consist of elements of order $P$ and $6$ consist of elements
of order $3P$, with $P = p$ if $p > 2$, and $P = 4$ if $p = 2$. If $\ell = p$ then none of them 
is an $\ell'$-class. If $\ell \neq 3,p$ then all of them are $\ell'$-classes. If $\ell = 3$ 
then exactly $3$ of them are $\ell'$-classes. So the statement follows from Corollary 
\ref{brauer2}(iii). 

\smallskip
(ii) Assume $\ell \neq p$ and $\varphi \in \IBRL(H)$. Here we show that {\it either $\varphi$ 
lifts to characteristic $0$, or there is a linear character $\lam$ such that one of the 
following holds:

{\rm (a)} $3 \neq \ell|(q^{2}-q+1)$, $\varphi\lam = \widehat{St}-1$, where $St$ is the Steinberg
character;

{\rm (b)} $q \equiv 1 (\mod 4)$, $\ell = 2$, $\varphi\lam = \hat{\chi}-1$ for some 
$\chi \in \Irr(H)$ of degree $q(q^{2}-q+1)$.\\}
To prove the claim, we identify $H$ with its dual group. We may assume $q > 2$ as 
$H$ is solvable and so $\varphi$ lifts if $q = 2$. Assume that $\varphi$ belongs to 
$\ECL(H,(s))$ for a semisimple $\ell'$-element $s \in H$, cf. \cite{BM}. There are $5$ 
possibilities for $C_{H}(s)$: $GU_{1}(q)^{3}$, $GL_{1}(q^{2}) \times GU_{1}(q)$, $GU_{1}(q^{3})$, 
$GU_{2}(q) \times GU_{1}(q)$, and $GU_{3}(q)$. Moreover, it is shown in the proof of 
\cite[Proposition 11.3]{GMST} that $\varphi$ lifts to characteristic $0$ in the first $3$ cases.

Consider the case $C_{H}(s) = GU_{2}(q) \times GU_{1}(q)$. According to Lusztig's classification 
\cite{L} of irreducible characters of $H$, $\EC(H,(s)) \cap \Irr(H) = \{\al,\beta\}$
with $\al(1) = q^{2}-q+1$ and $\beta = q(q^{2}-q+1)$. Moreover, by \cite{FS} and \cite{GH}, 
$\HAA$ and $\HAB$ form a basic set for $\ECL(H,(s)) \cap \IBRL(H)$. 
In particular, $\psi \in \ECL(H,(s)) \cap \IBRL(H)$ has degree divisible by $q^{2}-q+1$ and 
therefore $\HAA$ is irreducible, and $|\ECL(H,(s)) \cap \IBRL(H)| = 2$. It is easy to
see that $\al$ and $\beta$ are irreducible over $G$. If $3 \neq \ell|(q^{2}-q+1)$, then
$\al$ and $\beta$ have $\ell$-defect $0$. If $2 \neq \ell|(q-1)$, then $\HAA$ 
and $\HAB$ are irreducible over $G$ by \cite{Ge}. It remains to consider the case
$\ell|(q+1)$ and $\varphi \neq \HAA$. Clearly, $s$ is represented by $\diag(x,x,y)$ with 
$x^{q+1} = y^{q+1} = 1$ and $x \neq y$. Since $\ell|(q+1)$, $s$ is centralized by an element 
$t := \diag(z,1,1)$ with $|z| = \ell$. Observe that $C_{H}(st) = GU_{1}(q)^{3}$, so 
$\EC(H,(st)) \cap \Irr(H)$ contains a character $\gamma$ of degree $(q-1)(q^{2}-q+1)$; moreover,
$\gamma \in \ECL(H,(s))$. It follows that $\HAC = a\varphi + b\HAA$ with $0 \leq a,b \in \ZZ$. 
Since $(\gamma|_{U},1_{U})_{U} = 0$ and $(\al|_{U},1_{U})_{U} = 1$ for a Sylow $p$-subgroup $U$ of
$H$ (cf. \cite{Ge}), $b = 0$. Also, $\gamma(g) = -1$ for some element $g \in U$, so $a = 1$. 
Thus $\varphi = \HAC$.      

Now we suppose that $C_{H}(s) = GU_{3}(q)$, i.e. $s \in Z(H)$. Multiplying $\varphi$ by a 
linear character $\lam$ of $H$, we may assume that $s = 1$. In this case, 
$\EC(H,(s)) \cap \Irr(H) = \{1_{H},\zeta, St\}$ with $\zeta$ the cuspidal unipotent character
(of degree $q^{2}-q$) of $H$. It is well known that $\HAZ$ is irreducible over $G$. We may therefore
assume that $\varphi$ is a nontrivial constituent of $\widehat{St}$. If $2 \neq \ell|(q-1)$ then 
$\widehat{St}$ is irreducible over $G$, cf. \cite{Ge}. If $3 \neq \ell|(q^{2}-q+1)$ then we 
arrive at (a) by \cite{Ge}. If $2 \neq \ell|(q+1)$, then $\varphi$ lifts to a complex character
of degree $(q-1)(q^{2}-q+1)$ by the results of \cite{Ge} and \cite{OW2}. Assume $\ell = 2$ and
$q$ is odd. From the results of \cite{H2}, it follows that $\varphi$ again lifts to a complex 
character of degree $(q-1)(q^{2}-q+1)$ if $q \equiv 3 (\mod 4)$, whereas one arrives at (b) 
if $q \equiv 1 (\mod 4)$.    

\smallskip
(iii) We will assume $3|(q+1)$ and explicitly determine the Brauer characters of $G$ that do not 
extend to $H$. If $\ell = 0$, then they are the ones labeled as $\chi^{(u)}_{(q-1)(q^{2}-q+1)/3}$ and 
$\chi^{(u,v)}_{(q+1)^{2}(q-1)/3}$ with $0 \leq u,v \leq 2$ in \cite{Ge}. Observe that 
$(q-1)(q^{2}-q+1)/3 < (q+1)^{2}(q-1)/3$, $\chi^{(u)}_{(q-1)(q^{2}-q+1)/3}$ with different $u$ 
take distinct values at some $p$-element, and $\chi^{(u,v)}_{(q+1)^{2}(q-1)/3}$ with different 
$(u,v)$ take distinct values at some $\{3,p\}$-element. If $2 \neq \ell|(q-1)$,
then all of these $9$ complex characters are of $\ell$-defect $0$, so their reductions modulo 
$\ell$ yield the $9$ Brauer characters that do not extend to $H$. These $9$ characters also stay 
irreducible modulo $\ell$ when $2,3 \neq \ell |(q^{3}+1)$, cf. \cite{Ge}. If $\ell = 3$, then
$\chi^{(u)}_{(q-1)(q^{2}-q+1)/3}$ are irreducible modulo $3$ by \cite{OW2}. Finally, assume 
$\ell = 2|(q-1)$. Clearly, $\chi^{(u,v)}_{(q+1)^{2}(q-1)/3}$ are of $2$-defect $0$. On the 
other hand, $\chi^{(u)}_{(q-1)(q^{2}-q+1)/3}$ is labeled by the $PGU_{3}(q)$-conjugacy class of 
$s := \diag(1,\om,\om^{2})$ with $|\om| = 3$, whose centralizer in $PGU_{3}(q)$ is 
$GU_{1}(q)^{2} : \ZZ_{3}$. It follows by \cite{HM} that $\chi^{(u)}_{(q-1)(q^{2}-q+1)/3}$ are 
irreducible modulo $2$. Thus, {\it any $\varphi \in \IBRL(G)$ either lifts to characteristic $0$ 
or extends to $H$. Furthermore, if $q \geq 3$ and $\ell \neq p$ then the smallest degree of 
$\ell$-modular Brauer characters of $G$ which are neither trivial nor a Weil character is 
$(q-1)(q^{2}-q+1)/(3,q+1)$.}}      
\end{examp}

Arguing similarly, we obtain: 

\begin{examp}\label{sl30}
{\em Let $G := SL_{3}(q)$ and $H := GL_{3}(q)$ with $q = p^{f}$. 

(i) {\it The number of $\varphi \in \IBRL(G)$ that does not extend to $H$ equals} 
$$\left\{ \begin{array}{ll}0, & \mbox{if }(3,q-1) = 1 \mbox{ or }\ell = p,\\
 9, & \mbox{if }3|(q-1) \mbox{ and }\ell \neq 3,p,\\ 
 3, & \mbox{if }\ell = 3|(q-1).\end{array}\right.$$

(ii) Assume $\ell \neq p$ and $\varphi \in \IBRL(H)$. Then {\it either $\varphi$ 
lifts to characteristic $0$, or $\ell|(q^{2}+q+1)$ and there is a linear character $\lam$ such 
that $\varphi\lam = \widehat{\tau}-1$, where $\tau$ is the unipotent character of degree 
$q^{2}+q$ of $H$.}

(iii) Assume $\ell \neq p$ and $\varphi \in \IBRL(G)$. Then {\it $\varphi$ either lifts to 
characteristic $0$ or extends to $H$.}}
\end{examp}

\begin{lemma}\label{sl21}
{\sl Let $g \in G = SL_{2}(q)$ with $q>3$, $o(g) = p^{a}$ and $g^{q+1} = 1$. 
Assume $\Theta \in \IBRL(G)$, $(\ell,q) = 1$, and $1 < \deg(\Theta(g)) < o(g)$. 
Then one of the following holds.

{\rm (i)} $q$ is odd, $p\ell \neq 4$, $o(g)=(q+1)/2$ and 
$\dim(\Theta) = \deg(\Theta(g)) = (q-1)/2$. Moreover, if $\ell \neq p$ then
$(-1)^{p-1}\not \in \Spec(\Theta(g))$.

{\rm (ii)} $q \equiv 3 (\mod 4)$, $p = \ell = 2$, $o(g) = (q+1)/2 = 2^{a}$, and 
$\deg(\Theta(g)) = o(g)-1$. Furthermore, $\dim(\Theta) = (q-1)/2$ or $\dim(\Theta) = q-1$.

{\rm (iii)} $q$ is even, $o(g) = |g| = q+1$. If $\ell = p$ then   
$\dim(\Theta) = \deg(\Theta(g)) = q-1$. If $\ell \neq p$, then either $\dim(\Theta) = q$ and 
$1 \not\in \Spec(\Theta(g))$, or $\dim(\Theta) = q-1$ and there is a primitive $(q+1)$-root 
$\eps \neq 1$ of unity in $\FF$ such that 
$\Spec(\Theta(g)) = \{1, \eps,\eps^{2}, \ldots \eps^{q}\} \setminus \{\eps,\eps^{-1}\}$.}
\end{lemma}

\begin{proof}
First assume that $\ell \neq p$. It is well known that $\Theta$ lifts to characteristic $0$, 
so we may assume $\ell = 0$. Now the statement follows from \cite{Z3} if $p > 2$, and by
inspecting the character table of $G$ if $p = 2$. If $\ell = p > 2$ then again the statement 
follows from \cite{Z3}. Finally, let $\ell = p = 2$. Since $2 < o(g) = 2^{a} | (q+1)$,
$q \equiv 3 (\mod 4)$. By \cite{Bu1}, $\dim(\Theta) = q+1$, $(q-1)/2$, or $q-1$. In the first
case, $\Theta$ (considered as a representation of $PSL_{2}(q)$) is of $2$-defect $0$
and so $\deg(\Theta(g)) = o(g)$. In the second case, assertion (ii) follows from \cite[\S2]{GT4}. 
\end{proof}

\begin{lemma}\label{sl22}
{\sl Let $g \in G = GL_{2}(q)$ with $3 < q \equiv 3 (\mod 4)$, $g \in G \setminus Z(G)$ a 
$2$-element. Assume $\Theta \in \IBRL(G)$, $(\ell,q) = 1$, and 
$1 < \deg(\Theta(g)) < o(g)$. Then $q+1 = 2^{c} \geq 8$, and one of the following holds.

{\rm (i)} $\dim(\Theta) = q$, $o(g) = 2^{c}$, and $\deg(\Theta(g)) = o(g)-1$. 

{\rm (ii)} $\dim(\Theta) = q-1$, $o(g) = 2^{c}$, and $\deg(\Theta(g)) = o(g)-2$. 

{\rm (iii)} $\dim(\Theta) = q-1$, $o(g) = 2^{c-1}$, $\deg(\Theta(g)) = o(g)-1$. Furthermore, 
$\Theta$ is reducible over $SL_{2}(q)$ if $\ell \neq 2$.}
\end{lemma}

\begin{proof}
It is well known that $\Theta$ lifts to characteristic $0$. 
Since $o(g) > 2$, $|g|$ is divisible by $4$ and so $g$ is irreducible on $\FF_{q}^{2}$. Thus we 
may embed $g$ in a maximal torus $\ZZ_{q^{2}-1}$ of $G$. Now the statements follow by direct 
computation with the character table of $G$ if $\ell \neq 2$. Assume that $\ell = 2$. Then 
$\Theta$ can be viewed as a representation of $O_{2'}(G) \times PGL_{2}(q)$ and so the lemma follows 
from \cite{GT4} if $\dim(\Theta) = q-1$. It remains to consider the case $\dim(\Theta) = q+1$. In this
case $\Theta$ is induced from a Borel subgroup $B$ of $G$, and for any $x \in G$, 
$xBx^{-1} \cap \la g \ra \leq Z(G)$, whence $\deg(\Theta(g)) = o(g)$.
\end{proof}

Observe that any irreducible $p$-element of $GL_{2}(q)$ with $p$ odd is contained in 
$SL_{2}(q)$ and has order dividing $q+1$. On the other hand, any irreducible $2$-element $g$ of 
$GL_{2}(q)$ with $q \equiv 1 (\mod 4)$ has $o(g) = 2$. So Lemmas \ref{sl21} and \ref{sl22} 
have determined $d_{\Theta}(g)$ for all irreducible $p$-elements in $GL_{2}(q)$ and 
for all cross characteristic representations $\Theta$. Next we turn to $GU_{2}(q)$.

\begin{lemma}\label{su2}    
{\sl Let $g \in G := GU_{2}(q)$ be a semisimple $p$-element with $o(g) = p^{a}$ and let  
$\Phi \in \IBRL(G)$ with $(\ell,q) = 1$. Suppose that $1 < \deg(\Phi(g)) < o(g)$. 
Then $g^{q+1}=1$ and one of the following holds.

{\rm (i)} $\ell \neq p$, $q+1 = p^{a}$, and $\deg(\Phi(g)) = \dim(\Phi) = q$.

{\rm (ii)} $q+1 = p^{a}$, $\dim(\Phi) = q-1$, and $\deg(\Phi(g)) = q-1$.

{\rm (iii)} $q+1 = 2p^{a}$, $\dim(\Phi) = q-1$, and $\deg(\Phi(g)) = (q-1)/2$.}
\end{lemma}

\begin{proof}
It is easy to see that $g$ is always reducible on the natural module $V := \FF_{q^{2}}^{2}$ for $G$. 
If $g$ fixes a singular $1$-space of $V$, then $\deg(\Phi(g)) = o(g)$ by \cite{DZ1}
(one can also see it directly by looking at the action of $g$ on a unipotent subgroup of 
order $q$ of $G$ and using Lemma \ref{perm}.). We will assume that $g$ fixes no singular 
$1$-space of $V$. It follows that $g$ fixes a decomposition of $V$ into an orthogonal sum of
two nondegenerate $1$-spaces, whence $g^{q+1} = 1$ and $g$ belongs to the class $C_{3}(k,l)$ 
for some $k,l$ in the notation of \cite{Enn}. Suppose that $p = \ell > 2$. Then 
$g \in Z(G)S$ with $S := SU_{2}(q)$. But the $p$-Sylow subgroups of $S$ are cyclic, so we are 
done by \cite{Z3}. Thus we will assume that $p \neq \ell$ if $\ell \neq 2$. It is well known
that $\Phi$ can be lifted to characteristic $0$; in particular, $\dim(\Phi) = q,q \pm 1$.

If $\dim(\Phi) = q+1$ then $\Phi$ is induced from a one-dimensional representation of 
a Borel subgroup $B$ of $G$. It is easy to check that all elements $b \in B$ with $b^{q+1}=1$ 
are in $Z(G)$. Hence $\deg(\Phi(g)) = o(g)$ by Lemma \ref{perm}.

Assume $\dim(\Phi) = q$. Twisting $\Phi$ with a linear character of $G$, we may assume that 
$\Phi$ is obtained by reducing the Steinberg character of $G$ modulo $\ell$, and 
$(\ell,q+1) = 1$; in particular, $\ell \neq 2$. Inspecting the character table of $G$ \cite{Enn}, 
we arrive at (i) if $\ell \neq p$.

Assume $\dim(\Phi) = q-1$. If $\ell \neq p$, then by inspecting the character table of $G$ 
\cite{Enn}, we arrive at (ii) and (iii). It remains to consider the case $p = \ell = 2$. 
If $q \equiv 1 (\mod 4)$ then $o(g) = 2$ (as $g^{q+1} = 1$) and so $\deg(\Phi(g)) = 2$. So we may 
assume $q \equiv 3 (\mod 4)$. Clearly, $\Phi(G)$ can be considered as a representation of 
$O_{2'}(Z(G)) \times PGL_{2}(q)$, whence we arrive at (ii) and (iii) by \cite{GT4}.
\end{proof}

To efficiently restrict our problem to certain natural subgroups, we need the following statement.
  
\begin{lemma}\label{prod2}
{\sl Let $G$ be a finite Lie-type group of simply connected type defined over $\FF_{q}$, with
$q = p^{f}$. Assume $G$ contains a central product $Y := X_{1}*X_{2}$, where $X_{1}$ contains a 
long-root subgroup $U$, and $X_{2}$ contains a quasisimple subgroup $T$. 
Assume $\Phi \in \IBRL(G)$ with $\ell \neq p$ and $\dim(\Phi) > 1$. Assume in addition that 
$G \not\in \{Sp_{2n}(q),F_{4}(q),\ta F_{4}(q),\ta B_{2}(q)\}$ if $2|q$, and 
$G \not\in \{G_{2}(q),\ta G_{2}(q)\}$ if $3|q$. Then $\Phi|_{Y}$ contains an irreducible 
constituent $\Phi_{1} \otimes \Phi_{2}$ such that $\Phi_{i} \in \IBRL(X_{i})$ for $i = 1,2$, 
$\Phi_{1}|_{U}$ is nontrivial, and $\dim(\Phi_{2}) > 1$.}
\end{lemma}

\begin{proof}
Assume the contrary. Decomposing the representation space $V$ of $\Phi$ into 
$V = C_{V}(U) \oplus [U,V]$, we see that $[U,V] \neq 0$ as $\dim(V) > 1$. The conditions imposed 
on $G$ imply that $O_{p}(C_{G}(U)) = Q$ is a $p$-subgroup of symplectic type and 
$Z(Q) = U$ (i.e. $Z(Q) = [Q,Q] = \Phi(Q)$ and $[x,Q] = Z(Q)$ for all $x \in Q \setminus Z(Q)$),
cf. \cite{LS}, \cite{MT}. By our assumption, all 
$T$-composition factors of $[U,V]$ are of dimension $1$, whence the (perfect) subgroup $T$ acts 
trivially on $[U,V]$. It follows that $[Q,T] \leq C_{Q}([U,V])$. Since $C_{Q}([U,V])$ is normal
in $Q$ and does not contain $Z(Q) = U$, it is contained in $Z(Q)$. In particular, 
$[Q,T] \leq Z(Q) = U$. But $[U,T] = 1$ by our assumptions, so $[[Q,T],T] = [[T,Q],T] = 1$,
whence $[T,Q] = [[T,T],Q] = 1$ by the three-subgroup lemma. On the other hand,
$C_{G}(Q) \leq Z(Q)Z(G)$ is solvable, and we arrive at a contradiction.
\end{proof}

The next lemma considers the symplectic groups in characteristic $2$ separately. 
In what follows, if $G \leq I(V)$ is a 
classical group with natural module $V$, then by a {\it standard subgroup} $Y$ of $G$ we 
mean any subgroup $Y = X_{1} * X_{2}$, with $X_{i} = I(V_{i}) \cap G$, and $V = V_{1} \oplus V_{2}$ 
is an orthogonal sum of (nondegenerate) subspaces.   
 
\begin{lemma}\label{prod21}
{\sl Let $G = Sp_{2n}(q)$ with $q = 2^{f}$, $n \geq 2$, and $(n,q) \neq (2,2)$. 
Let $Y := X_{1} \times X_{2} = Sp_{2}(q) \times Sp_{2n-2}(q)$ be a standard subgroup of $G$. 
Assume $\Phi \in \IBRL(G)$ with $\ell \neq p$ and $\dim(\Phi) > 1$. Then $\Phi|_{Y}$ contains an 
irreducible constituent $\Phi_{1} \otimes \Phi_{2}$ such that $\Phi_{i} \in \IBRL(X_{i})$ for $i = 1,2$, 
$\dim(\Phi_{2}) > 1$ and $\Phi_{1}$ is nontrivial on a Sylow $p$-subgroup of $X_{1}$.}
\end{lemma}

\begin{proof}
Assume the contrary and consider a Sylow $2$-subgroup $U$ of $X_{1}$. Then $|U| = q$. 
Decomposing the representation space $W$ of $\Phi$ into $W = C_{W}(U) \oplus [U,W]$, we see that 
$[U,W] \neq 0$ as $\dim(W) > 1$. By our assumption, all $X_{2}$-composition factors of $[U,W]$ are 
of dimension $1$, whence the perfect subgroup $[X_{2},X_{2}]$ acts trivially on $[U,W]$.  
Observe that $N_{G}(U) = Q : (X_{2} \times \ZZ_{q-1})$, where $Q$ is 
elementary abelian of order $q^{2n-1}$. Observe that the set of linear characters of $Q$ that are 
nontrivial on $U$ partition into two $N_{G}(U)$-orbits $\OC^{\eps}$ of length 
$q^{n-1}(q^{n-1}+\eps)(q-1)/2$ for $\eps = \pm$. Obviously, the $Q$-module $[U,W]$ has to afford at 
least one of these two orbits $\OC^{\eps}$. Thus $(X_{2}/[X_{2},X_{2}])\times \ZZ_{q-1}$ must act 
transitively on that orbit, a contradiction as $2(q-1) < |\OC^{\eps}|$.
\end{proof} 

\begin{corol}\label{prod4}
{\sl Let $G$ be a classical group and $Y := X_{1} * X_{2}$ be a standard subgroup of $G$.
Assume $\Phi \in \IBRL(G)$ with $(\ell,q) = 1$ and $\dim(\Phi) > 1$. Then $\Phi|_{Y}$ contains an 
irreducible constituent $\Phi_{1} \otimes \Phi_{2}$ such that $\Phi_{i} \in \IBRL(X_{i})$ and  
$\dim(\Phi_{i}) > 1$ for $i = 1,2$, provided one of the following holds.  

{\rm (i)} $G := Sp_{2n}(q)$, $Y = Sp_{2k}(q) \times Sp_{2l}(q)$, $1 \leq k \leq l$, and 
$(k,q) \neq (1,2)$, $(1,3)$. 

{\rm (ii)} $G := SL_{n}(q)$, $Y = SL_{k}(q) \times SL_{l}(q)$, $2 \leq k \leq l$, and 
$(k,q) \neq (2,2)$, $(2,3)$.

{\rm (iii)} $G := SU_{n}(q)$, $Y = SU_{k}(q) \times SU_{l}(q)$, $2 \leq k \leq l$, and 
$(k,q) \neq (2,2)$, $(2,3)$, $(3,2)$.

{\rm (iv)} $G := Spin^{\eps}_{n}(q)$, $Y = Spin^{-}_{2k}(q) * Spin^{-\eps}_{l}(q)$, 
and $6 \leq 2k \leq l$.}
\end{corol}

\begin{proof}
Our assumptions imply that $[X_{i},X_{i}]$ is perfect and $X_{1}$ contains a long-root subgroup, 
whence we are done in the cases of (ii) -- (iv) by Lemma \ref{prod2}. Assume the claim is false in 
the case of (i). Then every 
irreducible constituent of $\Phi|_{Y}$ is either trivial on $[X_{1},X_{1}]$ or trivial on 
$[X_{2},X_{2}]$. In particular, if $A \leq [X_{1},X_{1}]$ and $B \leq [X_{2},X_{2}]$, then every 
irreducible constituent of $\Phi|_{A * B}$ is either trivial on $A$ or trivial on $B$.
Now we get a contradiction by choosing $A \simeq Sp_{2}(q)$ to be a direct factor of a standard subgroup
of $X_{1}$ and $B = X_{2}$ and applying Lemmas \ref{prod2} and \ref{prod21} to $A \times B$ inside 
$Sp_{2(l+1)}(q)$ if $(k,q) \neq (2,2)$. If $(k,q) = (2,2)$, then we can choose $A \times B$ to be a 
standard subgroup $Sp_{2}(4) \times Sp_{2}(4)$ inside a subgroup $Sp_{4}(4)$ of $G$.
\end{proof}

\begin{lemma}\label{prod5}
{\sl Let $G = \Omega(V) = \Omega^{\pm}_{2n}(q)$ with $q = 2^{f} \geq 4$ and $n \geq 4$. Decompose the 
natural module $\FF_{q}^{2n}$ of $G$ into the orthogonal sum $V_{1} \oplus V_{2}$ of nondegenerate 
subspaces $V_{1}$ of dimension $4$ and type $-$, $V_{2}$ of dimension $2n-4$, and let 
$X_{i} := \Omega(V_{i})$ for $i = 1,2$. Assume $\Phi \in \IBRL(G)$ with 
$\ell \neq 2$ and $\dim(\Phi) > 1$. Then $\Phi|_{X_{1} \times X_{2}}$ contains an irreducible 
constituent $\Phi_{1} \otimes \Phi_{2}$ such that $\Phi_{i} \in \IBRL(X_{i})$ and 
$\dim(\Phi_{i}) > 1$ for $i = 1,2$.}
\end{lemma}

\begin{proof}
It suffices to prove the lemma for $n = 4$. Assume the statement is false 
for $n = 4$. We can find $\eps = \pm$ such that $(\ell,q-\eps) = 1$. 
Decompose $V_{2}$ into an orthogonal sum of nondegenerate subspaces: $2$-dimensional 
$V_{21} = \la a,b \ra$ of type $\eps$, and $2$-dimensional $V_{22} = \la c,d \ra$, with $d$ being  
nonsingular. Then $Stab_{G}(d)$ has the commutator subgroup $Y \simeq Sp_{6}(q)$. Consider the 
commutator subgroups $Y_{1}$ of  $Stab_{Y}(a,b)$ (so $Y_{1} \simeq Sp_{4}(q)$), and $Y_{2}$ of 
$\{ y \in Y \mid y_{V_{1}} = 1_{V_{1}} \}$ (so $Y_{2} \simeq Sp_{2}(q)$). Clearly, 
$\Phi|_{Y}$ has an irreducible constituent $\Phi'$ of dimension $> 1$. By Lemma \ref{prod21}, 
$\Phi'|_{Y_{1} \times Y_{2}}$ has an irreducible constituent $\Phi^{1} \otimes \Phi^{2}$, with
$\Phi^{i} \in \IBRL(Y_{i})$ and $\dim(\Phi^{i}) > 1$ for $i = 1,2$. Since $Y_{1} > X_{1}$, 
$\Phi^{1}|_{X_{1}}$ has an irreducible constituent $\Phi^{11}$ of dimension $> 1$. Since 
$Y_{2} > T := \Omega(V_{21}) \simeq \ZZ_{q-\eps}$ and $(\ell,q-\eps) = 1$, $\Phi^{2}|_{T}$ contains 
a nontrivial irreducible constituent $\lam$. Thus $\Phi|_{X_{1} \times T}$ contains the 
irreducible constituent $\Phi^{11} \otimes \lam$ with $\dim(\Phi^{11}) > 1$ and $\lam \neq 1_{T}$. 
On the other hand, we have assumed that every irreducible constituent of $\Phi|_{X_{1} \times X_{2}}$
is either trivial on $X_{1}$ or trivial on $[X_{2},X_{2}]$. This is a contradiction, since 
$[X_{2},X_{2}] > T$.        
\end{proof}

Next we prove the counterpart of Lemma \ref{prod5} for odd $q$.

\begin{lemma}\label{prod6}
{\sl Let $V = \FF_{q}^{n}$ be a nondegenerate orthogonal space with $n \geq 8$ and $q \geq 5$ odd. 
Consider a standard subgroup $X_{1} * X_{2} = Spin^{-}_{4}(q) * Spin_{n-4}(q)$ of 
$G := Spin(V)$. Assume $\Phi \in \IBRL(G)$ with $(\ell,q) = 1$ and $\dim(\Phi) > 1$. Then 
$\Phi|_{X_{1} * X_{2}}$ contains an irreducible constituent $\Phi_{1} \otimes \Phi_{2}$ such 
that $\Phi_{i} \in \IBRL(X_{i})$ and $\dim(\Phi_{i}) > 1$ for $i = 1,2$.}
\end{lemma}

\begin{proof}
It suffices to prove the statement for $n = 8$. Let $Q$ denote the quadratic form on $V$. Consider 
an orthogonal decomposition 
$V = \la e_{1},f_{1} \ra \oplus \la e_{2},f_{2} \ra \oplus \la e_{3},f_{3} \ra \oplus 
   \la e_{4},f_{4} \ra$, such that $Q(e_{1}) = Q(f_{1}) = Q(e_{3}) = Q(f_{3}) = 0$,
$Q(e_{2}) = 1$, $e_{2} \perp f_{2}$, $e_{4} \perp f_{4}$, $Q(e_{4}) = -1$, and $\la e_{4},f_{4} \ra$ 
is of type $-$. Observe that $A_{1} := \la e_{1},f_{1},f_{2} \ra$ is nondegenerate, 
$A_{2} := \la e_{3},f_{3},e_{2},e_{4} \ra$ is nondegenerate of type $+$, and $A_{1} \perp A_{2}$.
Inside $H := Spin(A_{1} \oplus A_{2}) = Spin_{7}(q)$, $Spin(A_{2})$ contains a long-root subgroup 
$U$ of $H$, and $Spin(A_{1}) = Spin_{3}(q) \simeq SL_{2}(q)$ contains a quasisimple subgroup 
$T \simeq SL_{2}(q)$ if $q \neq 9$, and $T \simeq SL_{2}(5)$ if $q = 9$. In particular, 
$\Mult(T) = 1$. By Lemma \ref{prod2} applied to $H$, $\Phi$ contains an irreducible constituent
$\Phi_{11} \otimes \Phi_{12}$ with $\Phi_{1i} \in \IBRL(Spin(A_{i}))$ for $i = 1,2$, 
$\dim(\Phi_{11}) > 1$, and $\Phi_{12}|_{U}$ is nontrivial. Setting
$A_{3} := \la e_{3},f_{3},e_{2},e_{4},f_{4} \ra$, we see that $Spin(A_{3}) = Spin_{5}(q)$ is a
perfect subgroup containing $U$. It follows that $\Phi$ contains an irreducible constituent
$\Phi_{21} \otimes \Phi_{23}$ such that $\Phi_{2j} \in \IBRL(Spin(A_{j}))$ and $\dim(\Phi_{2j}) > 1$
for $j = 1,3$. Next, $A_{4} := \la e_{3},f_{3},e_{4},f_{4} \ra$ is a nondegenerate subspace of type 
$-$ of $A_{3}$, and $X_{1} = Spin(A_{4})$ is perfect. Hence $\Phi$ contains an irreducible constituent
$\Phi_{31} \otimes \Phi_{34}$ such that $\Phi_{3j} \in \IBRL(Spin(A_{j}))$ and $\dim(\Phi_{3j}) > 1$
for $j = 1,4$. Finally, $A_{5} := \la e_{1},f_{1},e_{2},f_{2} \ra$ is an $(n-4)$-dimensional 
nondegenerate subspace containing $A_{1}$, and $X_{2} = Spin(A_{5})$ contains the perfect subgroup 
$Spin(A_{1})$. Consequently, $\Phi$ contains an irreducible constituent
$\Phi_{44} \otimes \Phi_{45}$ such that $\Phi_{4i} \in \IBRL(Spin(A_{i}))$ and 
$\dim(\Phi_{4i}) > 1$ for $i = 4,5$, as stated.     
\end{proof}

Note that the subgroup $X_{1}$ in Lemmas \ref{prod5}, \ref{prod6} does not contain a long-root 
subgroup of $G$. 

\section{Regular semisimple elements}

In this section we study the semisimple elements which are minimal
in a certain sense. For instance, we consider $G = GL^{\eps}_{r}(q)$, with 
$\eps = +$ for $GL_{r}(q)$ and $\eps =  -$
for $GU_{r}(q)$ and $r$ is a prime (and $r \geq 3$ for $\eps =
-$).  Let $p \neq \ell$ be a {\it primitive prime
divisor} of $q^{r} - 1$, resp. $q^{2r}-1$, for $\eps = +$, resp.
$\eps = -$ (cf. \cite{Zs} for the definition of such prime
divisors), and let $g \neq 1$ be a $p$-element in $G$. It is easy to see that $C_{G}(g)$ is a 
maximal torus of $G$ of order $q^{r} -\eps$ (so $g$ is regular), 
and $g$ is not contained in any proper parabolic subgroup of $G$. In fact, the elements
we consider are the $p$-elements for which the Sylow $p$-subgroups $S$ of $G$ are cyclic
and $C_{G}(S)$ is abelian -- we will not use this observation in the sequel, though.

For the reader's convenience, we record the following statement, (see also 
\cite[Lemma 10.2]{GT2}: 

\begin{lemma}\label{senkr}
{\sl Let $\GC$ be a connected reductive algebraic group, $F$ a Frobenius map on $\GC$,
and let the pair $(\GCD,\FD)$ be dual to $(\GC,F)$. Assume $\chi \in \Irr(\GCF) \cap \ECS$
for some semisimple element $s \in \GCDF$. Then $\chi$ is a $\QQ$-linear combination of those 
$\RGT$ belonging to $\ECS$ and some class function that vanishes at semisimple elements of $\GCF$.
The same statement holds if one replaces the Lusztig series $\ECS$ by the rational series 
$\EC(\GCF,(s)_{\GCDF})$.}
\end{lemma}

\begin{proof}
Let $I := \Irr(\GCF) \cap \ECS$, $I' := \Irr(\GCF) \setminus I$, 
$J := \{\RGT \mid \RGT \in \ECS\}$, and $J' := \{\RGT \mid \RGT \notin \ECS\}$. Notice that
the scalar product of characters is positive definite on the space of rational-valued 
class functions on $\GCF$. Clearly,
$I \perp I'$, $J \subseteq \langle I \rangle_{\QQ}$, and 
$J' \subseteq \langle I' \rangle_{\QQ}$ (as Lusztig series form a partition of $\Irr(\GCF)$).
Write $\langle I \rangle_{\QQ}$ as the orthogonal sum $\langle J \rangle_{\QQ} \oplus S$.
Then any function $f \in S$ is orthogonal to $J$ and also to $J'$, so $f$ is an orthogonal
function; in particular it vanishes at semisimple elements of $\GCF$ by 
\cite[Lemma 10.2(i)]{GT2}. Now any $\chi \in I$ can be written as $\al + \beta$ where 
$\al \in \langle J \rangle_{\QQ}$ and $\beta \in S$, and so we are done. The same argument applies 
to the rational series $\EC(\GCF,(s)_{\GCDF})$ in place of $\ECS$, as the rational series also 
form a partition of $\Irr(\GCF)$ (cf. \cite[Prop. 14.41]{DM}).
\end{proof}

One of the main results of this section is the following theorem:

\begin{theor}\label{min2}
{\sl Let $\FF$ be an algebraically closed field of characteristic $\ell$, $\GC$ a connected 
reductive algebraic group in characteristic $r \neq \ell$,
$F$ a Frobenius map on $\GC$, and let $G := \GCF$. Let $\GCD$ be an algebraic group with
a Frobenius map $\FD$ such that $(\GCD,\FD)$ is dual to $(\GC,F)$. Let $p \neq \ell$ be
a prime with the property that any nontrivial $p$-element in $\GCDF$ is regular
semisimple in $\GCD$ and that $(p,|Z(\GC)/Z(\GC)^{\circ}|) = 1$.
Consider any irreducible $\FF G$-module $V$, with Brauer character $\psi$
and any nontrivial $p$-element $g \in G$.

\smallskip
{\rm (A)} Then one of the following holds.

{\rm (i)} $\psi(g) \in \ZZ$ and $\dvg \geq \varphi(|g|)$.

{\rm (ii)} $V$ lifts to characteristic $0$.

\smallskip
{\rm (B)} Assume in addition that $|T|^{3} \leq |G|_{r'}$ for any maximal torus $T$ of $G$
with order divisible by $p$, that any nontrivial $p$-element in $\GCF$ is regular
semisimple in $\GC$, and that $(p,|Z(\GCD)/Z(\GCD)^{\circ}|) = 1$. Then $\dvg = |g|$ in case 
of {\rm (ii)}.}
\end{theor}

\begin{proof}
Notice that in statement (i), the second conclusion follows from the first one (and the
assumption that any nontrivial power of $g$ acts nontrivially on $V$).

1) First we consider a generalized Deligne-Lusztig character $\RGT$, where $\TC$ is a
maximal torus of $\GC$ and $\theta$ a linear character of $\TCF$. By \cite[Cor. 12.18]{DM},
$\RGT \cdot St$ and $\Ind^{\GCF}_{\TCF}(\theta)$ are equal up to sign, where $St$ denotes the 
Steinberg character of $G$. Since $St(g)$ is a nonzero integer, it follows that $\RGT(g) = 0$ 
if no $G$-conjugate of $g$ is contained in $\TC$. Next assume that
$g \in \TC$ but the (multiplicative) order $M$ of $\theta$ is not divisible by $p$.
Consider any $G$-conjugate $xgx^{-1}$ of $g$ that lies in $\TCF$. Then
$(\theta(xgx^{-1}))^{k} = 1$ for $k = M$ and for $k = |g|$. Since $(M,|g|) = 1$, we conclude 
$\theta(xgx^{-1}) = 1$. Thus $\Ind^{\GCF}_{\TCF}(\theta)$, and so $\RGT(g)$, is rational.

2) By the results of Brou\'e and Michel \cite{BM}, $V$
belongs to a union $\ECL(\GCF,(s))$ for some semisimple $\ell'$-element $s \in \GCDF$.
It is well-known, cf. \cite[Theorem 61.6]{Do} for instance, that any Brauer character in an
$\ell$-block is an integral combination of the restrictions to $\ell'$-classes of $G$ of 
the complex irreducible characters belonging to the block. It now follows by Lemma \ref{senkr}
that $\psi(g)$ is a $\QQ$-linear combination of $\RGT(g)$ with $\RGT$ belonging to some 
$\ECST$. In particular, if $\RGT(g) \in \QQ$ for all such $\RGT$ then $\psi(g) \in \QQ$, as 
stated in (i).

3) From now on we assume $\psi(g) \not\in \QQ$, whence
there is a $\RGT$ belonging to some $\ECST$ such that $\RGT(g) \notin \QQ$. In this case, the
results of 1) imply that the order of $\theta$ is divisible by $p$. As recorded in
\cite[Remark 10.3]{GT2}, the order of $st$ is also divisible by $p$, whence $p$ divides $|s|$ as
$t$ is an $\ell$-element centralizing $s$ and $p \neq \ell$. Let $x$ denote
the $p$-part of $s$. By our assumptions, $C_{\GCD}(x)^{\circ}$ is a maximal torus. On the
other hand, by \cite[Lemma 13.14(iii), Remark 13.15(i)]{DM}, the exponent of 
$C_{\GCD}(x)/C_{\GCD}(x)^{\circ}$ divides both $|x|$ and $|Z(\GC)/Z(\GC)^{\circ}|$ and so
it is $1$ by our assumptions. Thus $C_{\GCD}(x)$ is a maximal torus. Since $C_{\GCD}(s)$ is 
contained in $C_{\GCD}(x)$ and it contains a maximal torus, we conclude that 
$C_{\GCD}(s) = C_{\GCD}(x)$. In fact, we have shown that $C_{\GCD}(st')$ is connected and it is 
equal to the maximal torus $C_{\GCD}(x)$ for any $\ell$-element $t'$ that centralizes $s$, whence 
the Lusztig series $\EC(\GCF,(st'))$ and the rational series $\EC(\GCF,(st')_{\GCDF})$ coincide. 
Hence, Lusztig's parametrization of irreducible complex characters ensures 
that the degree of any irreducible complex character $\chi$ in $\ECL(\GCF,(s))$ is just 
$D := (G:C_{\GCD}(s)^{\FD})_{r'}$ that is coprime to $q$. This in turn implies, as shown in
\cite{HM}, that the degree of any irreducible Brauer character in $\ECL(\GCF,(s))$, in
particular of $\psi$, is divisible by $D$. Now $\psi$ is an irreducible constituent of
the reduction modulo $\ell$ of some irreducible complex character $\chi$ in
$\ECL(\GCF,(s))$.  We can now conclude that $\psi = \chi (\mod \ell)$, i.e. $V$ lifts to
characteristic $0$. Thus the statement (A) has been proved.

4) Next we claim that $\chi$ is actually a Deligne-Lusztig character $\RGT$ (for
some $\TC$ and $\theta$ with $p$ dividing $|\TCF|$) up to sign. For consider any $\RGT$ in
$\ECL(\GCF,(s))$. Then $\RGT$ belongs to some $\EC(\GCF,(st))$ with $t$ an $\ell$-element
centralizing $s$. We have shown in 2) that $C_{\GCD}(st)$ is a maximal torus, whence the
relative Weyl group $W(st)/W^{0}(st)$ is trivial (cf. \cite[Remark 2.4]{DM}). Since
$(\RGT,\RGT)_{G} \leq |W(st)/W^{0}(st)|$ (see the proof of \cite[Prop. 14.43]{DM}), $\RGT$ is
irreducible (up to sign). Thus all Deligne-Luzstig characters in $\ECL(\GCF,(s))$ are
irreducible up to sign. Since $\chi$ is an irreducible constituent of some Deligne-Lusztig
character $\RGT$ in $\ECL(\GCF,(s))$, we get $\chi = \pm \RGT$. As $p$ divides $|\theta|$, 
$p$ also divides $|\TCF|$, and the claim follows.

5) Now we may assume that $\chi = \pm \RGT$ and $p$ divides $|\TCF|$. Consider any
nontrivial power $h$ of $g$. The assumptions in (B) and the arguments in 2) show that 
$C_{\GC}(h)$ is a maximal torus in $\GC$, and so $C_{G}(h)$ is a maximal torus in $G$ of order 
divisible by $p$. Among maximal tori of $G$ order divisible by $p$, choose $T$ of largest order. 
Then $|\TCF|,|C| \leq |T|$. Clearly,
$|\Ind^{\GCF}_{\TCF}(\theta)(h)| \leq (|\TCF|-1)|C|/|\TCF| < |T|$. Hence
\cite[Prop. 7.5.4]{C} implies $|\chi(h)| = |\RGT(h)| < |T|$, meanwhile
$\chi(1) = (G:\TCF)_{r'} \geq |G|_{r'}/|T| \geq |T|^{2}$. Now for any linear character
$\al$ of $A := \langle g \rangle$ we have
$$|A| \cdot (\chi|_{A},\al)_{A} = \sum_{x \in A}\chi(x)\overline{\al(x)}
  > |T|^{2} - |A| \cdot |T| > 0,$$
as $A \leq \TCF$. It follows that $\dvg = |g|$.
\end{proof}

\begin{lemma}\label{sylow}{\em \cite[\S 10.1]{GL}}
{\sl Let $\Phi_{m}(x)$ denote the $m^{\mathrm {th}}$ cyclotomic polynomial over $\QQ$,
and $\Pi \Phi_{m}^{r_{m}}(x)$ be the polynomial associated with the finite group of Lie
type $G$, see {\rm \cite[{\rm Table~   on~ p. ~111}]{GL}}. Then Sylow $p$-subgroups
of $G$ are cyclic if and only if there is exactly one $m$ such that $p$ divides
$\Phi_{m}(q)$ and $r_{m}=1$ for this $m$.\epf}
\end{lemma}

\begin{propo}\label{cyclic}
{\sl Let $\GC$ be a simple simply connected algebraic group with a Frobenius map $F$, and
let $G := \GCF$ be quasi-simple of exceptional type. Let $q$ denote the common
absolute value of eigenvalues of $F$. Assume $p$ is a prime dividing $|G|$
such that Sylow $p$-subgroups $S$ of $G$ are cyclic. Then there is an $F$-stable
simple simply connected algebraic subgroup $\LC$ of $\GC$ of rank $\geq 2$ such that
$\LCF$ is quasi-simple and contains $S$, $(p,q^{2}-1) = 1$, and moreover one of the following
holds.

{\rm (i)} $\LCF$ is exceptional, any nontrivial $p$-element in $\LCF$ is regular 
semisimple in $\LC$, and $(p,|Z(\LC)|) = 1$. 

{\rm (ii)} $\LCF = SL^{\eps}_{m}(q)$ with $m \geq 3$. Furthermore, $p$ is a primitive prime
divisor of $q^{m}-1$ if $\eps = +$ and of $q^{2m}-1$ if $\eps = -$.}
\end{propo}

\begin{proof}
1) First we show that the statements hold for $\LC = \GC$ if $p$ satisfies the
following three conditions:

(a) There is an integer $N$ with $\varphi(N) = r := \rank(\GC)$ such that
$p | (q^{N}-1)$;

(b) If $k$ is any integer with $k|N$ and $\varphi(k) < r$ then $(p,q^{k}-1) = 1$;

(c) If $\SC$ is a simple simply connected algebraic group of rank $s|r$, then
$(p,|Z(\SC)|) = 1$.

Indeed, assume (a) -- (c) hold but a nontrivial $p$-element $g \in G$ is not regular
semisimple. Since $\GC$ is simply connected, $\CL := C_{\GC}(g)$ is connected reductive.
Let $\ZC := (Z(\CL))^{o}$ and $\DC = \CL/\ZC$. By assumption, $g$ is not regular, so $\DC$
is semisimple of rank $> 0$. Also, $g \in \CL^{F}$ and
$|\CL^{F}| = |\ZC^{F}| \cdot |\DCF|$, so $p$ divides either $|\ZCF|$ or $|\DCF|$.

Notice that $p$ cannot divide $|\HCF|$ for any $F$-stable connected reductive subgroup
$\HC$ of $\GC$ of rank $< r$. Assume the contrary. Then there is an $F$-stable maximal
torus $\TC$ of $\HC$ such that $p$ divides $|\TCF|$. It is known that $F$ acts on the
character group of $\TC$ as $qF_{0}$ with $F_{0}$ of finite order $t$, and
$|\TCF| = \det(q-F)$. Any eigenvalue $\om$ of $F_{0}$ has minimal polynomial
$\Phi_{l}(x)$ over $\QQ$, with $\varphi(l) \leq \rank(\HC) < r$. It follows that
$p|(q^{l}-1)$ for some $l$ with $\varphi(l) < r$. But $p|(q^{N}-1)$, so $p|(q^{k}-1)$ for
$k := \gcd(l,N)$, contrary to (b).

The above observation implies that $\ZC = 1$. Let $\DC_{1}$ be a smallest
$F$-stable semisimple subgroup of $\CL$. Since $p$ divides either $|\DC_{1}|^{F}$
or $|(\CL/\DC_{1})^{F}|$, our observation implies that $\DC_{1} = \CL$. Thus
$\CL$ is a central quotient of $\SC^{m}$ with $\SC$ being simple simply connected
of rank $s|r$. Since $g \in Z(\CL)$, we conclude that $p||Z(\SC)|$, contrary to (c).

2) Now we can inspect the exceptional groups $G$ of Lie type.

Assume $G = G_{2}(q)$. By Lemma \ref{sylow}, $p | \Phi_{N}(q)$ with $N = 3$ or $6$, and
$(p,q^{2}-1) = 1$. In particular, $p > 3$, so (a) -- (c) hold. Thus we may take
$\LC = \GC$.

Assume $G = F_{4}(q)$. By Lemma \ref{sylow}, $p | \Phi_{N}(q)$ with $N = 8$ or $12$, and
$(p,(q^{6}-1)(q^{2}+1)) = 1$. In particular, $p > 5$, so (a) -- (c) hold. Thus we
may take $\LC = \GC$.

Assume $G = E_{6}(q)$. Here, $p | \Phi_{N}(q)$ with $N = 5$, $8$, $9$ or
$12$, and $(p,\Phi_{k}(q)) = 1$ for $k = 1,2,3,4,6$; in particular, $p > 7$.
Now we may take $\LC$ to be $\GC$ if $N = 9$, an $F$-stable subgroup of type $F_{4}$
if $N = 8$ or $12$, and an $F$-stable subgroup of type $A_{4}$ (such that
$\LCF \simeq SL_{5}(q)$) if $N = 5$.

Assume $G = \ta E_{6}(q)$. Then $p | \Phi_{N}(q)$ with $N = 8$, $10$, $12$
or $18$, and $(p,\Phi_{k}(q)) = 1$ for $k = 1,2,3,4,6$; in particular, $p > 7$.
Now we may take $\LC$ to be $\GC$ if $N = 18$, an $F$-stable subgroup of type $F_{4}$
if $N = 8$ or $12$, and an $F$-stable subgroup of type $A_{4}$ (such that
$\LCF \simeq SU_{5}(q)$) if $N = 10$.

Assume $G = E_{7}(q)$. In this case, $p | \Phi_{N}(q)$ with $N = 5$, $7$, $8$, $9$,
$10$, $12$, $14$ or $18$, and $(p,\Phi_{k}(q)) = 1$ for $k = 1,2,3,4,6$; in particular,
$p > 7$. Now we may take $\LC$ to be an $F$-stable subgroup of type $E_{6}$
if $N = 9$ or $18$, of type $F_{4}$ if $N = 8$ or $12$, of type $A_{4}$ if $N = 5$ or $10$,
and of type $A_{6}$ if $N = 7$ or $14$.

Assume $G = E_{8}(q)$. In this case, $p | \Phi_{N}(q)$ with $N = 7$, $9$, $14$, $15$,
$18$, $20$, $24$ or $30$, and $(p,\Phi_{k}(q)) = 1$ for $k = 1,2,3,4,6$; in particular,
$p > 7$. Now we may take $\LC$ to be an $F$-stable subgroup of type $E_{6}$
if $N = 9$ or $18$, of type $A_{6}$ if $N = 7$ or $14$, and $\GC$ otherwise.

Assume $G = \tb D_{4}(q)$. By Lemma \ref{sylow}, $p | \Phi_{12}(q)$ and
$(p,q^{6}-1) = 1$. In particular, $p > 5$, and we may take $\LC = \GC$.

If $G = \ta B_{2}(q^{2})$ one can take $\LC = \GC$ (cf. \cite{Bu2}). If
$G = \ta G_{2}(q^{2})$ one can take $\LC = \GC$ (cf. \cite[D.2]{H1}). Finally, assume
$G = \ta F_{4}(q^{2})$. Then $p$ divides $q^{4}-q^{2}+1$,
$q^{4}-\sqrt{2}q^{3}+q^{2}-\sqrt{2}q+1$, or $q^{4}+\sqrt{2}q^{3}+q^{2}+\sqrt{2}q+1$
(cf. \cite[D.3]{H1}), and one can take $\LC = \GC$ (cf. \cite[Prop. 4]{V}).
\end{proof}

Now we are ready to deduce two important consequences of Theorem \ref{min2}.

\begin{corol}\label{slprime}
{\sl {\rm (A)} Let $S = SL_{n}^{\eps}(q)$ with $n \geq 3$ and let $g \in S$ be any semisimple 
element of prime power order $p^{a} > 1$, with $p$ a primitive prime divisor of $q^{n}-1$ if 
$\eps = +$ and of $q^{2n}-1$ if $\eps = -$. If $V$ is a nontrivial 
irreducible $\FF S$-representation in cross characteristic $\ell$ then either 

{\rm (i)} $\ell \neq p$ and $\dvg \geq p^{a-1}(p-1)$, or

{\rm (ii)} $\ell = p$ and $\dvg \geq p^{a}-1$, or

{\rm (iii)} $\ell = p = (q^{n}-1)/(q-1)$, $S = SL_{n}(q)$, $n$ a prime, $o(g) = p$, $V$ is 
a Weil module, and $\dim(V) = \dvg = p-2$.

{\rm (B)} The conclusion of {\rm (A)} also holds if we assume that $n$ is a prime and 
$g \in S$ is any irreducible element of prime power order $p^{a}$ with $(p,q^{2}-1) = 1$.

{\rm (C)} Either {\rm (i)} or {\rm (ii)} holds, if we assume that $S = Sp_{2n}(q)$ or
$Spin^{-}_{2n}(q)$, and that $g \in S$ is any semisimple element of prime power order $p^{a} > 1$, 
with $p$ a primitive prime divisor of $q^{2n}-1$.}
\end{corol}

\begin{proof}
(A) Our assumptions imply that $p > n$ and Sylow $p$-subgroups of $S$ are cyclic. Moreover, 
$(\GC,F,p)$ satisfies all assumptions of Theorem \ref{min2}, with 
$\GC := SL_{n}(\overline{\FF}_{q})$ and $F$ chosen such that $\GCF = SL_{n}^{\eps}(q)$. 
Therefore, the result follows from \cite{Z3} if $\ell = 0$ or if $\ell = p$, and from
Theorem \ref{min2} if $\ell \neq p$ but $V$ does not lift to a complex representation. 

(B) Observe that the classification of irreducible elements (cf. Lemma \ref{icl} below) implies 
that the assumptions of (A) hold.

(C) Apply Theorem \ref{min2} to $\GC$ of type $C_{n}$ or $D_{n}$.
\end{proof}

\begin{theor}\label{min5}
{\sl Let $\GC$ be a simple simply connected algebraic group with a Frobenius map $F$, and
let $G := \GCF$ be quasi-simple of exceptional type. Assume $p$ is a prime dividing $|G|$
such that Sylow $p$-subgroups of $G$ are cyclic. Let $\FF$ be an algebraically closed
field of cross characteristic $\ell$ and let $V$ be any nontrivial irreducible
$\FF G$-module. Then for any element $g \in G$ of order $p^{b}$, $\dvg \geq p^{b-1}(p-1)$ if 
$\ell \neq p$, and $\dvg = p^{b}$ if $\ell = p$.}
\end{theor}

\begin{proof}
By the virtue of \cite{Z3} we may assume that $\ell \neq p$. Next we proceed by induction
on $|\GCF|$ and embed $g$ in $\LCF$ as in Proposition \ref{cyclic}. If we are in case (ii) 
of Proposition \ref{cyclic}, then we are done by Corollary \ref{slprime}. Otherwise $\LCF$ 
is of exceptional type; in particular, the simply connected algebraic group $(\LC^{*})_{sc}$ 
dual to $\LC$ is isomorphic to $\LC$, and $(p,|\LC|) = 1$. Let $(\LC^{*},\FD)$ be dual to 
$(\LC,F)$, and let $q$ denote the common absolute value of eigenvalues of $F$. We may now 
identify $\LC^{*}$ with $\LC/Z(\LC)$ and $\FD$ with $F$. Claim that any nontrivial $F$-stable 
$p$-element in $\LC/A$ is regular semisimple, for any $A \leq Z(\LC)$. Indeed, 
consider any nontrivial $F$-stable $p$-element $gA \in (\LC/A)^{F}$. We may assume $g \in \LC$ 
has the same order as the order $N$ of $gA$ in $\LC/A$. Since $gA$ is $F$-stable, we have 
$F(g) = gz$ for some $z \in A$. Hence $1 = F(g^{N}) = g^{N}z^{N} = z^{N}$, whence $z = 1$ as 
$(N,|A|) = 1$. Thus $g$ is $F$-stable, and so $C_{\LC}(g)$ is a maximal torus of $\LC$ by our 
assumptions. Now if $xA \in C_{\LC/A}(gA)$, then $xgx^{-1} = ga$ with $a \in A$. Comparing the 
order and using $(N,|A|) = 1$ once more, we get $a = 1$. Thus $C_{\LC/A}(g) = C_{\LC}(g)/A$ is a 
maximal torus in $\LC/A$ as stated. The case $\LCF$ is of type $\ta B_{2}$, $\ta G_{2}$, 
$\tb D_{4}(2)$, or $\ta F_{4}(2)'$, can be checked directly, so we will assume $\LCF$ is none of 
these types. Then it is easy to verify that $|\LCF|_{r'} > (q+1)^{3\rank(\LC)} \geq |T|^{3}$ for 
any maximal torus of $T$ of $\LC$, if $r$ is the prime dividing $q$. Now we can apply Theorem 
\ref{min2} to $\LC$.
\end{proof}

\section{Classical groups: Irreducible elements}

We begin with the following (well-known) description of irreducible elements in $GL_{n}(q)$:

\begin{lemma}\label{igl}
{\sl Let $G := GL_{n}(q)$, $S := SL_{n}(q)$, and let $V := \FF_{q}^{n}$ be the natural module for $G$. 
Assume $g \in G$ is an irreducible $p$-element with $o(g) = p^{a}$. Then $(p,q) = 1$, and the following
statements hold.

{\rm (i)} Up to $G$-conjugacy, $g$ is uniquely determined by its characteristic polynomial 
$p(t)$, and $p(t)$ is irreducible over $\FF_{q}$. If $1 \leq m < o(g)$ then no eigenvalue of 
$g^{m}$ can belong to $\FF_{q}$. 

{\rm (ii)} If $\lam \in \Spec(g,V)$, then $\FF_{q}(\lam) = \FF_{q^{n}}$ and 
$\Spec(g,V) = \{\lam,\lam^{q}, \ldots ,\lam^{q^{n-1}}\}$. Furthermore, 
$g$ is $G$-conjugate to $n$ distinct powers $g^{q^{i}}$, $0 \leq i \leq n-1$, of $g$.

{\rm (iii)} Let $k$ be the smallest positive integer such that $p|(q^{k}-1)$. Then 
$n/k = p^{b}$ for some $b \in \ZZ$. 

{\rm (iv)} In the notation of {\rm (iii)}, assume $b \geq 1$. Then $g$ is 
$S$-conjugate to an irreducible $p$-element $h$ of $H := GL_{p}(q^{n/p})$ (naturally 
embedded in $G$) and $p|(q^{n/p}-1)$. Assume in addition $(p,n) \neq (2,2)$. Then 
$1 \neq h^{p} \in Z(H)$ and $h^{p}$ belongs to a (proper) parabolic subgroup of $G$. Moreover, 
$g^{p} \notin Z(G)$ if $n > p$, and $o(g) = p$ if $n = p > 2$. If $n = p^{b} > 1$ and $p > 2$ 
then $g \notin SL_{n}(q)$.}
\end{lemma}

\begin{proof}
Since the unipotent part of $g$ has nonzero fixed points on $V$, $g$ must be semisimple, i.e.
$(p,q) = 1$.

(i) Inspecting the rational canonical form of $g$ on $V$, we see that the irreducibility of
$g$ on $V$ implies that $p(t)$ is also the minimal polynomial of $g$. Let $f(t)$ be any 
irreducible divisor of $p(t)$. Then $f(g)V$ is $g$-invariant, and the action of $g$ on $f(g)V$ is
annihilated by $p(t)/f(t)$. Hence $f(g)V \neq V$, and so $f(g)V = 0$ by irreducibility. It
follows that $p(t) = f(t)$ is irreducible. Since $p(t)$ is the unique invariant factor of $g$, 
$g$ is uniquely determined by $p(t)$ up to conjugacy. The second statement follows by 
irreducibility.

(ii) Clearly, $\lam$ is a root of $p(t)$. But $p(t)$ is irreducible of degree $n$ over
$\FF_{q}$, so $\FF_{q}(\lam) = \FF_{q^{n}}$. All the roots of $p(t)$ are 
$\lam^{q^{i}}$ with $0 \leq i \leq n-1$, and all of them belong to $\Spec(g,V)$ since 
$g \in GL_{n}(q)$. In a basis of $V \otimes \overline{\FF}_{q}$ we may assume that 
$g$ is represented by $\diag(\lam,\lam^{q}, \ldots ,\lam^{q^{n-1}})$. Then the $n$ powers
$g^{q^{i}}$ with $0 \leq i \leq n-1$ are distinct, and all of them are conjugate to
$g$ by (i), since they have the same characteristic polynomial $p(t)$.

(iii) In the notation of (ii) we see that $|g| = |\lam|$ divides $q^{n}-1$; in particular,
$p|(q^{n}-1)$ and so $k|n$. Write $n/k = p^{b}s$ with $s \in \ZZ$ coprime to $p$. Then 
$(q^{n}-1)/(q^{n/s}-1) \equiv s (\mod p)$. It follows that the multiplicative groups of 
$\FF_{q^{n}}$ and $\FF_{q^{n/s}}$ have the same Sylow $p$-subgroups, whence 
$\lam \in \FF_{q^{n/s}}$. By (ii) we must have $s = 1$.

(iv) Set $Q := q^{n/p}$. Since $\FF_{Q^{p}} = \FF_{q^{n}}$, there is an element $h \in H$ 
with one eigenvalue (on $W := \FF_{Q}^{p}$) equal to $\lam$. When one embeds $H$ in $G$ 
naturally (by viewing $W$ as an $\FF_{q}$-space and identifying it with $V$), $h$ has the 
characteristic polynomial $p(t)$ on $V$ and so $g$ and $h$ are $G$-conjugate by (i).
Next, $HS = G$, so $g$ is $S$-conjugate to some $H$-conjugate of $h$. So, up to 
$S$-conjugacy, we can embed $g$ in $H$. Since $g$ is irreducible on $V$, $h$ is irreducible 
on $W$. In a basis of $W \otimes \overline{\FF}_{q}$ we may assume that 
$h$ is represented by $\diag(\lam,\lam^{Q}, \ldots ,\lam^{Q^{p-1}})$. Since 
$\lam$ is primitive in $\FF_{Q^{p}}$, $\lam^{p} \neq 1$ and so $h^{p} \neq 1$. Assume 
$(p,n) \neq (2,2)$. We show that $\lam^{(Q-1)p} = 1$, which implies that $h^{p}$ is scalar and 
so $h^{p} \in Z(H)$. Indeed, if $p$ is odd, then $(Q^{p}-1)/(Q-1) \equiv p (\mod p^{2})$,
whence $(Q^{p}-1)_{p} = p(Q-1)_{p}$ and so $|\lam|$ divides $p(Q-1)$. If $p = 2$ and 
$b \geq 2$ then $4|(Q-1)$ and we again have $(Q^{p}-1)_{p} = p(Q-1)_{p}$. If $p = 2$ and 
$b = 1$, then $k = 1$ (by the choice of $k$) and so $n = 2$, contrary to our assumption.
Thus $h^{p} \in Z(H)$. In particular, $h^{p}$ stabilizes a $1$-space of $W$ and so it 
fixes an $n/p$-dimensional subspace of $V$, whence it belongs to a parabolic subgroup of $G$.
Next we assume $g^{p} \in Z(G)$. Then $|\lam| = |g|$ divides $p(q-1)_{p}$. If $k > 1$, then 
$(p,q-1) = 1$ by the choice of $k$, so $|\lam| = p$ and $\FF_{q}(\lam) \subseteq \FF_{q^{n/p}}$,
contrary to (ii). Thus $k = 1$. In this case, $p(q-1)$ divides $q^{p}-1$, so (ii) implies
$n = p$. Conversely, assume $n = p > 2$. Then $p(q-1)_{p} = (q^{p}-1)_{p}$, whence 
$\lam^{p} = \lam^{pq}$ and $g \notin Z(G) \ni g^{p}$, i.e. $o(g) = p$. Finally, assume 
that $n = p^{b} > 1$ and $p > 2$ but $g \in SL_{n}(q)$. Then $k = 1$, $p|(q-1)$,  
$((q^{n}-1)/(q-1))_{p} = p^{b}$ and $1 = \det(g) = \lam^{(q^{n}-1)/(q-1)}$, yielding 
$\lam^{p^{b}} = 1$. In this case $\FF_{q}(\lam) \subseteq \FF_{q^{n/p}}$, contrary to (ii).
\end{proof}

In what follows, $Sp_{n}(q)$ is considered only for even $n$, $GO_{2m+1}(q)$ is considered
only for odd $q$, and $GU_{n}(q^{1/2})$ is considered only for $q$ a square. (The reason
we consider $GU_{n}(q^{1/2})$ instead of $GU_{n}(q)$ in what follows is that we want to use 
$q$ for the size of the defining field. This allows us to have a uniform argument for both 
$SL$- and $SU$-cases.) 
 
\begin{propo}\label{icl}
{\sl Let $G \in \{ GU_{n}(q^{1/2}), CSp_{n}(q), GO^{\eps}_{n}(q)\}$, $S := [G,G]$, and let
$V := \FF^{n}_{q}$ be the natural module for $G$. Assume $g \in G$ is an irreducible 
$p$-element with $o(g) = p^{a}$. Then $(p,q) = 1$ and the following statements hold.

{\rm (i)} Let $k$ be the smallest positive integer such that $p|(q^{k}-1)$. Then 
$n/k = p^{b}$ for some $b \in \ZZ$. Set $G_{1} = Sp_{n}(q)$ if $G = CSp_{n}(q)$ and 
$G_{1} = G$ otherwise. Then up to $G_{1}$-conjugacy, $g$ is uniquely determined by its 
characteristic polynomial $p(t)$, and $p(t)$ is irreducible over $\FF_{q}$. Furthermore, 
$g$ is $G_{1}$-conjugate to $n$ distinct powers $g^{q^{i}}$, $0 \leq i \leq n-1$, of $g$.

{\rm (ii)} Assume $G = GU_{n}(q^{1/2})$. Then $n$ is odd. Assume in addition that $b \geq 1$. 
Then $g$ is $S$-conjugate to an irreducible $p$-element $h$ of $H := GU_{p}(q^{n/2p})$ 
(naturally embedded in $G$), $p > 2$, and $p|(q^{n/2p}+1)$. Moreover, $1 \neq h^{p} \in Z(H)$, 
$h^{p}$ belongs to a parabolic subgroup of $G$. If $n > p$ then $o(g) > p$, and if $n = p$ then 
$o(g) = p$. If $n = p^{b} > 1$ then 
$g \notin SU_{n}(q^{1/2})$. 

{\rm (iii)} Assume $G = CSp_{n}(q)$ with $n > 2$. Then $n = 2m$, $p > 2$, $p|(q^{m}+1)$, and 
$g$ is $S$-conjugate to an element in a torus of order $q^{m}+1$ of $SL_{2}(q^{m})$ 
(naturally embedded in $Sp_{n}(q)$).

{\rm (iv)} Assume $G = GO^{\eps}_{n}(q)$ with $n > 1$. Then $n$ is even, $\eps = -$, 
$g^{q^{n/2}+1} = 1$, and either $p > 2$, or $p = n = 2$. Assume in addition that $n > 2$. 
If $b \geq 1$, then $g$ can be embedded as an irreducible $p$-element $h$ of 
$H := GU_{p}(q^{n/2p}) < GO^{-}_{2p}(q^{n/2p})$ (naturally embedded in $G$) 
and $p|(q^{n/2p}+1)$; moreover, $h^{p} \in Z(H) \setminus Z(G)$ and $h^{p}$ belongs to a parabolic 
subgroup of $G$. Furthermore, $g$ is $S$-conjugate to an irreducible $p$-element in a subgroup 
$GU_{n/2}(q)$ of $G$ if $n \equiv 2 (\mod 4)$, and in a torus $\ZZ_{(q^{n/2}+1)/(2,q-1)}$ of a 
subgroup $PSL_{2}(q^{n/2})$ of $G$ if $n \equiv 0(\mod 4)$.}
\end{propo}

\begin{proof}
(i) Clearly, $g$ is an irreducible element in $A := GL_{n}(q)$. So $(p,q) = 1$. 
Semisimple elements in classical groups are classified up to conjugacy in 
\cite{Wa}. In particular, two irreducible elements of $G_{1}$ are conjugate if and only if
they are conjugate in $A$, cf. \cite[pp. 34 - 39, 59]{Wa}. We will show in (iii) that 
all irreducible $p$-elements of $CSp_{n}(q)$ are contained in $S$. Hence the statements 
follow from Lemma \ref{igl}.

(ii) It is well known that $n$ is odd and $|g|$ divides 
$q^{n/2}+1$. Assume $b \geq 1$; in particular, $p$ is odd. Note that $p|(q^{k/2}+1)$ 
and so $p|(q^{n/2p}+1)$ as $p$ is odd. (Indeed, if $p$ does not divide $q^{k/2}+1$ then $k$ 
divides $q^{k/2}-1$ and $q^{n/2}+1$, whence $p = 2$, a contradiction.) Fixing 
$\lam \in \Spec(g,V)$, we have $\lam^{q^{n/2}+1} = 1$. Using the natural embeddings
$GU_{1}(q^{n/2}) < H < G$, we can find an element $h \in H$ with an eigenvalue equal to $\lam$. 
Arguing as in the proof of Lemma \ref{igl}(ii) we obtain that $g$ is $S$-conjugate to an 
$H$-conjugate of $h$. Embedding $g$ in $GL_{p}(q^{n/p})$ we get $1 \neq h^{p} \in Z(H)$ by Lemma 
\ref{igl}(iv). Clearly, the central element $h^{p}$ of $H$ stabilizes a singular $1$-space $L$ 
of the Hermitian space $\FF_{q^{n/p}}^{p}$. Notice that $L$ becomes a totally singular 
$n/p$-dimensional subspace of the Hermitian space $V$, whence $h^{p}$ belongs to a parabolic 
subgroup of $G$. The last two claims follow from Lemma \ref{igl}(iv).     

(iii) In our case $n = 2m$ is even. It is shown in \cite[pp. 2127, 2130]{TZ1} that $|g|$ 
divides $(q-1)(q^{m}+1)$. Notice that $(p,q-1) = 1$. (Otherwise we would have $k = 1$ and so 
$n = p^{b}$ by (i). If $p > 2$ then $|\lam|$ divides $(q-1)(q^{m}+1)_{p} = q-1$. If $p = 2$,
then $q$ is odd; furthermore, $m$ is even as $n = 2^{b} > 2$, whence $|\lam|$ divides 
$(q-1)(q^{m}+1)_{2} = 2(q-1)$. In both of the cases, $|\lam|$ divides $q^{2}-1$, 
violating the condition $\FF_{q}(\lam) = \FF_{q^{n}}$.) This in turn implies that 
$|g|$ divides $q^{m}+1$ and it is shown in \cite[pp. 2127, 2130]{TZ1} that $g$ is 
$S$-conjugate to an element in a torus of order $q^{m}+1$ of a natural subgroup 
$SL_{2}(q^{m})$ of $Sp_{2n}(q)$. Finally, we notice that $p > 2$. (Indeed, if $p = 2$ then 
$q$ is odd and $p|(q-1)$, a contradiction.)   

(iv) If $n$ is even, then the irreducibility of $g$ implies $\eps = -$, cf. 
\cite[pp. 38, 59]{Wa}. In general, the map $\lam \mapsto \lam^{-1}$ acts on $\Spec(g,V)$ as 
an involution without fixed points (as $g$ is irreducible), so $n$ is even. So we may assume 
$G = GO^{-}_{2m}(q)$. It is well known that $|g|$ divides $q^{m}+1$. From 
now on assume $n > 2$. Arguing as in (iii) we get $p > 2$. Fixing 
$\lam \in \Spec(g,V)$, we have $\lam^{q^{m}+1} = 1$. 

Suppose $b \geq 1$. Using the natural embeddings
$GU_{1}(q^{m}) < H < GO^{-}_{2p}(q^{m/p}) \leq G$, we can find an element $h \in H$ with an 
eigenvalue equal to $\lam$. Arguing as in (ii) we obtain that $h$ and $g$ are 
$G$-conjugate, $p|(q^{m/p}+1)$ and $1 \neq h^{p} \in Z(H)$. Clearly, the central element 
$h^{p}$ of $H$ stabilizes a singular $1$-space $L$ of the Hermitian space $\FF_{q^{2m/p}}^{p}$.
Notice that $L$ becomes a totally singular $2$-space of the orthogonal space $\FF_{q^{m/p}}^{2p}$,
which in turns is a totally singular $2m/p$-dimensional subspace of the orthogonal space $V$.
Hence $h^{p}$ belongs to a parabolic subgroup of $G$. Since 
$n = 2m \geq 2p$, $g^{p} \notin Z(G)$ by Lemma \ref{igl}(iv). 

If $n \equiv 2 (\mod 4)$, then using the natural embeddings 
$GU_{1}(q^{m}) < K := GU_{m}(q) < G$ we can find an element $t \in K$ with an eigenvalue 
equal to $\lam$, and so $g$ and $t$ are $G$-conjugate as above. We claim that $g$ is 
$S$-conjugate to $t$ or $t^{q}$. Firstly, the Frobenius map (that raises every element 
of $K$ to its $q$-power) normalizes $K$, sends $t$ to (a $K$-conjugate of) $t^{q}$, and it can 
be realized as an element of $G \setminus SO^{-}_{n}(q)$ if $q$ is odd, and of 
$G \setminus S$ if $q$ is even (as $n/2$ is odd). Secondly, if $q$ is odd then any generator 
of $Z(K)$ belongs to $SO^{-}_{n}(q) \setminus S$ by \cite[Lemma 8.14]{TZ4} and centralizes $t$,
so we are done. 
 
If $n \equiv 0 (\mod 4)$, we can argue similarly, using the embeddings 
$\ZZ_{(q^{m}+1)/e} < PSL_{2}(q^{m}) \simeq \Omega^{-}_{4}(q^{m/2}) < G$ with $e := (2,q-1)$. 
In more details, one can choose a primitive prime divisor $r$ of $q^{2m}-1$ and a Sylow 
$r$-subgroup $R$ of $PSL_{2}(q^{m})$. Then $R$ is also a Sylow $r$-subgroup of $G$ (and $S$), and
$C_{G}(R) = \ZZ_{q^{m}+1}$. Hence the $S$-conjugacy follows. 
\end{proof}

A crucial rule in the treatment of irreducible $p$-elements is played by the following 
statement.

\begin{lemma}\label{root2} 
{\sl Let $p > 2$ be a prime, $b \geq 1$, $n = p^{b}$, and assume that either  
$SL_{n}(q) \lhd G \leq GL_{n}(q)$, or $SU_{n}(q^{1/2}) \lhd G \leq GU_{n}(q^{1/2})$. Let 
$g \in G$ be an irreducible $p$-element.

{\rm (i)} Then $o(g) = p^{b}$.

{\rm (ii)} Let $\ell \neq p$, $H$ a finite group with a central $p'$-subgroup $C$ such that 
$H/C \simeq G$ and let $h$ be a preimage of $g$ in $H$. Let $\Phi$ be an irreducible 
$\FF H$-representation of degree $>1$. Assume in addition that $(n,q) \neq (3,4)$ in the 
$SU$-case. Then $\deg(\Phi(h)) = p^{b} = o(h)$. In fact, $\Spec(\Phi(h))$ consists of all
$p^{b}$-roots of the unique eigenvalue of $\Phi(h^{p^{b}})$.}
\end{lemma}

\begin{proof} 
{\rm (i)} Clearly, $(p,q) = 1$ and so $p|(q^{p-1}-1)$. By Lemma \ref{igl}(iv) and 
Proposition \ref{icl}(i), $p|(q-1)$. Let $V = \FF_{q}^{n}$ denote the natural 
$G$-module, and let $\lam$ be any eigenvalue of $g$ on $V$. Then $|\lam|$ divides
$(q^{n}-1)_{p} = p^{b+c}$ if $p^{c} := (q-1)_{p}$. In fact $|\lam| = p^{b+c}$ as otherwise 
$\FF_{q}(\lam) \leq \FF_{q^{n/p}}$, contradicting Lemma \ref{igl}(ii). Thus 
$\al := \lam^{n} \in \FF_{q}$. In the $SU$-case we even have $\al^{q^{1/2}+1} = 1$
as $p|(q^{1/2}+1)$ by Proposition \ref{icl}(ii). Moreover, by Lemma \ref{igl}(i) and
Proposition \ref{icl}(i), up to $\tG$-conjugacy, $g$ is uniquely determined by 
the minimal polynomial $t^{n}-\al$ of $\lam$ over $\FF_{q}$, where we set 
$\tG = GL_{n}(q)$ in the $SL$-case and $\tG = GU_{n}(q^{1/2})$ in the $SU$-case. It is 
easy to see that $o(g) = p^{b}$.  

Now we choose a basis $(v_{1}, \ldots ,v_{n})$ of $V$, which is orthonormal in the 
$SU$-case. Define a map $x \in GL(V)$ via $x(v_{i}) = v_{i+1}$ for $1 \leq i \leq n-1$
and $x(v_{n}) = \al v_{1}$. Then $x$ belongs to $\tG$ and is annihilated by $t^{n}-\al$, 
whence $g$ and $x$ are $\tG$-conjugate. But $G \lhd \tG$, so $x \in G$. Since 
$\det(x) = \al$, we conclude that $G$ contains all matrices that lie in $\tG$ and have 
determinant a power of $\al$. Observe that $Z(\tG)G = \tG$. Indeed, 
$G \geq [\tG,\tG]$ and the $p'$-part of $\tG/[\tG,\tG]$ is covered by $Z(\tG)$. On the other hand, 
$\det(g) = \al$ generates the $p$-part of $\FF_{q}^{\bullet}$, whence the $p$-part of 
$\tG/[\tG,\tG]$ is covered by $x$. It follows $x$ and $g$ are $G$-conjugate, and we will
identify $g$ with $x$.  

Next we consider the subgroup 
$$A = \ZZ_{p^{c}}^{n} = \left\{ \diag(x_{1}, \ldots ,x_{n}) \mid x_{i} \in \FQ, x_{i}^{p^{c}} = 1
  \right\}.$$ 
Clearly, $A \leq G$ and $A$ is normalized by $g$. Moreover, $g^{p^{b}} \in C_{G}(A)$ but 
$g^{p^{b-1}} \notin C_{G}(A (\mod Z(G)))$.  

(ii) Our assumptions imply that $G$ and $H$ are not solvable, and 
$\Ker(\Phi)\leq Z(H)$. Since $C$ is a $p'$-subgroup,
without loss we may assume that $g$ and $A$ are contained in $H$, $h = g$, and 
$g \in N_{H}(A)$, $g^{p^{b}} \in C_{H}(A)$ but 
$g^{p^{b-1}} \notin C_{H}(A (\mod Z(H)))$. Clearly, $g$ permutes the $A$-eigenspaces on
the representation space of $\Phi$, and any cycle of this permutation has length dividing 
$p^{b}$. If no cycle has length $p^{b}$ then $g^{p^{b-1}}$ centralizes $A$ modulo $\Ker(\Phi)$,
a contradiction. Hence some cycle has length $p^{b}$, and so $\deg(\Phi(g)) = p^{b}$ and 
$\Spec(\Phi(g))$ consists of all $p^{b}$-roots of the unique eigenvalue of $\Phi(g^{p^{b}})$.      \end{proof}
 
\begin{theor}\label{slsup} 
{\sl Let $SL^{\eps}_{p}(q) \leq G \leq GL^{\eps}_{p}(q)$ and let $g \in G$ be an 
irreducible $p$-element with $p > 2$. Let $\Theta \in \IBRL(G)$ with 
$\dim(\Theta) >1$ and $(\ell,q) = 1$. Then $\deg(\Theta(g)) = o(g) = p$, except when 
$(p,q,\eps) = (3,2,-)$.} 
\end{theor}
 
\begin{proof}
1) By virtue of Lemma \ref{root2}(ii), we may assume that $\ell = p$. Since 
$g \notin SL^{\eps}_{p}(q)$, by Lemma \ref{easy} without loss we may 
assume that $G = GL^{\eps}_{p}(q)$. Furthermore, the proof of Lemma \ref{root2} shows that $g$ 
can be embedded in the subgroup $X := \ZZ_{q-\eps} \wr \SSS_{p}$ (consisting of all the monomial
matrices in the basis $(e_{1}, \ldots ,e_{p})$ specified in that proof). In fact we may and 
will identify $g$ with the element $x$ constructed in that proof. 
First we consider the case where $q-\eps$ is not a $p$-power. Then we can find a prime 
$r \neq p$ that divides $q-\eps$, and consider the diagonal subgroup $R := \ZZ_{r}^{p}$ inside
$X$ (that is, any element in $R$ multiplies $e_{i}$ by an $r$-root of unity). Observe that
$g$ normalizes $R$ but $gZ(G)$ does not centralize $R$ in $G/Z(G)$. A standard consideration
of $\Theta(R)$-eigenspaces shows that $\deg(\Theta(g)) = o(g) = p$.

2) Now we may assume that $q-\eps = p^{c}$ for some integer $c$. Here we consider
the case $p = 3$ and $(q,\eps) \neq (2,-)$. Setting $K := \Ker(\Theta) \geq O_{3}(Z(G))$, we get 
$G/K = Z \times PGL^{\eps}_{3}(q)$ with $Z$ a central $3'$-subgroup. Applying \cite{Ch} to 
$G/K$, we get $d_{\Theta}(g) > 2$ and so $d_{\Theta}(g) = o(g)$. 

3) From now on we may assume that $p \geq 5$. By way of contradiction,
assume that $\deg(\Theta(g)) < p$, and let $\Phi$ be any irreducible constituent of 
$\Theta|_{X}$ of degree $> 1$. Then $\Phi|_{O_{p}(X)}$ is trivial and so $\Phi$ is an 
irreducible representation of $\SSS_{p}$; moreover, $\Phi(g)$ is just the image of 
some $p$-cycle in $\SSS_{p}$ under $\Phi$. Since $\deg(\Phi(g)) < p$, the main result of
\cite{KZ} implies that $\dim(\Phi) = \deg(\Phi(g)) = p-2$ and $\Phi|_{\AAA_{p}}$ is the heart 
of the natural permutation module. 

4) Assuming furthermore that $(p,q,\eps) \neq (5,4,-)$, we see that $p \geq 7$.  
Now we choose $A$ to be a standard subgroup $\AAA_{p-2}$ of $\SSS_{p}$ (that permutes 
the basis vectors $e_{2}, \ldots ,e_{p-1}$ and fixes $e_{1}$ and $e_{p}$). Then 
$\Phi|_{A}$ is the sum of the trivial module and the natural permutation module for $A$. 
We have shown that any irreducible constituent of $\Theta|_{A}$ is either trivial, or 
the heart of the natural permutation module for $A$ (and has dimension $p-3$). Let $\rho$ 
denote the Brauer character of the latter.

Consider the parabolic subgroup $P$ of $G$, which fixes $\la e_{1} \ra_{\FF_{q}}$ and 
$\la e_{1}, \ldots ,e_{p-1}\ra_{\FF_{q}}$ in the case $\eps = +$, and fixes a singular $1$-space 
inside $\la e_{1},e_{p}\ra_{\FF_{q^{2}}}$ if $\eps = -$. Then $A$ embeds in the subgroup 
$K := SL^{\eps}_{p-2}(q)$ in the Levi subgroup of $P$. One can show (see e.g. \cite{TZ2}) that 
there is a constituent $\Lambda$ of $\Theta|_{P}$ such that the Brauer character $\lam$ of 
$\Lambda|_{K}$ is the product of the so-called {\it reducible Weil character} of $K$ and 
another Brauer character $\om$ of $K$. That is, $\lam(h) = \eps(\eps q)^{e(h-1)}\om(h)$, for 
any $p'$-element $h \in K$ and with $e(y)$ being the dimension of the kernel of the 
transformation $y$ acting on the natural module of $K$. 

On the other hand, we have shown above that there are integers $a, b \geq 0$ such that 
$\lam(h) = a + b\rho(h)$. Equating the two formulae for $\lam(h)$, where $h$ is the identity and 
a $3$-cycle $h_{1}$ (recall that $A = \AAA_{p-2}$ with $p \geq 7$), we get 
$$a + b(p-3) = q^{p-2}\om(1),~~~a + b(p-6) = q^{p-4}\om(h_{1})\,.$$     
Since $p \geq 7$, we get $q^{2} \geq (p-1)^{2} \geq 5p+1$. It follows that 
$$a + b(p-3) = q^{p-2}\om(1) \geq q^{2} \cdot |q^{p-4}\om(h_{1})| \geq q^{2}(a+b(p-6)) 
 \geq (5p+1)(a+b),$$
whence $a = b = 0$ and $\om(1) = 0$, a contradiction.

5) From now on we will assume that $G = GU_{5}(4)$ and that $\deg(\Theta(g)) < 5$. Let 
$\vartheta$ denote the Brauer character of $\Theta$. Recall
we have embedded $g$ in the subgroup $X := \ZZ_{5} \wr \SSS_{5}$ consisting of all the monomial 
matrices in an orthonormal basis $(e_{1}, \ldots ,e_{5})$ of the natural module 
$V = \FF_{16}^{5}$ of $G$. Choosing $B$ to be the standard 
subgroup $\AAA_{5}$ in $X$ and arguing as in 3), we see that every 
irreducible constituent of $\vartheta|_{B}$ is either trivial, or equal to $\beta$, the unique 
irreducible $5$-Brauer character of degree $3$ of $B$.

6) Now we will construct a chain of embeddings $B < C < S < G$, with $C = Sp_{4}(4)$ and 
$S = SU_{4}(4)$. First we consider the permutation $\FF_{2}B$-module 
$M_{1} := \la v_{1}, \ldots ,v_{5}\ra_{\FF_{2}}$ (with $B$ permuting the $v_{i}$'s 
naturally), and let 
$M := \{\sum^{5}_{i=1}x_{i}v_{i} \mid x_{i} \in \FF_{2},\sum^{5}_{i=1}x_{i} = 0\}$. Then $M_{1}$ 
supports a $B$-invariant bilinear form $(\cdot,\cdot)$: $(v_{i},v_{j}) = \delta_{ij}$, 
whose restriction to $M$ is nondegenerate. Choosing $W := M \otimes_{\FF_{2}}\FF_{4}$ and 
extending $(\cdot,\cdot)$ to $W$, we get an embedding $B < C = Sp(W)$. 
Claim that, under this embedding, the involution 
$t := (12)(34) \in B$ belongs to the conjugacy class $2C$ of $C$ in the notation of \cite{Atlas}.
Indeed, fix the basis 
$$f_{1} := v_{1}+v_{2}+v_{3}+v_{4},~~f_{2} := v_{1}+v_{2},~~f_{3} := v_{1}+v_{3},~~
  f_{4} := v_{4}+v_{5}$$
of $M$. Then $t$ is represented in this basis by the matrix
$\begin{pmatrix}1 & 0 & 1 & 1\\0 & 1 & 0 & 1\\0 & 0 & 1 & 0\\0 & 0 & 0 & 1 \end{pmatrix}$,
which is an element of class $A_{32}$ in the notation of \cite{Eno}. Using the character 
table of $C$ as given in \cite{Eno}, we see that $\chi(t) = 2$ if $\chi \in \Irr(C)$ has degree 
$18$. Inspecting the character table of $C$ as given in \cite{Atlas}, we conclude that $t$ 
belongs to class $2C$ as stated. Also, by checking the trace of a $3$-cycle $y$ in $B$ while 
acting on $M$ and $W$, we may assume that $y$ belongs to class $3A$, in the notation of 
\cite{Atlas}.   

Next we define the nondegenerate Hermitian form $\circ$ on $U := W \otimes_{\FF_{4}}\FF_{16}$ by 
$(\alpha u) \circ (\beta v) = \alpha^{4}\beta (u,v)$ for $u,v \in W$ and 
$\alpha,\beta \in \FF_{16}$. This yields an embedding $C < S = SU(U)$. Finally, 
identifying $U$ with the orthogonal complement to $\sum^{5}_{i=1}e_{i}$ in $V$, as well as
$f_{1}$ with $e_{1}+e_{2}+e_{3}+e_{4}$, $f_{2}$ with $e_{1}+e_{2}$, $f_{3}$ with 
$e_{1}+e_{3}$, and $f_{4}$ with $e_{4}+e_{5}$, we get an embedding $S < G$. The identifications
we made ensure that the actions of $B$ on $M$ and on $V$ are compatible.      

7) Let $\varphi$ be any irreducible constituent of $\vartheta|_{C}$. According to 5),
$\varphi|_{B} = a + b\beta$ for some $0 \leq a, b \in \ZZ$. Since 
$\beta(1) = 3$, $\beta(t) = -1$, $\beta(y) = 0$, we arrive at the equality 
$\varphi(1) + 3\varphi(2C) - 4\varphi(3A) = 0$, where we have deliberately used the notations
of some classes in $C$ to denote representatives of the classes. Inspecting the $5$-Brauer
characters of $C$ \cite{JLPW}, we can verify that this equality forces $\varphi = 1_{C}$. It 
follows that $C \leq \Ker(\Theta)$ and so $\Theta(1) = 1$.   
\end{proof}

Observe that the condition $G \not\leq GU_{3}(2)$ in Lemma \ref{root2} and Theorem 
\ref{slsup} cannot be removed; cf. Proposition \ref{su31} below. 

For the next section we will need the following assertion about Weil modules:

\begin{lemma}\label{weil1}
{\sl Let $G := GU_{p^{b}}(q)$ with $q+1 = p^{c}$ for some odd prime $p$, $g \in G$ an 
irreducible $p$-element, and let $V$ be an (irreducible) Weil module of $G$ in characteristic 
$\ell$ coprime to $q$. Assume $(p^{b},q) \neq (3,2)$. Then $\dvg = o(g)$.}
\end{lemma}

\begin{proof}
By Lemma \ref{root2}(i), $o(g) = p^{b}$. Furthermore, by Lemma \ref{root2} and Theorem 
\ref{slsup} we may assume that $\ell = p$ and $b > 1$.    
Let $\zeta$ denote the Brauer character of $V$. 
Since $q+1 = p^{c}$, we see that $\zeta$ is just the restriction of the complex Weil character
$\zeta^{0}_{p^{b},q}$ as defined in \cite{TZ2} to $p'$-elements (cf. \cite{DT}). 
In particular, for any $p'$-element $x \in G$ we have
$$\zeta(x) = -\frac{1}{q+1}\sum^{q}_{l=0}(-q)^{e(x-\delta^{l})},$$
where $\delta$ is a fixed primitive $(q+1)^{\mathrm {th}}$-root of unity in $\FF_{q^{2}}$ and 
$e(y)$ is the dimension of the kernel of the transformation $y$ acting on $\FF_{q^{2}}^{p^{b}}$.

Setting $Q := q^{p^{b-1}}$, by Proposition \ref{icl}(ii) 
we may embed $g$ in a natural subgroup $H := GU_{p}(Q)$ of $G$. First we want to find all 
degree $1$ composition factors of $\zeta|_{H}$. Recall (cf. \cite{TZ2}) that the reducible complex 
Weil character of $G$, resp. of $H$, is given by $\om(x) = -(-q)^{e(x-1)}$, resp. by 
$\om'(x) = -(-Q)^{E(x-1)}$, where $E(y)$ is the dimension of the kernel of the transformation 
$y$ acting on the $\FF_{Q^{2}}$-space $\FF_{Q^{2}}^{p}$. Observe that 
$\om' = \om|_{H}$. Next, let $t$ be a generator of $Z(H) \simeq \ZZ_{Q+1}$ and 
$z := t^{(Q+1)/(q+1)}$. Also fix a primitive $(Q+1)^{\mathrm {th}}$-root $\tau$ of unity in 
$\CC$ and let $\xi := \tau^{(Q+1)/(q+1)}$. By its definition (cf. \cite{TZ2}), for 
$0 \leq i \leq q$, the Weil character $\zeta^{i}_{p^{b},q}$ of $G$ is afforded by the 
$\xi^{i}$-eigenspace of $z$ in the representation space of $\om$. Similarly, for 
$0 \leq j \leq Q$, the Weil character $\zeta^{j}_{p,Q}$ of $G$ is afforded by the 
$\tau^{j}$-eigenspace of $t$ in the representation space of $\om'$. It follows that 
$$\zeta^{0}_{p^{b},q}|_{H} = \sum_{0 \leq j \leq Q,~(q+1)|j}\zeta^{j}_{p,Q}\,.$$
Notice that $Q+1 = p^{b+c-1}l$ with $(l,p) = 1$. By \cite{DT, HM}, 
$\zeta^{j}_{p,Q} (\mod p)$ is irreducible if $l \not{|}j$ or if $j = 0$, whereas if 
$l|j$ and $j \neq 0$ then $\zeta^{j}_{p,Q} (\mod p)$ is the sum of $1_{H}$ and an irreducible 
Brauer character of degree $> 1$. Thus the number of degree $1$ composition factors of 
$\zeta|_{H}$ is exactly $(Q+1)/l(q+1)-1 = p^{b-1}-1$. 

4) By Proposition \ref{icl}(ii) we can apply the main result of \cite{DZ1} to $g^{p}$
and obtain $d_{V}(g^{p}) = o(g^{p}) = p^{b-1} =: m$. Consider the filtration 
$V = V_{m} \supset V_{m-1} \supset \ldots \supset V_{1} \supset 0$, where 
$V_{k} := \Ker((g^{p}-1)^{k})$. Also by Proposition \ref{icl}(ii), $H \leq C_{G}(g^{p})$. 
If the $H$-module $V_{m}/V_{m-1}$ has a composition factor of degree $1$, then so does every 
quotient $V_{k}/V_{k-1}$ by Lemma \ref{filtr}(i), and so the number of degree $1$ composition 
factors of $\zeta|_{H}$ is at least $m$, contradicting the result of 3). Thus all composition 
factors of the $H$-module $V_{m}/V_{m-1}$ are of degree $> 1$. Notice that 
$(p,Q) \neq (3,2)$, $(5,4)$ as $b \geq 2$. Now $g$ is an irreducible $p$-element in $H$ of order
$p$ modulo $Z(H)$, so $d_{V_{m}/V_{m-1}}(g) \geq p$ by Theorem \ref{slsup}. It follows by
Lemma \ref{filtr}(ii) that $\dvg \geq (m-1)p + p = p^{b} = o(g)$.               
\end{proof}

Next we prove the following analogue of Lemma \ref{root2} for $p = 2$.

\begin{lemma}\label{root3} 
{\sl Let $n = 2^{b} \geq 2$, $SL_{n}(q) \lhd G \leq GL_{n}(q)$, $(n,q) \neq (2,3)$, and $\ell \neq 2$.
Assume that $\Phi$ is an irreducible $\FF G$-representation of degree $>1$ and $g \in G$ is an 
irreducible $2$-element with $o(g) = 2^{a}$. Then 
$\deg(\Phi(g)) \geq 2^{a-1}+1$. In fact, $\deg(\Phi(g)) = 2^{a}$ if $q \equiv 1 (\mod 4)$.}
\end{lemma}

\begin{proof} 
(i) Here we consider the case $q \equiv 1 (\mod 4)$.
Let $V = \FF_{q}^{n}$ denote the natural $G$-module, and let $\lam$ be any eigenvalue of 
$g$ on $V$. Then $|\lam|$ divides $(q^{n}-1)_{p} = 2^{b+c}$, where $2^{c} := (q-1)_{2} \geq 4$. 
In fact $|\lam| = 2^{b+c}$ as otherwise $\FF_{q}(\lam) \leq \FF_{q^{n/2}}$, contradicting 
Lemma \ref{igl}(ii). Thus $\al := \lam^{n} \in \FF_{q}$. Moreover, by Lemma \ref{igl}(i), up to 
$\tG$-conjugacy $g$ is uniquely determined by the minimal 
polynomial $t^{n}-\al$ of $\lam$ over $\FF_{q}$, where we set $\tG = GL_{n}(q)$. It is easy 
to see that $o(g) = 2^{b}$ in this case.  

Arguing as in the proof of Lemma \ref{root2}(i), we can find a basis 
$(v_{1}, \ldots ,v_{n})$ of $V$, such that $g(v_{i}) = v_{i+1}$ for $1 \leq i \leq n-1$
and $g(v_{n}) = \al v_{1}$. Since $\det(g) = -\al$, we conclude that $G$ contains all matrices 
that lie in $\tG$ and have determinant a power of $-\al$. Now we can follow the proof of Lemma 
\ref{root2} to obtain that $\deg(\Phi(h)) = 2^{a}$.

(ii) Now we will assume that 
$q \equiv 3 (\mod 4)$. Since $(n,q) \neq (2,3)$, we may apply Corollary \ref{two}(i) and get 
$\deg(\Phi(g)) \geq 2$. In particular we are done if $a = 1$. We will therefore assume that 
$a \geq 2$ and so $g^{2} \notin Z(G)$. 
By Lemma \ref{igl}(iv), $g^{2}$ belongs to a parabolic subgroup of 
$G$. Applying the main theorem of \cite{DZ1}, we see that 
$\deg(\Phi(g^{2})) = o(g^{2}) = 2^{a-1}$. Let $\mu$ be the unique eigenvalue of 
$\Phi(g^{2^{a}})$. Then $\Spec(\Phi(g^{2}))$ consists of all $2^{a-1}$-roots of $\mu$ in $\FF$. 
It suffices to show that $\deg(\Phi(g)) > \deg(\Phi(g^{2}))$.
  
Assume $b = 1$. If $Z(\tG)G = \tG$, then the statement follows from
Lemma \ref{sl22}. Otherwise $SL_{2}(q) \leq G \leq Z(\tG) \cdot SL_{2}(q)$ and $g^{q+1} = 1$,
and we can apply Lemma \ref{sl21}.

From now on we may assume that $b \geq 2$. Setting 
$2^{c} := (q+1)_{2} \geq 4$ and arguing as in (i), we see that $|\lam| = 2^{b+c}$, where $\lam$ 
is any eigenvalue of $g$ on $\overline{\FF}_{q}^{n}$. First we notice that $Z(\tG)G = \tG$. 
(Indeed, $G \geq [\tG,\tG]$ and the $2'$-part of $\tG/[\tG,\tG]$ is covered by $Z(\tG)$. 
On the other hand, $\det(g) = \lam^{(q^{n}-1)/(q-1)}$ generates the $2$-part of 
$\FF_{q}^{\bullet}$, whence the 
$2$-part of $\tG/[\tG,\tG]$ is covered by $g$.) In particular, $g$ and $g^{q}$ are 
$G$-conjugate by Lemma \ref{icl}(ii). Observe that $o(g) = 2^{b+c-1}$
and $|g| = 2^{b+c}$.  First suppose that $\mu \neq 1$, i.e. $\mu = -1$. Then any 
$\beta \in \Spec(\Phi(g^{2}))$ is a $2^{b+c-2}$-root of $-1$. Obviously, $\Spec(\Phi(g))$ 
contains a square root $\gamma$ of $\beta$, and such a $\gamma$ is a primitive $2^{b+c}$-root 
of unity. As $g$ and $g^{q}$ are conjugate, $\gamma^{Q} \in \Spec(\Phi(g))$ for 
$Q := q^{2^{b-1}}$. Since $(Q-1)_{2} = 2^{b+c-1}$, $\gamma^{Q} = -\gamma$. It follows that   
$\Spec(\Phi(g))$ contains both square roots of each $\beta \in \Spec(\Phi(g^{2}))$, whence 
$\deg(\Phi(g)) = 2^{a}$. Finally, we consider the case $\mu = 1$. In this case we show that 
for at least one value $\beta_{0} \in \Spec(\Phi(g^{2}))$, $\Spec(\Phi(g))$ contains both 
square roots of $\beta_{0}$. Since $\Spec(\Phi(g^{2}))$ consists of all $2^{b+c-2}$-roots of 
$\mu = 1$ and $b,c \geq 2$, we may choose $\beta_{0} = -1$. Now $\Spec(\Phi(g))$ contains a 
square root $\gamma$ of $\beta_{0}$, as well as $\gamma^{q}$. Since $q \equiv 3 (\mod 4)$, 
$\gamma^{q} = -\gamma$, and so we are done.       
\end{proof}

\begin{propo}\label{sp1}
{\sl Let $G := Sp_{2n}(q)$ with $n > 1$ and let $g \in G$ be an irreducible
$p$-element. Assume that $1 < \deg (\Theta(g)) < o(g)$ for some $\Theta \in \IBRL(G)$
with $(\ell,q) = 1$. Then $p > 2$, $o(g) = |g| = (q^{n}+1)/(2,q+1)$, 
and one of the following holds: 

{\rm (i)} $q$ odd, $n$ is a $2$-power, $\deg(\Theta(g)) = o(g)-1$, and the Sylow $p$-subgroups 
of $G$ are cyclic.

{\rm (ii)} $(n,q,|g|) = (3,2,9)$. Furthermore, either $\ell \neq 2$ and 
$\dim(\Theta) = \deg(\Theta(g)) = 7$, or $\ell = 3$, $\dim(\Theta) = 21$ and 
$\deg(\Theta(g)) \geq 7$.

{\rm (iii)} $(n,q,|g|) = (2,2,5)$. Furthermore, $\ell = 3$, and 
$\dim(\Theta) = \deg(\Theta(g)) = 4$.}   
\end{propo}

\begin{proof} 
1) By Proposition \ref{icl}(iii), $g$ can be embedded in a subgroup $X \simeq SL_{2}(q^{n})$
of $G$, $p > 2$, and $g^{q^{n}+1} = 1$. As $Z(G)=Z(X)$, the values of $o(g)$ for $G$ and for $X$ 
coincide; furthermore, $o(g) = |g|$. Let $\Phi$ be a non-trivial irreducible constituent of 
$\Theta|_{X}$. Then $\deg(\Phi(g)) < o(g)$. By Lemma \ref{sl21}, $|g|=(q^{n}+1)/(2,q+1)$ and 
either $\deg(\Phi(g)) = o(g)-1$, or $q$ is even and $\deg(\Phi(g)) = o(g)-2$.
We claim that $p$ is a primitive prime divisor of $q^{2n}-1$ and so the Sylow $p$-subgroups of 
$G$ are cyclic, unless $(n,q,p) = (3,2,3)$. Assume the contrary. Then $p|n$ and $p|(q^{2n/p}-1)$ 
by Proposition \ref{icl}(iii). Since $p > 2$ also divides $q^{n}+1$, we see that 
$p|(q^{n/p}+1)$. This in turn implies $(q^{n}+1)/(q^{n/p}+1) \equiv p (\mod 2p^{2})$ as 
$p > 2$. Since $(q^{n}+1)/(2,q+1)$ is a $p$-power, we conclude that $(q^{n}+1)/(q^{n/p}+1) = p$.
Thus $q^{n}+1 \leq (q^{n/p}+1)^{2}$, which is possible only when $(n,q,p) = (3,2,3)$.     

2) Here we show that $n$ is a $2$-power if $(n,q) \neq (3,2)$. Indeed, the condition 
$q^{n} \neq 8$ implies by \cite{Zs} that there is a primitive prime divisor $r$ of $q^{2n}-1$.
Since $p^{a} = |g| = (q^{n}+1)/(2,q+1)$, we have $r = p$. Thus $p$ is the unique primitive prime 
divisor of $q^{2n}-1$. Now assume $n$ is divisible by an odd prime $s$. Then $(q^{n}+1)/(2,q+1)$ is 
divisible by $q_{1} := (q^{n/s}+1)/(2,q+1)$ with $q_{1} > 1$ and $(q_{1},p) = 1$ 
(as $(p,q^{2n/s}-1) = 1)$, a contradiction.
 
So we are done if $q$ is odd. We may now assume that $q$ is even, and let $\theta$ be the 
Brauer character of $\Theta$.

3) Assume $(n,q,p) \neq (3,2,3)$. Then $|g| = p$ by Lemma \ref{zgm}, and by \cite{Z3} we are
also done if $\ell = p$. Assume $\ell \neq p$. Suppose there is an irreducible constituent $\Phi$ 
of $\Theta|_{X}$ of degree $q^{n}-1$. Then by Lemma \ref{sl21}, there is a primitive $|g|$-root 
of unity $\eps \in \FF$ such that
$\Spec(\Phi(g)) = \{1, \eps,\eps^{2}, \ldots ,\eps^{q^{n}}\} \setminus \{\eps,\eps^{-1}\}$. 
Then 
$\Spec(\Phi(g^{q})) = \{1, \eps,\eps^{2}, \ldots ,\eps^{q^{n}}\} \setminus \{\eps^{q},\eps^{-q}\}$.
However, $g$ and $g^{q}$ are $G$-conjugate by Proposition \ref{icl}(i). 
Hence $\Spec(\Theta(g))$ contains all $\eps^{j}$ and $\deg(\Theta(g)) = |g|$, a contradiction.
We have shown that any irreducible constituent $\Phi$ of $\Theta|_{X}$ is of degree $q^{n}$,
and $\Spec(\Theta(g)) = \{\eps,\eps^{2}, \ldots ,\eps^{q^{n}}\}$. In particular,
$q^{n}|\theta(1)$ and $\theta(g) = -\theta(1)/q^{n}$. 

The case $(n,q) = (2,2)$ leads to the conclusion (iii) by inspecting \cite{JLPW}. So we will 
assume $q^{n} > 4$. Claim that this assumption leads to a contradiction. Recall that $n > 1$ has
been shown to be a $2$-power. By embedding $SL_{2}(q^{n})$ in a subgroup $S \simeq Sp_{4}(Q)$ 
of $G$ with $Q := q^{n/2} > 2$, it suffices to prove the claim for $Sp_{4}(Q)$. 
Assuming $\deg(\Theta(g)) < |g|$, we have shown that there is $\varphi \in \IBRL(S)$ with
\begin{equation}\label{phi1}
  Q^{2}|\varphi(1)
\end{equation} 
and
\begin{equation}\label{phi2}
  \varphi(g) = -\varphi(1)/Q^{2}.
\end{equation}    
Suppose $\varphi$ lifts to characteristic $0$. Then (\ref{phi1}) implies by \cite{Eno} that
$\varphi$ is the (reduction modulo $\ell$) of the Steinberg character of $S$, but in this case
$\varphi(1) = Q^{4}$ and $\varphi(g) = 1$ violating (\ref{phi2}). Thus $\varphi$ does not
lift to characteristic $0$. This conclusion implies by \cite{Wh} that $\ell|(Q+1)$ and 
$\varphi$ is either of degree $Q(Q^{2}+1)/2 -1$, or 
$(Q-1)^{2}(Q^{2}+1) - (\al-2)Q(Q-1)^{2}/2$ with $1 \leq \al \leq Q/2$. In fact, 
$\al = 1$ or $2$ by \cite{OW1}. Since $Q > 2$, none of these two degrees is divisible by $Q^{2}$,
contradicting (\ref{phi1}).    

4) Assume $(n,q,p) = (3,2,3)$. Then $|g| = 9$, and $g$ is rational. First assume $\ell \neq 3$ 
and let $\eta$ be a primitive $9$-root of unity. Then the multiplicity of $1$, resp. $\eta^{3}$, 
as an eigenvalue of $\Theta(g)$ is equal to $a := (\theta(1) + 2\theta(g^{3}) + 6\theta(g))/9$, resp.  
$b := (\theta(1) + 2\theta(g^{3}) - 3\theta(g))/9$. Inspecting the Brauer characters of 
$G$ \cite{Atlas, JLPW}, we see that $a, b > 0$, except for the case $\dim(\Theta) = 7$, for
which $a = 1$, $b = 0$, $\deg(\Theta(g)) = 7$. 

Finally, we assume that $(n,q,p) = (3,2,3)$ and $\ell = 3$. Since $G$ has a unique class 
of elements of order $9$, we can embed $g$ in a subgroup $SL_{2}(8)$ of $G$ and obtain 
$\deg(\Theta(g)) \geq 7$ by Lemma \ref{sl21}. To find $\dim(\Theta)$, next we
embed $g$ in a subgroup $Y \simeq SU_{4}(2)$ and restrict $\Theta$ to $Y$. Direct 
computation done by F. L\"ubeck shows that the condition $\deg(\Theta(g)) \leq 8$ implies 
$\theta|_{Y} = a\varphi_{1} +  b\varphi_{5} + c\varphi_{10}$ for some integers $a,b,c \geq 0$ 
and $\varphi_{i} \in \IBR_{3}(Y)$ of degree $i = 1,5,10$. Let $h$ be an element of class $4A$ 
in $Y$ (in the notation of \cite{Atlas}). Restricting the unique complex irreducible character
$\chi$ of degree $7$ of $G$ to $Y$ we see that $h$ also belongs to the class $4A$ in $G$. Now
$\theta(h) = a+b+2c \geq \theta(1)/5$. Inspecting the $3$-Brauer characters of $G$ 
\cite{JLPW}, we now see that $\theta(1) = 7$, $21$, or $35$. Assume that $\theta(1) = 35$. 
Then $a = \theta(z) = 0$, whence $b+2c = \theta(h) = 7$. This in turn implies 
$\theta(h^{2}) = -3b+2c = -21$, $-13$, or $-5$, contradicting \cite{JLPW}. Thus
$\theta(1) = 7$ or $21$ as stated. 
\end{proof}

\begin{remar}
{\em The example of the $7$-dimensional representation $\Theta$ in Proposition \ref{sp1}(ii)
shows that there is no analogue of Lemma \ref{root2}(ii) for irreducible $p$-elements in 
$Sp_{2p}(q)$.}
\end{remar}

\begin{corol}\label{sp-gl}
{\sl Let $SL_{2n}(q) \lhd G \leq GL_{2n}(q)$ with $n > 1$ and let $g \in G$ be an irreducible
$p$-element with $p > 2$. Assume that $1 < \deg (\Theta(g)) < o(g)$ for some 
$\Theta \in \IBRL(G)$ with $(\ell,q) = 1$. Then $o(g) = |g| = (q^{n}+1)/(2,q+1)$, and one of
the following holds:

{\rm (i)} $n$ is a $2$-power, $\deg(\Theta(g)) = o(g)-1$, and the Sylow $p$-subgroups of $G$ are 
cyclic.

{\rm (ii)} $(n,q,|g|) = (3,2,9)$ and $\deg(\Theta(g)) \geq 7$.}
\end{corol}

\begin{proof}
Consider any eigenvalue $\lam$ of $g$ on $\FF_{q}^{2n}$. By Lemma \ref{igl}(ii), $|g| = |\lam|$ 
divides $q^{2n}-1$ but not $q^{n}-1$. If $p|(q^{n}-1)$ then, since $p > 2$, $p$ does not divide 
$q^{n}+1$ and so $p^{a} = |g|$ divides $q^{n}-1$, a contradiction. Thus $(p,q^{n}-1) = 1$ and 
so $|\lam|$ divides $q^{n}+1$. It is easy to see that in this case $p(t)$ is a multiple of 
$t^{2n}p(t^{-1})$, where $p(t)$ denotes the minimal polynomial of the matrix $g$. By \cite{Wa},
it follows that $g$ can be embedded in a subgroup $S \simeq Sp_{2n}(q)$ of $G$. Notice that the 
values of $o(g)$ are the same for $G$ and for $S$. Consider any nontrivial irreducible 
constituent $\Phi$ of $\Theta|_{S}$. Applying Proposition \ref{sp1} to $\Phi$ and excluding
the cases $(n,q,|g|) = (3,2,9)$ or $(2,2,5)$, we obtain all the conclusions of Corollary 
\ref{sp-gl}(i), except possibly the last one. Next, if $(n,q,p,|g|) \neq (3,2,3,9)$ then p. 1) of the 
proof of Proposition \ref{sp1} shows that $p$ is a primitive prime divisor of $q^{2n}-1$, whence 
the Sylow $p$-subgroups of $G$ are cyclic. It remains to consider the two exceptions.
Assume $(n,q,|g|) = (2,2,5)$. Then $\ell = 3$ by Proposition \ref{sp1}(iii), and we easily
obtain a contradiction by inspecting \cite{JLPW}. Consider the other exception $(n,q,|g|) = (3,2,9)$.
Note that any element of order $9$ in $G$ is irreducible. Next, $G$ contains a subgroup 
$X \simeq Sp_{6}(2)$, and $G$ has a unique conjugacy class of elements of order $9$. So we may
assume $g \in X$. Applying Proposition \ref{sp1} we obtain $\deg(\Theta(g)) \geq 7$. 
\end{proof}

\begin{lemma}\label{so1}
{\sl Let $V = \FF_{q}^{4n}$ be an orthogonal space of type $-$ with $n > 1$, 
$Spin(V) \leq G \leq \Gamma(V)$, and let $g \in G$ be an 
irreducible $p$-element. Assume that $1 < \deg (\Theta(g)) < o(g)$ for some 
$\Theta \in \IBRL(G)$ with $(\ell,q) = 1$. Then $p > 2$, $o(g) = |g| = (q^{2n}+1)/(2,q+1)$, 
$n$ is a $2$-power, $\deg(\Theta(g)) = o(g)-1$, and the Sylow $p$-subgroups of $G$ are cyclic.}   
\end{lemma}

\begin{proof} 
By Proposition \ref{icl}(iv), $g$ can be embedded in a subgroup $X \simeq SL_{2}(q^{2n})$
of $G$, $p > 2$, and $g^{q^{2n}+1} = 1$. As $Z(X)$ and $Z(G)$ are $2$-groups, $o(g) = |g|$
for both $X$ and $G$. Let $\Phi$ be a non-trivial 
irreducible constituent of $\Theta|_{X}$. Then $\deg(\Phi(g)) < o(g)$. By Lemma \ref{sl21}, 
$|g|=(q^{2n}+1)/(2,q+1)$ and either $\deg(\Phi(g)) = o(g)-1$, or $q$ is even and 
$\deg(\Phi(g)) = o(g)-2$. Considering a primitive prime divisor of $q^{4n}-1$, we see that in fact
$p$ is the unique primitive prime divisor of $q^{4n}-1$ and so the Sylow $p$-subgroups of $G$ are 
cyclic.

Arguing as in p.2) of the proof of Proposition \ref{sp1}, we see that $n$ is a $2$-power. 
So we are done if $q$ is odd. We may now assume that $q$ is even. Then $|g| = p$ by Lemma 
\ref{zgm}, and by \cite{Z3} we are also done if $\ell = p$. Assume $\ell \neq p$. Suppose there 
is an irreducible constituent $\Phi$ of $\Theta|_{X}$ of degree $q^{2n}-1$. Then by Lemma 
\ref{sl21}, there is a primitive $|g|$-root of unity $\eps \in \FF$ such that
$\Spec(\Phi(g)) = \{1, \eps,\eps^{2}, \ldots ,\eps^{q^{2n}}\} \setminus \{\eps,\eps^{-1}\}$.
Set $\tG := GO^{-}_{4n}(q)$. Recall that by Proposition \ref{icl}(i) $g$ is $\tG$-conjugate to 
the $4n$ powers $g^{q^{i}}$, $0 \leq i \leq 4n-1$. Since $(\tG :G) \leq 2$, we conclude that
$g$ and $q^{q^{2}}$ are $G$-conjugate. Clearly,   
$\Spec(\Phi(g^{q^{2}})) = \{1, \eps,\eps^{2}, \ldots ,\eps^{q^{2n}}\} 
 \setminus \{\eps^{q^{2}},\eps^{-q^{2}}\}$. Thus $\Spec(\Theta(g))$ contains all $\eps^{j}$ and 
$\deg(\Theta(g)) = |g|$, a contradiction. We have shown that any irreducible constituent $\Phi$ of 
$\Theta|_{X}$ is of degree $q^{2n}$, and 
$\Spec(\Theta(g)) = \{\eps,\eps^{2}, \ldots ,\eps^{q^{2n}}\}$. 
\end{proof}

The main goal of this section is to prove the following theorem:

\begin{theor}\label{clas1}
{\sl Under the assumptions of Theorem \ref{main2}, assume that the element $g$ is irreducible. 
Then one of the following holds.

{\rm (i)} $p^{a} \geq \deg(\Theta(g)) > p^{a-1}(p-1)$. 

{\rm (ii)} $p > 2$, $\deg(\Theta(g)) = p^{a-1}(p-1)$ and Sylow $p$-subgroups of $G/Z(G)$ are 
cyclic. Furthermore, either $a = 1$, or $\ell \neq p$. If $a \geq 2$ then 
$S = PSL^{\eps}_{2m+1}(q)$ or $P\Omega^{-}_{4m+2}(q)$. 

{\rm (iii)} $o(g) = (q^{n}-1)/(q-1)$, $S = PSL_{n}(q)$, $n > 2$ a prime, 
$\Theta$ is a Weil representation of degree $o(g)-1$ or $o(g)-2$, and 
$\dim(\Theta) = \deg(\Theta(g))$. Furthermore, either $n > 2$ and $\ell = o(g) = p$, or 
$n = 2$ and $q$ is even.}
\end{theor}

\begin{proof}
1) Without loss we may assume that $\Theta$ is an $\FF G$-representation, with $\FF$ an 
algebraically closed field of characteristic $\ell$. By Lemma \ref{easy}, we may assume that
$g$ is a $p$-element. Denote $Z := Z(G)$ and $L = G^{(\infty)}$. For any $x \in G$, $\bar{x}$ 
denotes the coset $xZ$. Let $o(g) = p^{a}$, $g^{p^{a}} = z \in Z$, and let 
$\Theta(z) = \mu \cdot \Id$. 

First we handle the case $p = 2$. By Lemma \ref{icl} and Proposition \ref{igl}, in this case
$SL_{n}(q) \leq G \leq GL_{n}(q)$ and $n = 2^{b}$. If $\ell \neq 2$ then we are done by 
Lemma \ref{root3}(ii). Notice (see \cite{Atlas}) that ${\rm Inndiag}(S)/S$ 
is either elementary abelian of order $4$ or cyclic for any finite simple group of Lie type; 
moreover the first case can happen only when $S = P\Omega^{+}_{4n}(q)$ with $q$ odd and $n \geq 2$. 
Thus the assumptions of Corollary \ref{two} are satisfied, and so $d_{\Theta}(h) > 1$ for 
any $h \in G \setminus Z$, and $d_{\Theta}(g) \geq 2^{a-1}+1$ if $\ell = 2$.
(In fact, this argument works for any $2$-elements and $\ell = 2$.) From now on we may assume 
$p > 2$.

Recall we are assuming $g$ is irreducible (on the natural module $V$ for the classical group
corresponding to $S$). Let $k$ be the smallest positive integer such that $p|(q^{k}-1)$. 
By Proposition \ref{icl}, $S$ is not of types $PSU_{2m}(q)$, $P\Omega_{2m+1}(q)$, or
$P\Omega^{+}_{2m}(q)$. 

\smallskip
2) Assume $S = PSp_{2n}(q)$. Then $S \leq G/Z \leq PCSp_{2n}(q)$. By Proposition \ref{icl}(iii), 
$\bg \in S$. Thus we can pick an element $h \in L$ such that $g \in hZ$. Since $o(g) = o(h)$ and 
$d_{\Theta}(g) = d_{\Theta}(h)$, we may assume $g = h \in L = Sp_{2n}(q)$. Hence we are done 
by Proposition \ref{sp1}.

\smallskip
3) Assume $S = PSL_{n}(q)$ or $PSU_{n}(q^{1/2})$, respectively. By 
Lemma \ref{igl} and Proposition \ref{icl}, $n/k = p^{b}$ for some $b \in \ZZ$. 

First we suppose that $b = 0$, that is, $p$ is a primitive prime divisor of $q^{n}-1$.
By virtue of \cite{Z3} and Corollary \ref{sp-gl}, 
we arrive at (i), (ii), or (iii) if $\ell = p$. So we may assume $\ell \neq p$. 
Then $(p,q-1) = 1$, whence $\bg \in S$. As in 2), we may assume $g \in L$ and $L = SL_{n}(q)$ or 
$SU_{n}(q^{1/2})$. Assume in addition that $r|n$ for some odd prime $r$. By comparing the $p$-part in 
the group order, we may embed $g$ in a subgroup $R = SL_{r}(Q)$ or $SU_{r}(Q^{1/2})$ with 
$Q := q^{n/r}$. Then the values of $o(g)$ are the same in $L$ and in $R$. Applying Corollary 
\ref{slprime}, we are done. Next assume that $n > 2$ is a $2$-power. Then $S = PSL_{n}(q)$, and we 
are done by Corollary \ref{sp-gl}. Finally, let $n = 2$; in particular,
$S = PSL_{2}(q)$ and $g^{q+1} = 1$. In this case one can just apply Lemma \ref{sl21}.

Next we suppose that $b \geq 1$ but $n > p$. By Lemma \ref{igl}(iv) and Proposition 
\ref{icl}(ii), $g^{p} \notin Z(G)$ and $g^{p}$ belongs to a parabolic subgroup of $G$. Moreover,
no eigenvalue of $g^{p}$ can belong to $\FF_{q}$ by Lemma \ref{igl}(i). 
By the main result of \cite{DZ1}, $d_{\Theta}(g^{p}) = o(g^{p})= p^{a-1}$. On the other hand,
$d_{\Theta}(g^{p}) > 1$ by Corollary \ref{two}(i); in particular, $a > 1$. If $\ell = p$, then 
Lemma \ref{root1}(ii) implies 
$d_{\Theta}(g) \geq p(d_{\Theta}(g^{p})-1)+1 = p^{a}-p+1 > p^{a-1}(p-1)$. Assume $\ell \neq p$.
If $k = 1$, then $a = b$ and $d_{\Theta}(g) = p^{a}$ by Lemma \ref{root2}(ii). We
will now assume that $k > 1$. 
In this case, by Lemma \ref{igl}(iv) and Proposition \ref{icl}(ii) we may view $g$ as an 
irreducible element of $GL_{p}(q^{n/p})$, resp. $GU_{p}(q^{n/2p})$. Notice that 
$L$ contains a subgroup $L_{1}$ isomorphic to $SL_{p}(q^{n/p})$, resp. $SU_{p}(q^{n/2p})$. Then 
$L_{2} := \langle L_{1},g \rangle$ is contained in $GL_{p}(q^{n/p})$, resp. $GU_{p}(q^{n/2p})$.
By Lemma \ref{root2} applied to $L_{2}$, if $\Phi$ is any irreducible constituent of 
$\Theta|_{L_{2}}$ of degree $> 1$, then $\Phi(g)$ consists of all $p$-roots of the unique 
eigenvalue of $\Phi(g^{p})$. Lemma \ref{div} implies $a = b+c$ if 
$p^{c} := (q^{k}-1)_{p}$. Let $\Phi$ be any constituent of $\Theta|_{L_{2}}$ of degree $1$. 
Since $L_{1}$ is perfect, $\Phi$ is trivial on $L_{1}$. Observe that $g^{q^{n/p}-1} \in L_{1}$ and 
$(q^{n/p}-1)_{p} = p^{b+c-1} = p^{a-1}$. As $g$ is a $p$-element, we conclude that 
$g^{p^{a-1}} \in L_{1}$; in particular, if $\beta$ is the unique eigenvalue of $\Phi(g^{p})$ then
$\beta^{p^{a-2}} = 1$. We have shown that $\Spec(\Theta(g))$ contains all $p$-roots of
$\beta \in \Spec(\Theta(g^{p}))$, except possibly for the $\beta$'s with $\beta^{p^{a-2}} = 1$,
for which we can say only that $\Spec(\Theta(g))$ contains at least one $p$-root of
$\beta$. Consequently, 
$\deg(\Theta(g)) \geq p(p^{a-1}-p^{a-2}) + p^{a-2} = p^{a-2}(p^{2}-p+1)$.

Finally, assume $n = p$. If $\ell \neq p$ then $o(g) = p = d_{\Theta}(g)$ by Lemma \ref{root2}.
The case $p = \ell$ is handled by Theorem \ref{slsup}.  

\smallskip
4) Here we handle the case $S = P\Omega^{-}_{2n}(q)$ with $n > 3$. Then 
$S \leq G/Z \leq PGO^{-}_{2n}(q)$.  Notice that $o(g) = |g|$ here, so $\mu = 1$.
Recall that $p > 2$, but the index of $S$ in $PGO^{-}_{2n}(q)$ divides $4$. 
Hence $\bg \in S$ and as in 2) we may assume $g \in L = Spin^{-}_{2n}(q)$. If $2|n$
then we are done by Lemma \ref{so1}. If $b = 0$ and $n$ is odd then we can apply 
Corollary \ref{slprime}.

From now on we suppose $n$ is odd and $b \geq 1$. 
Proposition \ref{icl}(iv) shows that $g^{p}$ belongs to a 
parabolic subgroup of $L$. Applying the main result of \cite{DZ1} to a nontrivial irreducible 
constituent of $\Theta|_{L}$, we obtain that $\deg(\Theta(g^{p})) = p^{a-1}$. Now if $\ell = p$ then 
$\deg(\Theta(g)) \geq p(p^{a-1}-1)+1$, so we will assume that $\ell \neq p$. Also by Proposition 
\ref{icl}(iv) we may embed $g$ in a subgroup $Y$ of $L$, with $SU_{n}(q) \leq Y \leq GU_{n}(q)$. 
First we assume that $k > 2$. Then $o(g)$ is the same for $G$ and for $Y$, and 
$|G|_{p} = |Y|_{p}$. Restricting $\Theta$ to $Y$ and applying the results of p. 3) above, we are done.

Notice that $k \neq 1$ as $2n/k = p^{b}$. It remains therefore to consider the case 
$k = 2$ and $n = p^{b}$. Let $p^{c} := (q+1)_{p}$. Then $o(g) = p^{b}$ in 
$Y$ by Lemma \ref{root2}, and $a = b+c$. We have already known that $\Spec(\Theta(g^{p}))$ 
consists of all $p^{a-1}$-roots of unity in $\FF$. We will prove that 
$\deg(\Theta(g)) \geq p^{a}-p^{c}+p^{c-1}$, which completes the proof of Theorem \ref{clas1}. 

First we suppose $b \leq c$. Consider any $p^{a-1}$-root $\beta$ of unity in $\FF$ and let 
$E(\beta) := \{\delta \in \FF \mid \delta^{p} = \beta\} \cap \Spec(\Theta(g))$. Clearly, 
$E(\beta) \neq \emptyset$. Claim that $|E(\beta)| = p$ if $\beta^{p^{c-1}} \neq 1$. For consider
such a $\beta$ and an irreducible constituent $\Phi$ of $\Theta|_{Y}$ such that 
$\beta \in \Spec(\Phi(g^{p}))$. It is easy to see that $g^{p^{c}} \in X := SU_{n}(q)$. Hence the 
assumption $\beta^{p^{c-1}} \neq 1$ implies that $\Phi|_{X}$ is nontrivial; in particular,
$\dim(\Phi) > 1$. By Lemma \ref{root2}(ii), $|E(\beta) \cap \Spec(\Phi(g))| = p$, whence the claim
follows. Thus $\deg(\Theta(g)) \geq p(p^{a-1}-p^{c-1}) + p^{c-1} = p^{a}-p^{c}+p^{c-1}$.

Finally, we suppose $b > c$. Consider any $p^{c}$-root $\gamma$ of unity in $\FF$. 
Recall that $g^{p^{b}} \in Z(Y)$. The condition $\deg(\Theta(g^{p})) = o(g^{p})$ implies that 
$\Theta|_{Y}$ contains an irreducible constituent $\Phi_{\gamma}$, such that 
$\Phi_{\gamma}(g^{p^{b}}) = \gamma \cdot \Id$ if $\gamma \neq 1$, and $\Spec(\Phi_{\gamma}(g^{p}))$ 
contains some $\beta$ with $1 = \beta^{p^{b-1}} \neq \beta^{p^{c-1}}$ if 
$\gamma = 1$. Consider the former case. If $\beta \in \Spec(\Phi_{\gamma}(g^{p}))$, then 
$\beta^{p^{b-1}} = \gamma \neq 1$, whence $\beta^{p^{c-1}} \neq 1$ as $b > c$. Arguing as in the
previous paragraph, we see that $\dim(\Phi_{\gamma}) > 1$. Moreover, 
$\Phi_{\gamma}(g^{p^{b}}) = \gamma \cdot \Id \neq \Id$. By Lemma \ref{root2}(ii), 
$\Spec(\Phi_{\gamma}(g))$ consists of all $p^{b}$-roots of $\gamma$. In the latter case, 
the condition $\beta^{p^{c-1}} \neq 1$ again implies that $\dim(\Phi_{1}) > 1$. So by Lemma 
\ref{root2}(ii), $\Spec(\Phi_{1}(g))$ consists of all $p^{b}$-roots of 
$\gamma=1$. Obviously, the sets $\Spec(\Phi_{\gamma}(g))$ for distinct $\gamma$ 
are disjoint. Consequently, $\deg(\Theta(g)) = p^{b+c} = p^{a}$.               
\end{proof}

The proof of Theorem \ref{clas1} yields the following consequence:

\begin{corol}\label{p-bound}
{\sl Under the assumptions of Theorem \ref{main2}, suppose that $p > 2$ divides the dimension of 
the natural module $V$ for $G$ and that $g$ is irreducible on $V$. Then one of the following holds.

{\rm (i)} $\deg(\Theta(g)) = o(g) = p$.

{\rm (ii)} $o(g) = p^{a} > p$ and $\deg(\Theta(g)) \geq p^{a-2}(p^{2}-p+1)$.
\hfill $\Box$}
\end{corol}

\section{Classical groups: Reducible elements}

\begin{propo}\label{su31}    
{\sl Let $S := SU_{3}(q) \leq G \leq H := GU_{3}(q)$ and $g \in G$ be a semisimple $p$-element with 
$o(g) = p^{a}$. Suppose $\Phi \in \IBRL(G)$ with $(\ell,q) = 1$ and $\dim(\Phi) > 1$. If 
$q = 2$, assume in addition that either $g \notin O_{3}(S)$ or $\Ker(\Phi)$ does not contain 
$O_{3}(S)$. Then either $\deg(\Phi(g)) = o(g)$, or one of the following holds.

{\rm (i)} $q+1 = p^{a}$, $\Phi|_{S}$ is a Weil representation, $g$ is a pseudoreflection, 
and $\deg(\Phi(g)) = o(g)-1$.

{\rm (ii)} $p^{a} = q^{2}-q+1$, $(3,q+1) = 1$, $\Phi|_{S}$ is the Weil representation of degree 
$q^{2}-q$, and $\deg(\Phi(g)) = o(g)-1$.

{\rm (iii)} $\ell = p > 2$, $o(g) = q+1$, $\Phi|_{S}$ is the Weil representation of degree 
$q^{2}-q$, and $g$ is represented by $\diag(\al^{i},\al^{i+1},\al^{i+k})$ for some 
$\al \in \FF_{q^{2}}^{\bullet}$ with $|\al| = q+1$ and some $i,k \in \ZZ$ such that 
$\al \neq \al^{k} \neq 1$. Furthermore, either $p = q+1$ is a Fermat prime and 
$\deg(\Phi(g)) \geq o(g)-1$, or $q = 8$ and $\deg(\Phi(g)) \geq 7$.

{\rm (iv)} $q = 3$, $o(g) = |g| = q+1$, $g$ is not a pseudoreflection, 
$\Phi|_{S}$ is the Weil representation of degree $q(q-1)$, and $\deg(\Phi(g)) = o(g)-1$.

{\rm (v)} $q = 2$, $|g| = 9$, $o(g) = 3$, $\Phi|_{S}$ is a Weil representation, and 
$\deg(\Phi(g)) = o(g)-1$.

{\rm (vi)} $2 = \ell = p|(q+1)$ and $g$ belongs to a maximal torus of
order $(q+1)^{3}$ of $GU_{3}(q)$.}
\end{propo}

\begin{proof}  
Let $V = \FF_{q^{2}}^{3}$ denote the natural module for $H := GU_{3}(q)$ and $\varphi$ be
the Brauer character of $\Phi$. Observe that Weil representations of $S$ extend to $H$.
We will assume $q > 2$ in pp. 1) -- 3) of this proof.

1) Suppose $g$ fixes a singular $1$-space of $V$. Applying the main result of \cite{DZ1} and 
\cite[Theorem 3.2]{GMST}, we arrive at (i). (More precisely, let $\chi = \zeta^{i}_{n}$ be the 
Weil character of $GU_{n}(q)$ as described in \cite[Lemma 4.1]{TZ2}, $0 \leq i \leq q$, and let 
$g$ be a pseudoreflection of order $q+1$. Then 
$$\chi|_{A} = \frac{q^{n-1}+(-1)^{n}}{q+1} \cdot \sum_{1_{A} \neq \lam \in \Irr(A)}\lam 
  - \left(1-\delta_{i,0}\right)(-1)^{n}\al_{i},$$
where $A := \la g \ra$ and $\al_{i}$ is the linear character of $A$ that sends $g$ to 
$\delta^{i}$, with $\delta = \exp(2\pi i/(q+1))$.) So we will assume that $g$ fixes no nonzero 
singular subspaces of $V$.

Suppose the $p$-Sylow subgroups of $G$ are cyclic; i.e. $2 \neq p |(q-1)$ or 
$3 \neq p |(q^{2}-q+1)$. Then $g \in S := SU_{3}(q)$. In the case 
$\ell = 0$ or $p$, applying the main result of \cite{Z3} to $S$ we arrive at (ii). Assume 
$\ell \neq 0,p$ and let $\Psi$ be an irreducible constituent of $\Phi|_{S}$. If $\Psi$ lifts to 
characteristic $0$ then we again arrive at (ii). In the remaining case, $\Psi$ extends to
$H$ by Example \ref{su30}. Since $H/S$ is cyclic, we see by \cite[Theorem III.2.14]{Fe}
that $\Phi = \Theta|_{G}$ for some $\Theta \in \IBRL(H)$, and $\Theta$ does not lift to 
characteristic $0$. According to Example \ref{su30}, we may assume $\varphi = \hat{\chi}-1$,
where $\chi$ is either the Steinberg character, or $\chi(1) = q(q^{2}-q+1)$. We will use the 
character table of $S$ as well as the notation for conjugacy classes of $S$ as given in \cite{Ge}.
In the case $2 \neq p|(q-1)$, any power $g^{a} \neq 1$ belongs to the class $C_{7}^{(k)}$ and so  
$\varphi(g^{a}) = -1$, whence $d_{\Phi}(g) = |g|$. In the case $3 \neq p|(q^{2}-q+1)$, any
power $g^{a} \neq 1$ belongs to the class $C_{8}^{(k)}$ and so $\varphi(g^{a}) = -2$, 
resp. $-1$, if $\chi = St$, resp. if $\chi(1) = q(q^{2}-q+1)$. Since $q \geq 3$, it follows
that $d_{\Phi}(g) = |g|$.

2) From now on we may assume $p|(q+1)$. Here we consider the case $g$ is reducible on $V$.
Under our assumptions, this implies by Proposition \ref{icl} that $g$ belongs to a torus 
$GU_{1}(q)^{3}$ of $H$. In particular, $C_{H}(g)G = H$, and $Z(G) \leq Z(H)$, so by Lemma 
\ref{easy} we may replace $G$ by $H$. 
Multiplying $g$ by a suitable central element of $H$, we may assume that $g$ is 
represented by $\diag(1,\al^{j},\al^{k})$ for some $\al \in \FF_{q^{2}}^{\bullet}$ of order
$q+1$ and some $j \neq k \in \ZZ$; in particular, $o(g) = |g|$. Since $o(g) = p^{a}$, we may assume 
that $|\al^{j}| = p^{a}$. Setting $h := \diag(1,\al^{j},1)$ and $t := \diag(1,1,\al^{k})$, we 
see that $g = ht$, $h$ is contained in a standard subgroup $X = GU_{2}(q)$ of $H$, and 
$o(h) = |h| = |g|$. By Example \ref{su30}, either $\varphi$ lifts to 
$\chi \in \Irr(H)$, or we may assume that $\varphi = \hat{\chi}-1$ with $\chi = St$ or 
$\chi(1) = q(q^{2}-q+1)$. 

Assume $\ell \neq p$ and $\chi$ is not a Weil character. By \cite[Remark 4.18]{TZ2}, $\chi|_{X}$ 
contains an irreducible constituent $\eta$ of degree $q+1$, whence 
$d_{\eta}(h) = o(h)$ by Lemma \ref{su2}. Restricting $\chi$ to $X \times \la t \ra$, we 
conclude that $d_{\chi}(g) = d_{\eta}(h) = |g|$; in particular, we are done if 
$\varphi = \hat{\chi}$. Assume $\varphi = \widehat{St}-1$. Then 
$3 \neq \ell|(q^{2}-q+1)$ and so $(\ell,q^{2}-1) = 1$. Hence $\hat{\eta}$ is irreducible over 
$X$, and so $\eta$ is a constituent of $\varphi|_{X}$, whence we again have  
$d_{\Phi}(g) = d_{\eta}(h) = |g|$. One can also check that 
$d_{\varphi}(g) = |g|$ in the other case where $\chi(1) = q(q^{2}-q+1)$. 

Assume $\ell = p > 2$ and $\chi$ is not a Weil character. Then again $\chi|_{X}$ contains an 
irreducible constituent $\eta$ of degree $q+1$. As $p > 2$, we can write $h = zh_{1}$ with 
$z \in Z(X)$ and $h_{1} \in X_{1} := SU_{2}(q)$. Since $2 < \ell|(q+1)$, all irreducible 
constituents of $\eta|_{X_{1}}$ are of $\ell$-defect $0$ and they are trivial at the 
$p$-element $z$; furthermore, $|h| = |h_{1}|$. It follows that $\hat{\eta}$ is irreducible over 
$X$ and it is an irreducible constituent of $\varphi|_{X}$; moreover, 
$\hat{\eta}|_{\la h_{1} \ra}$ is free and so $d_{\hat{\eta}}(h) = |h|$. Again restricting $\Phi$ 
to $X \times \la t \ra$, we conclude that $d_{\Phi}(g) = d_{\eta}(h) = |g|$. 

Assume $\chi$ is a Weil character; in particular, $\varphi = \hat{\chi}$. Multiplying 
$\chi$ by a linear character of $H$, we may assume that $\chi = \zeta^{i}_{3}$ where 
$0 \leq i \leq q$. In particular, $|\chi(x)| < q+1$ for all 
$1 \neq x \in A := \la g \ra$. First we consider the case $\ell \neq p$. If $|g| \leq (q+1)/2$, 
then the multiplicity of any linear character of $A$ in $\chi|_{A}$ is at least 
$(q^{2}-q - (q^{2}-1)/2)/|A| > 0$, whence $d_{\chi}(g) = |g|$. So we must have $|g| = q+1$. 
Assume $q \geq 7$. Then one can show that $|\chi(x)| < 4$ for at least $(q+1)/2 \geq 4$ elements 
$x \in A$. Hence the multiplicity of any linear character of $A$ in $\chi|_{A}$ is at least 
$(q^{2}-q - (q+1)(q-4) - 4 \cdot 4)/|A| > 0$, whence $d_{\chi}(g) = |g|$. In the remaining
cases $q = 3,4$ one can check directly that $d_{\chi}(g) = |g|$, except for the additional
case recorded in (iv), where $S = SU_{3}(3)$, $|g| = o(g) = 4$ but $g$ is not a pseudoreflection, 
and $d_{\chi}(g) = 3$. Since $d_{\Phi}(g) \geq 3$ by Corollary \ref{two} if $\ell = 2$,
we get $d_{\Phi}(g) = 3$ for any $\ell \neq 3$ in this case.    

Assume $\ell = p \neq 2$ and $\chi$ is a Weil character. By \cite{T}, the branching rule for 
restricting Weil characters to $X$ \cite{T} is as follows: 
$\zeta^{i}_{3}|_{X} = \sum_{0 \leq j \neq i \leq q}\zeta^{j}_{2}$. 
Hence $\zeta^{i}_{3}|_{X}$ always contains $\zeta^{j}_{2}$ for some $j \neq 0, (q+1)/2$. For such a 
$j$, $\widehat{\zeta^{j}_{2}}$ is irreducible over $X_{1} = SU_{2}(q)$. Applying Lemma \ref{su2} to 
the constituents of $\Phi|_{X}$ and the element $h$, we see that $\Phi|_{X}$ contains an 
irreducible constituent $\Psi$ such that either $d_{\Psi}(h) = |g|$, or $p^{a} = q+1$ and 
$d_{\Psi}(h) \geq |g|-2$. Restricting $\Phi$ to $X \times \la t \ra$, we conclude that 
either $d_{\Phi}(g) = |g|$, or $p^{a} = q+1$ and $d_{\Phi}(g) \geq |g|-2$. Consider the second 
possibility. Since $q+1 = \ell^{a}$, $\widehat{\zeta^{i}_{3}} = \widehat{\zeta^{0}_{3}} + 1$ if 
$i \neq 0$, cf. \cite[Theorem 7.2]{DT} and its proof. It follows that 
$\varphi = \widehat{\zeta^{0}_{3}}$ has degree $q^{2}-q$. Since $q+1 = p^{a}$ and $p > 2$, by 
Lemma \ref{zgm} we have either $a = 1$ or $p^{a} = 9$. In the former case, the rationality of 
$\zeta^{0}_{3}$ implies by Lemma \ref{unram} that $d_{\Phi}(g) \geq o(g)-1$. In the latter case 
where $p^{a} = 9$ we have $d_{\Phi}(g^{3}) = o(g^{3}) = 3$ (as $o(g^{3}) < q+1$) and so 
$d_{\Phi}(g) \geq 7$ by Lemma \ref{root1}. Thus we arrive at (iii).   
    
The open case $\ell = p = 2$ is recorded in (vi). 

3) Suppose that $g$ is irreducible and $p|(q+1)$. By Proposition \ref{icl}(ii), 
$p = 3 = o(g)$. In this case,  $d_{\Phi}(g) = 3$ by Theorem \ref{slsup}.

4) Finally, we consider the case $q = 2$, so $S = 3^{1+2}_{+} : Q_{8}$ and 
$H = 3^{1+2}_{+}:SL_{2}(3)$. Clearly, $|g| = 3$ or $9$, and $o(g) = 3$. If $|g| = 9$ then $g$ is 
irreducible by Proposition \ref{icl}(ii); moreover, $g \notin S$. Claim that if 
$x \in G \setminus Z(G)$ has order $3$ then either $x \in Q := O_{3}(S) = 3^{1+2}_{+}$ or $x$ is a 
pseudoreflection. Indeed, assume $|x| = 3$ but $x \notin S$. Notice that 
$GU_{1}(2)^{3}:\ZZ_{3}$ is a Sylow $3$-subgroup of $H$. If $x \in GU_{1}(2)^{3} \setminus S$ then 
$x$ is a pseudoreflection in $H$. Otherwise $x$ permutes $3$ orthonormal vectors 
$e_{1},e_{2},e_{3}$ in $V$ cyclically, in which case $x$ fixes the nonsingular vector 
$e_{1}+e_{2}+e_{3}$ and $x$ belongs to another subgroup $GU_{1}(2)^{3}$ of $H$ and so we are done 
again.         

Since $H$ is solvable, $\varphi$ lifts to $\chi \in \Irr(G)$. It suffices to consider two cases: 
$\ell = 0$ and $\ell = 3$. Suppose $\ell = 3$. Then $\Phi$ is an irreducible representation of 
$SL_{2}(3)$. By our assumption, $g \notin S$. Hence, $\deg(\Phi(g)) = 2$ if $\Phi$ has degree $2$ 
(in which case $\Phi|_{S}$ is a Weil representation and we arrive at (i) and (v)), and 
$\deg(\Phi(g)) = 3$ if $\dim(\Phi) = 3$.

From now on we may assume $\ell = 0$. Observe that $H$ has an irreducible character $\om$ of degree 
$3$ which is faithful on $O_{3}(S)$. We can list all irreducible characters of $H$ as follows:
$3$ of degree $1$ (we denote them by $1a$, $1b$, and $1c$), $3$ of degree $2$ (denoted by $2a$, 
$2b$, and $2c$), $7$ of degree $3$ ($\om \otimes 1abc$ and their complex conjugates, and one more 
denoted by $3d$ which is trivial on $O_{3}(S)$), $6$ of degree $6$ ($\om \otimes 2abc$ and their 
complex conjugates), $3$ of degree $8$ (denoted by $8a$, $8b$, and $8c$), and $2$ of degree $9$ 
($\om \otimes 3d$ and its complex conjugate). Among them, the seven characters $1abc$, $2abc$,
and $3d$ are trivial on $O_{3}(S)$. The ones of degree $2$ and the ones of degree $3$ except $3d$ 
are the $9$ Weil characters of $H$ -- they lead to the conclusions (i) and (v). It is easy
to check that $d_{\chi}(g) = o(g)$ if $\chi = 3d$. Assume $\chi = \om \otimes \gamma$ for some
$\gamma \in \Irr(SL_{2}(3))$ with $\gamma(1) > 1$. If $g \in O_{3}(S)$ then $\Spec(g,\om)$ consists of
all $3$ cubic roots of unity, so $d_{\chi}(g) = 3 = o(g)$. If $g \notin O_{3}(S)$, then 
$\Spec(g,\om)$ and $\Spec(g,\gamma)$ each contain at least $2$ (distinct) cubic roots of unity, whence 
$d_{\chi}(g) = 3 = o(g)$. It remains to consider the case $\chi(1) = 8$. Then 
$\chi|_{O_{3}(S)}$ is the sum of all $8$ nontrivial linear characters of $O_{3}(S)$, and so 
$d_{\chi}(g) = 3$ for $g \in O_{3}(S)$. On the other hand, if $g \notin O_{3}(S)$ then $g$ permutes
cyclically $3$ nontrivial linear characters of $O_{3}(S)$, yielding $d_{\chi}(g) = 3$.        
\end{proof}

\begin{corol}\label{su32}
{\sl Let $S := SU_{n}(q) \leq G \leq H := GU_{n}(q)$ with $n \geq 4$, and let $g \in G$ be a 
$p$-element with $o(g) = p^{a}$ and $g^{q+1} = 1$. Suppose $\Phi \in \IBRL(G)$ with $(\ell,q) = 1$ and 
$1 < \deg(\Phi(g)) < o(g)$. Assume in addition that $(\ell,p) \neq (2,2)$. Then $q+1 = p^{a}$, 
$g$ is a pseudoreflection, $\Phi|_{S}$ is a Weil representation, and 
$\deg(\Phi(g)) = o(g)-1$.}
\end{corol}

\begin{proof}
1) If $g$ belongs to a parabolic subgroup of $H$, then we are done by \cite{DZ1} and 
\cite[Theorem 3.2]{GMST}. So we may assume that $g$ fixes no nonzero singular subspace of 
the natural module $V = \FF_{q^{2}}^{n}$. Since $g^{q+1} = 1$, it follows that $g$ is contained 
in a torus $GU_{1}(q)^{n}$, i.e. $g$ can be diagonalized in an orthonormal basis of $V$. In 
particular, $C_{H}(g)G = H$, and $Z(G) \leq Z(H)$, so by Lemma \ref{easy} we may replace $G$ by 
$H$. Since $g$ does not fix any singular $1$-space of $V$, all diagonal entries of $g$ are 
distinct; in particular, $q+1 \geq n$ and so $q \geq 3$. Multiplying $g$ by a central element 
and reordering the diagonal entries of $g$, we may assume that 
$g = \diag(1,\al,\al^{j_{1}}, \ldots ,\al^{j_{n-2}})$ for some 
$\al \in \FF_{q^{2}}^{\bullet}$ with $|\al| = p^{a} > 2$; in particular, $o(g) = |g|$.
Setting $h := \diag(1,\al,\al^{j_{1}},1, \ldots ,1)$ and 
$t := \diag(1,1,1,\al^{j_{2}}, \ldots ,\al^{j_{n-2}})$, we see that $g = ht$, $h$ is contained in a 
standard subgroup $X = GU_{3}(q)$ of $H$, and $o(h) = |h| = |g|$. Moreover, 
$h$ is not a pseudoreflection in $X$. Recall we are assuming $\deg(\Phi(g)) < o(g)$. 
Restricting $\Phi$ to the subgroup $X \times \la t \ra$, we see that $\deg(\Psi(h)) < o(h)$ for 
every irreducible constituent $\Psi$ of $\Phi|_{X}$. Applying 
Proposition \ref{su31} to $h$ and any such constituent $\Psi$ of degree $> 1$, we conclude that 
$q+1 = p^{a}$ and $\Psi$ is a Weil representation. Moreover, either $\ell = p$ and 
$\deg(\Psi(h)) \geq o(g)-1$, or $\ell = 3$, $q = 8$ and $\deg(\Psi(h)) = 7$, or 
$q = 3$ and $\deg(\Psi(h)) = 3$. By \cite[Theorem 2.5]{GMST}, this 
in turn implies that $\Phi$ is a Weil representation, and either $\deg(\Phi(g)) \geq o(g)-1$, or
$q = 8$ and $\deg(\Phi(g)) = 7$.  

2) Consider the case $q \geq 4$. By Proposition \ref{su31}, $\ell = p$; moreover, either $q+1$ is a 
Fermat prime or $q = 8$. Write $g = h_{1}h_{2}$, with 
$h_{1} = \diag(1,\al,\al^{j_{1}}, \ldots ,\al^{j_{n-4}},1,1)$ and 
$h_{2} := \diag(1, \ldots ,1,\al^{j_{n-3}},\al^{j_{n-2}})$. Then we can embed $g$ into a 
standard subgroup $Y = Y_{1} \times Y_{2}$, with $h_{1} \in Y_{1} = GU_{n-2}(q)$ and 
$h_{2} \in Y_{2} = GU_{2}(q)$. Since $\Phi$ is a Weil representation, $\Phi|_{Y}$ contains an 
irreducible constituent $\Phi_{1} \otimes \Phi_{2}$, with $\Phi_{i}$ a Weil representation of $Y_{i}$. 

Assume $q \neq 8$ and $n \geq 5$. By the results of 1) (or Proposition \ref{su31}) 
applied to $Y_{1}$ and by Lemma \ref{su2} applied to $GU_{2}(q)$, 
$\deg(\Phi_{1}(h_{1})) \geq o(h_{1}) - 1 = o(g) -1$, and 
$\deg(\Phi_{2}(h_{2})) \geq 2$. By Lemma \ref{prod1}(ii), 
$\deg(\Phi_{1}(h_{1}) \otimes \Phi_{2}(h_{2})) = o(g)$, whence $d_{\Phi}(g) = o(g)$.

Assume $q = 8$. By the results of 1) (or Proposition \ref{su31}) 
and by Lemma \ref{su2} applied to $Y_{1}$ and $Y_{2}$, 
$\deg(\Phi_{1}(h_{1})) \geq o(h_{1}) - 2 = 7$, and 
$\deg(\Phi_{2}(h_{2})) \geq \min\{3,9-2\} = 3$ (notice that $o(h_{2}) = 3$ or $9$). Direct 
computation shows that $\deg(\Phi_{1}(h_{1}) \otimes \Phi_{2}(h_{2})) = 9$, whence 
$d_{\Phi}(g) = o(g)$. 

Assume $q \neq 8$ and $n = 4$. By Lemma \ref{su2} applied to $Y_{1}$ and $Y_{2}$, 
$\deg(\Phi_{1}(h_{1})) \geq o(h_{1}) - 2 = o(g) -2$, and 
$\deg(\Phi_{2}(h_{2})) \geq (q+1)-2 \geq 3$ (notice that $o(h_{2}) = p = q+1$ as $q+1$ is a 
prime). By Lemma \ref{prod1}(ii), $\deg(\Phi_{1}(h_{1}) \otimes \Phi_{2}(h_{2})) = o(g)$, whence 
$d_{\Phi}(g) = o(g)$.     

3) Finally, assume $q = 3$. Since $q+1 \geq n$, we get $n = 4$, and 
$g = \diag(1,\al,\al^{2},\al^{3})$ with $|\al| = 4$. Furthermore, $\ell \neq 2$ by our assumption.
Direct computation using \cite{Atlas} and \cite{JLPW} shows that $d_{\Phi}(g) = o(g) = 4$.
\end{proof}

\begin{lemma}\label{weil2}
{\sl Let $G := GU_{p^{b}+1}(q)$ with $q+1 = p^{c}$ for some odd prime $p$, and let $g \in G$ be a 
$p$-element such that $h := g^{p^{b}}$ is a pseudoreflection with $o(h) = q+1$. Let $V$ be  
an (irreducible) Weil module of $G$ in characteristic $\ell$ coprime to $q$. Then either 
$\dvg = p^{b}q$, or $(p^{b},q) = (3,2)$ and $\dvg = p^{b}q-1$.}
\end{lemma}

\begin{proof}
The case $(p^{b},q) = (3,2)$ follows by direct check, so we will assume $(p^{b},q) \neq (3,2)$.
It is well known (see e.g. \cite{DZ1}) that $d_{V}(h) = o(h)-1 = q$. Without loss we may assume 
that the natural $G$-module has an orthonormal basis $(e_{1}, \ldots ,e_{p^{b}+1})$ in which 
$h = \diag(\al, \ldots, \al,\beta)$ with $\al \neq \beta \in \FF_{q^{2}}$. It follows that in 
the same basis $g = \diag(g_{1},\gamma)$ with $\gamma \in \FF_{q^{2}}$ and 
$g_{1} \in H := GU_{p^{b}}(q)$. Replacing $g$ by $\gamma^{-1}g$ 
we may assume that $\beta = \gamma = 1$. Claim that $g_{1}$ is an irreducible $p$-element of $H$,
of order $p^{b}$ modulo $Z(H)$. (Indeed, since $o(h) = q+1$, $\alpha$ has order $q+1$. But 
$g^{p^{b}} = h$, hence $g_{1}$ is annihilated by $t^{p^{b}}-\alpha$. Now the proof of Lemma 
\ref{root2} establishes the claim.)

Consider the case $\ell \neq p$. Then we may assume that the Brauer character of $V$ 
is the reduction modulo $\ell$ of a Weil character $\zeta^{i}_{p^{b}+1,q}$ as given in 
\cite{TZ2} for some $i$, $0 \leq i \leq q$. The branching rule for Weil representations \cite{T} 
yields $\zeta^{i}_{p^{b}+1,q}|_{H} = \sum_{0 \leq j \leq q,~j \neq i}\zeta^{j}_{p^{b},q}$. 
Moreover, if $\delta$ is a primitive $(q+1)^{\mathrm {th}}$-root of unity in $\CC$ then
the central element $h$ acts as the scalar $\delta^{j}$ on the representation space $W_{j}$ of 
$\zeta^{j}_{p^{b}+1,q}$. By Lemma \ref{weil1}, the spectrum of $g$ on $W_{j}$ consists of all
$p^{b}$-roots of $\delta^{j}$. It follows that $\dvg = p^{b}q$ as stated. 

From now on we assume that $\ell = p$. Consider the filtration 
$V = V_{q} \supset V_{q-1} \supset \ldots \supset V_{1} \supset 0$, where $V_{k} := \Ker((h-1)^{k})$.
Also notice that $C_{G}(h) = H \times Z(G)$ and $Z(G)$ acts trivially on $V$. Furthermore, since 
$H$ is a standard subgroup of $G$, the branching rule for Weil representations 
\cite{T} implies that any $H$-composition factor of $V$, in particular of 
$V_{q}/V_{q-1}$, is either of degree $1$ or a Weil module. Claim that $V_{q}/V_{q-1}$ has 
at least one $H$-composition factor of degree $> 1$. Assume the contrary. Since $V_{q} \neq V_{q-1}$,
by Lemma \ref{filtr}(i) we see that the $H$-module $V$ has at least $q$ composition factors of
degree $1$. On the other hand, since $q+1 = p^{c}$ and $p > 2$ we may assume that the Brauer 
character $\zeta$ of $V$ is obtained by reducing the complex Weil character $\zeta^{1}_{p^{b}+1,q}$
modulo $p$. Again by \cite{T},
$\zeta^{1}_{p^{b}+1,q}|_{H} = \zeta^{0}_{p^{b},q} + \sum^{q}_{i=2}\zeta^{i}_{p^{b},q}$.
Moreover, $\zeta^{0}_{p^{b},q} (\mod p)$ is irreducible and $(\zeta^{i}_{p^{b},q} (\mod p)-1_{H})$ 
is irreducible for $i > 0$ by \cite{DT, HM}. Thus $V|_{H}$ contains exactly $q-1$ composition 
factors of degree $1$, a contradiction. 

Applying Lemma \ref{weil1} to the element $g_{1} \in H$ and an $H$-composition factor of 
degree $> 1$ (which is a Weil module) of $V_{q}/V_{q-1}$, we get    
$d_{V_{q}/V_{q-1}}(g) = d_{V_{q}/V_{q-1}}(g_{1}) = p^{b}$. Hence we are done by Lemma 
\ref{filtr}(ii). 
\end{proof}

\begin{propo}\label{su33}
{\sl Let $S := SU_{n}(q) \leq G \leq H := GU_{n}(q)$ with $n \geq 4$, and let $g \in G$ be a 
reducible $p$-element with $o(g) = p^{a}$ and $p|(q+1)$. Suppose $g$ does not belong to any 
parabolic subgroup of $H$, and $g^{q+1} \neq 1$. 

{\rm (i)} Then $p > 2$.

{\rm (ii)} Assume in addition that $\Phi \in \IBRL(G)$ with $(\ell,q) = 1$ and 
$1 < \deg(\Phi(g)) \leq p^{a-1}(p-1)$. Then $q+1 = p$ is a Fermat prime, $n = p^{b}+1$, 
$o(g) = p^{b+1}$, $g^{p^{b}}$ is a pseudoreflection in $H$, and $\Phi|_{S}$ is a Weil 
representation. Furthermore, either $\deg(\Phi(g)) = p^{b}(p-1)$, or $(p^{b},q) = (3,2)$ and 
$\deg(\Phi(g)) = p^{b}(p-1)-1$.} 
\end{propo}

\begin{proof}
1) Consider the action of $g$ on the natural module $V = \FF_{q^{2}}^{n}$ for $H$. Since $g$ does 
not belong to any parabolic subgroup, by Lemma \ref{decom} we can decompose $V$ into 
an orthogonal sum $V = \oplus^{s}_{i=1}V_{i}$ of nondegenerate subspaces such that $g$ is 
irreducible on each $V_{i}$. Let $n_{i} = \dim(V_{i})$, $n_{1} \geq \ldots \geq n_{s} \geq 1$. 
Since $g^{q+1} \neq 1$, $n_{1} > 1$, and $s > 2$ as $g$ is reducible. By Proposition \ref{icl}(ii), 
$n_{i} = p^{b_{i}}$ with $b_{1} \geq 1$ and $p > 2$. In particular, we are done with (i). Now we 
proceed to prove (ii). Let $g_{i} := g|_{V_{i}}$ and $p^{c} := (q+1)_{p}$. 
Notice that $C_{H}(g)G = H$. (Indeed, we have shown in the proof of Lemma 
\ref{root2} that $GU(V_{1}) = \la SU(V_{1}), Z(GU(V_{1})), g_{1} \ra$. Clearly, $Z(GU(V_{1}))$ 
and $g_{1}$ centralize $g$, and $H = \la SU(V),GU(V_{1}) \ra$.) So by Lemma \ref{easy} we may 
assume $G = H$. 

2) Denoting $b = b_{1}$, we choose largest $k$ such that $b_{k} = b$. Also, fix an eigenvalue 
$\lam_{i}$ for each $g_{i}$. Observe that $\lam_{i}^{p^{b}}$ with $1 \leq i \leq k$ are distinct;
in particular, $a > b$. Indeed, suppose $\lam_{1}^{p^{b}} = \lam_{2}^{p^{b}}$. By Lemma 
\ref{igl}(ii), the $p^{b}$ eigenvalues of $g_{1}$ on $V_{1}$ yield $p^{b}$ roots to the equation 
$x^{p^{b}} = \lam_{1}^{p^{b}}$. So without loss we may assume $\lam_{2} = \lam_{1}$. This in turn 
implies by Proposition \ref{icl}(i) that $g_{1}$ and $g_{2}$ are conjugate in $GU_{p^{b}}(q)$. 
Choosing bases in $V_{1}$ and $V_{2}$ suitably, we achieve the effect that the matrices of $g_{1}$ 
and $g_{2}$, as well as the Gram matrices, relative to these bases are the same. But in this case 
$g$ belongs to a parabolic subgroup of $G$ by Lemma \ref{equal}, a contradiction. 

3) Here we consider the case $\ell = p$. By Proposition \ref{icl}(ii), $g^{p}$ belongs to a 
parabolic subgroup of $G$, and $a \geq 2$ according to 2). By Lemma \ref{root1}(ii) and by
our assumption, $\deg(\Phi(g^{p})) \leq p^{a-2}(p-1)$, whereas $o(g^{p}) = p^{a-1}$. Hence
by \cite{DZ1} and \cite{GMST}, $q+1 = p$ is a Fermat prime, $g^{p^{e}} = (g^{p})^{p^{e-1}}$ is 
a pseudoreflection of order $q+1$ for some $e \geq 1$, and $V$ is a Weil module. Since $g$ does 
not belong to any parabolic subgroup of $G$, the arguments in \cite{GMST} imply that $e = b$, $k = 1$, 
$n = p^{b}+1$. Applying Lemma \ref{weil2}, we are done.
    
4) From now on we assume $\ell \neq p$. Here we consider the case $k = s$, i.e. 
$n_{1} = \ldots = n_{s} = p^{b}$. Then we can view
$g$ as the element $\diag(\lam_{1}, \ldots ,\lam_{s})$ in the subgroup $GU_{1}(Q)^{s}$ of 
$X := GU_{s}(Q)$ naturally embedded in $G$, where $Q := q^{n_{1}} \geq 8$. Observe that the 
values of $o(g)$ in $G$ and in $X$ are the same. For, assume $g^{p^{a-1}} \in Z(X)$. Then
$\lam_{1}^{p^{a-1}} = \ldots  = \lam_{s}^{p^{a-1}}$. Observe that $p^{a-1}(q+1)$ is divisible 
by $p^{b+c}$ (as $a > b$), and $|g| = |\lam_{i}| = p^{b+c}$ (as $p$ is odd). It follows that 
all eigenvalues of $g_{i}^{p^{a-1}}$ are equal to $\lam_{i}^{p^{a-1}}$, whence $g_{i}^{p^{a-1}}$
is scalar on $V_{i}$. Therefore, our assumption $\lam_{1}^{p^{a-1}} = \ldots  = \lam_{s}^{p^{a-1}}$ 
implies that $g^{p^{a-1}}$ is scalar on $V$, a contradiction.

Now we apply Lemma \ref{su2}, Proposition \ref{su31}, and Corollary \ref{su32} to the element $g$ 
in $X$ and an irreducible constituent $\Psi$ of degree $> 1$ of $\Phi|_{X}$. It follows that 
$d_{\Psi}(g) = o(g)$, except possibly when $Q+1 \in \{p^{a},2p^{a}\}$, and moreover $g$ is a 
pseudoreflection in $X$ if $s > 2$. By the conclusion of p. 2), $g$ cannot be a a pseudoreflection 
in $X$, so $s = 2$. Furthermore, if $Q > 8$ then by considering a primitive prime divisor $r$ of 
$q^{2n_{1}}-1$ we see that $r$ divides $Q+1$ but not $2(q+1)$, whence 
$Q+1 \notin \{p^{a},2p^{a}\}$. Thus $Q = 8$, $G = GU_{6}(2)$, and 
$\Spec(g,V) = \{\lam_{1}^{j} \mid j = 1,2,4,5,7,8\}$, i.e. $g$ belongs to the class $9C$ in 
$SU_{6}(2)$ (in the notation of \cite{Atlas}). In particular, $g$ is rational in $SU_{6}(2)$. 
Let $\phi$ be the Brauer character of any irreducible constituent of $\Phi|_{S}$. Assuming
$d_{\Phi}(g) \leq 7$, resp. $d_{\Phi}(g) = 8$, we obtain $\phi(1) + 2\phi(g^{3}) - 3\phi(g) = 0$,
resp. $\phi(1) + 2\phi(g^{3}) + 6\phi(g) = 0$. By inspecting \cite{Atlas} and \cite{JLPW},
we see that $\phi$ is trivial, whence $\dim(\Phi) = 1$, a contradiction. We conclude that
$d_{\Phi}(g) = o(g)$ in this case.    
   
5) Now we may assume that $k < s$, i.e. $b_{s} < b$. For $i \leq k$ we have  
$g_{i}^{p^{b+c-1}} = \om_{i} \cdot \Id$ is scalar, where $\om_{i} := \lam_{i}^{p^{b+c-1}}$ has 
order $p$. Meanwhile, $g_{i}^{p^{b+c-1}} = \Id$ for $i > k$. Hence $a = b+c \geq 2$
and $o(g) = |g|$. Following the proof of Lemma \ref{root2}(i), for each $i \leq k$ we can find 
$d_{i} \in GU(V_{i})$ such that $d_{i}^{p} = \Id$, $[d_{i},g_{i}] = \om_{i} \cdot \Id$. Setting 
$d_{i} = \Id$ for $i > k$ and $d = \diag(d_{1}, \ldots, d_{s})$, we see that $d^{p} = 1$, 
$[d,g] = g^{p^{a-1}}$. 

Decompose $\Phi = \oplus^{p-1}_{i=0}\Phi_{i}$, with $\Phi_{i}([d,g]) = \eps^{i} \cdot \Id$ for a 
primitive $p$-root $\eps$ of unity. Since $o([d,g]) = p$, $\Phi_{i} \neq 0$ for some $i \neq 0$. 
For any such $i$, Lemma \ref{root1}(i) applied to the actions of $\Phi_{i}(d)$ and $\Phi_{i}(g)$ 
yields that $\Spec(\Phi_{i}(g))$ consists of all $p$-roots of the elements in 
$\Spec(\Phi_{i}(g^{p}))$. Also observe that if $\beta \in \Spec(\Phi_{i}(g^{p}))$, then 
$\beta^{p^{a-2}} = 1$ if $i = 0$ and 
$\beta^{p^{a-2}} = \eps^{i} \neq 1 = \beta^{p^{a-1}}$ if $i \neq 0$, since
$[d,g] = g^{p^{a-1}}$. Similarly, if $\alpha \in \Spec(\Phi_{i}(g))$, then 
$\alpha^{p^{a-1}} = 1$ if $i = 0$ and $\alpha^{p^{a-1}} \neq 1 = \alpha^{p^{a}}$ if $i \neq 0$.

6) Assume $d_{\Phi}(g^{p}) = o(g^{p}) = p^{a-1}$. Then 
$d_{\Phi}([d,g]) = d_{\Phi}(g^{p^{a-1}}) = o(g^{p^{a-1}}) = p$, whence $\Phi_{j} \neq 0$ for all 
$j$. The argument in 5) shows that $\Spec(\Phi(g))$ contains all the $p$-roots of 
$\beta \in \Spec(\Phi(g^{p}))$ with $\beta^{p^{a-2}} \neq 1$, as well as at least one 
$p$-root of $\beta \in \Spec(\Phi(g^{p}))$ with $\beta^{p^{a-2}} = 1$. Thus 
$d_{\Phi}(g) \geq p(p^{a-1}-p^{a-2}) + p^{a-2} = p^{a-2}(p^{2}-p+1)$, contrary to our assumption. 

7) It remains to consider the case $d_{\Phi}(g^{p}) < p^{a-1}$. By Proposition \ref{icl}(ii), 
$g^{p}$ belongs to a parabolic subgroup of $G$; moreover, $\rank((g^{p}-z)|_{V_{1}}) \geq 4$ for 
any $z \in Z(G)$ if $b > 1$. So the inequality $d_{\Phi}(g^{p}) < p^{a-1}$ implies by the main 
result of \cite{DZ1} that either $g^{p^{e}} = (g^{p})^{p^{e-1}}$ is a pseudoreflection of order 
$q+1 = p^{a-e}$ for some $e \geq 1$, or $b = 1$, $(q,p^{a-1}) = (2,9)$ and $\rank(g^{p}-z) = 3$ 
for some $z \in Z(G)$. In the latter case, $c = 1$, $b = 1$, and $a = 3$, violating the equality 
$a = b+c$. So the former case must occur. As shown in \cite{GMST}, in this case $b = e$, 
$n = p^{b}+1$ (since $g$ does not belong to any parabolic subgroup), and $\Phi$ is a Weil 
representation. Hence we can apply Lemma \ref{weil2}. Assume furthermore that 
$(p^{b},q) \neq (3,2)$. Then $d_{\Phi}(g) = p^{b}q = p^{a} - p^{a-c}$. Since we are assuming 
$d_{\Phi}(g) \leq p^{a} - p^{a-1}$, we get $c = 1$, $p = q+1$, and so $p$ is a Fermat prime.
\end{proof}

\begin{remar}\label{csp}
{\em Stricly speaking, the main result of \cite{DZ1} in the symplectic case was proved only for 
$Sp_{2n}(q)$. However, one can also handle any subgroup $G$ such that 
$Sp_{2n}(q) \leq G \leq H := CSp_{2n}(q)$ and $n \geq 2$. Indeed, let $g \in G$ be a semisimple 
$p$-element with $o(g) = p^{a}$ and $g$ stabilizes a totally singular $m$-dimensional subspace $U$ 
of the natural module $\FF_{q}^{2n}$ for $G$. Let $\Theta \in \IBRL(G)$ be afforded by 
an $\FF G$-module of dimension $> 1$, where $(\ell,q) = 1$. If $p > 2$ then 
$g \in Z(H) \cdot Sp_{2n}(q)$, and so we are done by \cite{DZ1}. In the case $\ell = p = 2$, 
$d_{\Theta}(g) > 2^{a-1}$ by Corollary \ref{two}(ii). So we will assume $p = 2$ and $\ell \neq 2$. 
By Lemma \ref{easy} we can replace such a $G$ by $Z(H)G$, so it suffices to consider the case 
$G = H$ (and $q$ odd). Let $Q$ be the unipotent radical of $Stab_{G}(U)$ 
and $Z := Z(Q)$. One can check that $T := \la g \ra/\la g^{2^{a}}\ra$ acts faithfully on $Q$.  
Since $Z$ acts faithfully on $V$, $V_{1} := [Z,V] \neq 0$. Note that $T$ acts on the set $\XC$ of 
linear characters of $Z$ afforded by $V_{1}$. Let $V_{\lam}$ denote the $\lam$-eigenspace for $Z$ 
on $V$ when $\lam \in \XC$.

First assume $m = n$,  whence $Z = Q$. Claim that $T$ has an orbit of length $2^{a}$ on $\XC$,
which implies that $d_{\Theta}(g) = o(g)$. Otherwise, $h := g^{2^{a-1}}$ fixes every $\lam \in \XC$.
Since $\Theta(x)$ is scalar on each $V_{\lam}$ and on $C_{V}(Z)$ when $x \in Q$, it follows that 
$[Q,h] = 1$ on $V$, a contradiction.

Now we may assume that $1 \leq m \leq n-1$. Observe that $Q$ acts faithfully on $V_{1}$, and, 
since $T$ acts faithfully on $Q$, $QT$ acts faithfully (projectively) on $V_{1}$. Let $\OC_{i}$, 
$1 \leq i \leq s$, be the $T$-orbits on $\XC$. For each $i$, let 
$2^{a_{i}} := |\OC_{i}|$, $\lam_{i} \in \OC_{i}$, $\Phi_{i}$ be the representation of 
$\la g \ra Q$ on $V_{\lam_{i}}$, and $2^{b_{i}}$ be the order of $\Phi_{i}(g^{2^{a_{i}}})$ modulo 
$\FF^{\bullet}$. Since $\ell \neq 2$, $\Phi_{i}$ is a direct sum of irreducible summands 
$\Phi_{ij}$, $1 \leq j \leq t_{i}$, where $\Phi_{ij}$ is a faithful irreducible representation of
the extraspecial group $Q/\Ker(\lam_{i})$. We may assume $2^{b_{i}}$ is the order of 
$\Phi_{i1}(g^{2^{a_{i}}})$ modulo $\FF^{\bullet}$. Clearly, $a_{i} + b_{i} \leq a$. Moreover, 
$\deg(\Theta(g)) \geq 2^{a_{i}} \cdot \deg(\Phi_{i}(g^{2^{a_{i}}}))$ by Lemma \ref{perm}, so
we may assume that $a_{i} < a$. Observe that there is some $i$ such 
that $a_{i} + b_{i} = a$ (otherwise $[h,Q] = 1$ on $V_{1}$, a contradiction). For this $i$, 
suppose $\Phi_{i}(g^{2^{a}}) = \beta \cdot \Id$. Multiplying $\Theta(g)$ by $\gamma \in \FF$ 
with $\gamma^{2^{a}} = \beta^{-1}$, we may assume that $\Phi_{i}(g^{2^{a}}) = \Id$. Applying 
\cite[Theorem 2.6]{DZ1} to the representation $\Phi_{i1}$ of $\la g^{2^{a_{i}}},Q/\Ker(\lam_{i})\ra$,
we obtain that $\deg(\Phi_{i1}(g^{2^{a_{i}}})) \geq \max\{2,2^{b_{i}}-1\}$. Consequently,
$$\deg(\Theta(g)) \geq 2^{a_{i}} \cdot \deg(\Phi_{i1}(g^{2^{a_{i}}}))
   \geq 2^{a_{i}} \cdot \max\{2,2^{b_{i}}-1\} > 2^{a-1}.$$}    
\end{remar}

\begin{lemma}\label{sp2}
{\sl Let $Sp_{2n}(q) \lhd G \leq CSp_{2n}(q)$ with $n \geq 2$ and $q \equiv 3 (\mod 4)$. Assume 
$g \in G$ with $o(g) = 2^{a} \geq 2$, and $\Theta \in \IBRL(G)$ with $(\ell,2q) = 1$ and 
$\dim(\Theta) > 1$. Then $\deg (\Theta(g)) \geq 2^{a-1}+1$.}
\end{lemma}

\begin{proof}
1) By Corollary \ref{two}(i) we may assume that $a \geq 2$. Furthermore, by \cite{DZ1} and Remark
\ref{csp} we may assume that $g$ is a $2$-element that does not belong to any parabolic subgroup 
of $G$. Arguing as in the proof of Lemma \ref{decom}, we can decompose the natural module
$V = \FF_{q}^{2n}$ into a direct sum $\oplus^{s}_{i=1}V_{i}$ of irreducible $\la g \ra$-submodules.
By Proposition \ref{icl}(iii), $\dim(V_{i}) = 2$ for all $i$, and so $s = n$. Let $\lam_{i}$ be 
an eigenvalue of $g_{i} := g|_{V_{i}}$, $2^{a_{i}} := |\lam_{i}|$, $2^{c} := (q+1)_{2}$. Note
that $a_{i} \leq c+1$. We reorder the $V_{i}$'s such that $a_{1} \geq a_{2} \geq \ldots \geq a_{n}$. 
We also assume that $g$ changes the symplectic form on $V$ by the scalar 
$\tau \in \FF_{q}^{\bullet}$. In particular, $\tau = \det(g_{i}) = \lam_{i}^{q+1} = \pm 1$. 
Let $X := X_{1} \times \ldots \times X_{n}$, where $X_{i} := Sp(V_{i})$.

2) Here we consider the case $\tau = 1$, i.e. $g \in Sp_{2n}(q)$. Observe that 
$g_{i}^{2^{a_{i}-1}} = -1_{V_{i}}$, hence $a = a_{1}-1$ if $a_{1} = a_{n}$, and 
$a = a_{1}$ if $a_{1} > a_{n}$. Furthermore, $g_{i}^{q+1} = 1_{V_{i}}$ for all $i$ and 
$a_{i} \leq c$.

Assume $a = a_{1} -1$. Then $q \geq 7$ as $a \geq 2$. Since $g_{1} \notin \Ker(\Theta)$, we can 
find an irreducible constituent $\Phi = \Phi_{1} \otimes \ldots \otimes \Phi_{n}$ of $\Theta|_{X}$ 
with $\Phi_{i} \in \IBRL(X_{i})$ and $\dim(\Phi_{1}) > 1$. By Lemma \ref{sl21}, 
$\deg(\Phi_{1}(g_{1})) \geq 2^{a_{1}-1}-1 = 2^{a}-1$. It follows that 
$\deg(\Theta(g)) \geq \deg(\Phi(g)) \geq 2^{a}-1$.

Assume $a = a_{1} = \ldots = a_{k} > a_{k+1}$. Then $g_{i}^{2^{a-1}} = -1_{V_{i}}$
if $i \leq k$ and $g_{i}^{2^{a-1}} = 1_{V_{i}}$ for $i > k$. It follows that $h := g^{2^{a-1}}$ 
has $o(h) = 2$ and belongs to a parabolic subgroup of $G$. By \cite{DZ1}, $\deg(\Theta(h)) = 2$, 
and so $\Theta = \Theta_{+} \oplus \Theta_{-}$ (as $C_{G}(h)$-modules), with 
$\Theta_{\eps}(h) = \eps \cdot \Id$ for $\eps = \pm$. Clearly, $C_{G}(h) \geq X$. Consider any 
irreducible constituent $\Phi = \Phi_{1} \otimes \ldots \otimes \Phi_{n}$ of $\Theta_{-}|_{X}$ 
with $\Phi_{i} \in \IBRL(X_{i})$. Since 
$-\Id = \Phi(h) = \otimes^{k}_{i=1}\Phi_{i}(-1_{V_{i}}) \otimes \Id$, there is some $i \leq k$
such that $\Phi_{i}(-1_{V_{i}}) = -\Id$; in particular, $\dim(\Phi_{i}) \geq (q+1)/2$. If 
$q \geq 7$ then $\deg(\Phi(g)) \geq \deg(\Phi_{i}(g_{i})) \geq 2^{a-1}$ by Lemma \ref{sl21}. 
The same is true for $q = 3$, since $SL_{2}(3) = Q_{8}:3$ and 
$g_{i} \in Q_{8} \setminus Z(Q_{8})$. Moreover, if $\gamma \in \Spec(\Phi(g))$, then 
$\gamma^{2^{a-1}} = -1$ as $\Phi(h) = -\Id$. Obviously, $\delta^{2^{a-1}} = 1$ for any eigenvalue 
$\delta$ of $\Theta_{+}(g)$. It follows that $\deg(\Theta(g)) \geq 2^{a-1}+1$.   

3) From now on we assume that $\tau = -1$. In particular, $\lam_{i}^{q+1} = -1$ and so 
$|g_{i}| = 2^{c+1}$ for all $i$, and $a = c$. Here we consider the case $q \geq 7$ and set 
$Y := \la X,g \ra$, $Y_{i} := \la X_{i},g_{i} \ra$. Then we can find an irreducible constituent 
$\Phi$ of $\Theta|_{Y}$ such that the $\FF X$-representation $\Phi$ contains an irreducible 
constituent $\Phi_{1} \otimes \ldots \otimes \Phi_{n}$ with $\Phi_{i} \in \IBRL(X_{i})$ and 
$\dim(\Phi_{1}) > 1$. It is easy to see that $g_{i} \notin X_{i} \ni g_{i}^{2}$ and 
$o(g_{i}) = 2^{c}$ (as an element in any subgroup between $Sp(V_{i})$ and $CSp(V_{i})$), since 
$q \equiv 3 (\mod 4)$. By Lemma \ref{sl21}, either $\deg(\Phi_{1}(g_{1}^{2})) = 2^{c-1}$, or 
$\deg(\Phi_{1}(g_{1}^{2})) = 2^{c-1}-1$ and $(q+1)/2 = 2^{c-1}$. Since $q \geq 7$, $c \geq 3$ in 
the latter case. Thus in both cases $\deg(\Phi_{1}(g_{1}^{2})) > 2^{c-2}$.

Assume $\otimes^{n}_{i=1}\Phi_{i}$ is not $g$-stable. By Lemma \ref{perm}, 
$\deg(\Phi(g)) \geq 2\deg(\otimes^{n}_{i=1}\Phi_{i}(g_{i}^{2})) \geq 
 2\deg(\Phi_{1}(g_{1}^{2})) > 2^{c-1}$, whence $\deg(\Theta(g)) > 2^{a-1}$ as desired.

Assume $\otimes^{n}_{i=1}\Phi_{i}$ is $g$-stable. As $X = X_{1} \times \ldots \times X_{n}$, it 
follows that $\Phi_{i}^{g_{i}} \simeq \Phi_{i}$ for each $i$. Since $Y_{i}/X_{i} = \ZZ_{2}$ and 
$\ell \neq 2$, there are exactly two $\FF Y_{i}$-representations $\Phi^{\pm}_{i}$ with 
$\Phi^{\pm}_{i}|_{X_{i}} = \Phi_{i}$. Then $\Phi^{+} = \otimes^{n}_{i=1}\Phi^{+}_{i}$ is a
representation of (the external direct product) $Y_{1} \times \ldots \times Y_{n}$ with 
$\Phi^{+}|_{X} = \otimes^{n}_{i=1}\Phi_{i} = \Phi|_{X}$. Now 
$Y/X = \ZZ_{2}$, and $\Phi^{+}|_{Y}$ and $\Phi$ are two extensions of $\otimes^{n}_{i=1}\Phi_{i}$ 
to $Y$. Hence, multiplying $\Phi$ by an $\FF Y$-representation of degree $1$, we obtain 
$\Phi = \Phi^{+}|_{Y}$. It follows that 
$\deg(\Phi(g)) = \deg(\Phi^{+}(g)) \geq \deg(\Phi^{+}_{1}(g_{1}))$. Observe that 
$SL_{2}(q) \simeq Sp(V_{i}) < Y_{i} \leq CSp(V_{i}) \simeq GL_{2}(q) = 
 (\ZZ_{(q-1)/2} \times SL_{2}(q)) \cdot 2$. So we may view $\Phi^{+}_{1}$ as a representation of 
degree $>1$ of $GL_{2}(q)$. By Lemma \ref{sl22}, either $\deg(\Phi_{1}(g_{1})) \geq 2^{c}-1$, or 
$\deg(\Phi_{1}(g_{1})) = 2^{c}-2$ and $q+1 = 2^{c}$. Since $q \geq 7$, $c \geq 3$ in 
the latter case. Thus in both cases $\deg(\Phi_{1}(g_{1})) > 2^{c-1}$. 
Consequently, $\deg(\Theta(g)) > 2^{a-1}$ as desired. 

4) Finally, we consider the case $\tau = -1$ and $q = 3$. In this case, $\lam_{i}^{4} = -1$ 
and $\Spec(g_{i}) = \{\lam_{i},\lam_{i}^{3}\}$. Consider the element $h = \diag(g_{i},g_{j})$ of 
$CSp(U) \simeq CSp_{4}(3)$, where $U := V_{i} \oplus V_{j}$ with $i \neq j$. If 
$\Spec(g_{i}) = \Spec(g_{j})$, then $g$ is centralized by an element of order $3$, whence $g$ 
belongs to a parabolic subgroup, a contradiction. So 
$\Spec(g_{i}) \cap \Spec(g_{j}) = \emptyset$ for any $i \neq j$. It follows that $n = 2$ and 
$\Spec(g) = \{\lam_{1}^{j} \mid j = 1,3,5,7 \}$. In particular, $g$ belongs to class $4D$ in 
$PCSp_{4}(3)$ (in the notation of \cite{Atlas}) and $g$ is rational in $G = CSp(V)$. 
Since $o(g^{2}) = 2$, $\deg(\Theta(g^{2})) = 2$ by Corollary \ref{two}(i). Recall that 
$o(g) = 4$. If $\Theta(g^{4}) = \Id$, then $\Spec(\Theta(g^{2})) = \{1,-1\}$, and the 
rationality of $g$ implies that $\Spec(\Theta(g))$ contains both two nonreal $4$-roots 
$\pm \mu$ of unity and at least one of the two real $4$-roots of unity. If 
$\Theta(g^{4}) = -\Id$, then $\Spec(\Theta(g^{2})) = \{\mu,-\mu\}$, and the 
rationality of $g$ implies that $\Spec(\Theta(g))$ contains all four primitive $8$-roots of 
unity. Thus $\deg(\Theta(g)) \geq 3$ as desired.
\end{proof}

Next we work with the Clifford algebra $C(V)$ and Clifford group $\Gamma(V)$, where 
$V = \FF_{q}^{n}$ is endowed with a nondegenerate quadratic form for $q$ odd and $n \geq 2$.
We refer to \cite{TZ5} for basic notation and facts about them.

\begin{lemma}\label{so2}
{\sl Let $V = \FF_{q}^{n}$ be endowed with a nondegenerate quadratic form, where $n \geq 7$ and 
$q \equiv 3 (\mod 4)$. Assume that $Spin(V) \lhd G \leq H := \Gamma^{+}(V)$, 
$g \in G$ with $o(g) = 2^{a} \geq 2$, and that $\Theta \in \IBRL(G)$ with 
$(\ell,2q) = 1$ and $\dim(\Theta) > 1$. Then $\deg (\Theta(g)) \geq 2^{a-1}+1$.}
\end{lemma}

\begin{proof}
1) By Corollary \ref{two}(i) we may assume that $a \geq 2$. Furthermore, by \cite{DZ1} and Lemma 
\ref{easy} we may assume that $g$ is a $2$-element that does not belong to any parabolic subgroup
of $G$. By Lemma \ref{decom} we can decompose $V$ into a direct 
sum $\oplus^{s}_{i=1}V_{i}$ of irreducible $\la g \ra$-submodules. Let $\lam_{i}$ be an 
eigenvalue of $g_{i} := g|_{V_{i}}$, $2^{a_{i}} := |\lam_{i}|$, $2^{c} := (q+1)_{2}$. By 
Proposition \ref{icl}(iv), $n_{i} := \dim(V_{i})$ is $1$ (in which case $\lam_{i} = \pm 1$)  
or $2$ (in which case $V_{i}$ is of type $-$, $\lam_{i}^{q+1} = \det(g_{i}) = 1$). Thus 
$a_{i} \leq c$ for all $i$. We reorder the $V_{i}$'s such that 
$n_{1} \geq n_{2} \geq \ldots \geq n_{s}$, and $a_{i} \geq a_{j}$ whenever $n_{i} = n_{j}$. 
Denote by $\chi$ the canonical (surjective) map $H \to SO(V)$. Notice that 
$Z := \Ker(\chi) \simeq \ZZ_{q-1}$ and $Z \leq Z(H)$. In fact, $(Z(H):Z) = 2$ if $n$ is even:
if $z \in H$ satisfies $\chi(z) = -1_{V}$, then $zvz^{-1} = \chi(z)v = -v$ for any $v \in V$,
hence $C^{+}(V) = C_{C(V)}(z)$. Thus $Z(G) = Z(H) \cap G$. This in turn implies that 
$Z(G) \leq Z(Y)$ for any subgroup $Y$ between $G$ and $H$.  

2) Next we define some nondegenerate subspace $A$ of $V$. Observe that $n_{1} = 2$. Otherwise 
we can find distinct indices $i,j,k$ such that $n_{i} = n_{j} = n_{k} = 1$ and 
$g_{i} = g_{j} = g_{k} = \pm 1$, and so $g$ fixes a nonzero singular vector of the $3$-space 
$V_{i} \oplus V_{j} \oplus V_{k}$, a contradiction. By the same reason, the total number of 
$1$-dimensional $V_{i}$'s is at most $4$, whence $n_{1} = n_{2} = 2$ as $n \geq 7$. Moreover, if 
$n_{i} = n_{j} = 1$ and $g_{i} = g_{j} = \pm 1$ for $i \neq j$, then the $2$-space 
$V_{i} \oplus V_{j}$ must be of type $-$ (otherwise $g$ would fix a nonzero singular vector of 
$V_{i} \oplus V_{j}$). Now,

2a) if there are two indices $i,j$ such that $n_{i} = n_{j} = 1$ and $g_{i} = g_{j}$, then 
we define $A := V_{1} \oplus V_{2} \oplus V_{i} \oplus V_{j}$ and $k := j$;

2b) otherwise, define $A := V_{1} \oplus V_{2} \oplus V_{s}$ and $k := s$.\\ 
Thus $A$ is either a $5$-dimensional nondegenerate subspace, or a $6$-dimensional nondegenerate 
subspace of type $-$ of $V$. Moreover, since $\chi(g) \in SO(V)$, it is easy to check that 
the action of $g$ on $A$ belongs to $SO(A)$. 

Let $B := A^{\perp}$. For any $f \in K := \chi^{-1}(SO(A) \times SO(B)) \cap H$, let $\bfa$, 
resp. $\bfb$, denote the action of $f$ on $A$, resp. on $B$; in particular, 
$\chi(f) = \diag(\bfa,\bfb)$. By \cite[Lemma 6.2(i)]{TZ5} there are elements 
$\fa \in C^{+}(A) \cap \Gamma(V)$ and $\fb \in C^{+}(B) \cap \Gamma(V)$ such that 
$\chi(\fa) = \diag(\bfa,1_{B})$ and $\chi(\fb) = \diag(1_{A},\bfb)$. Then  
$\chi(\fa\fb) = \chi(f)$ and so $\fa\fb \in fZ$. Moreover, $\fa \in C^{+}(A)$. It follows 
by \cite[Lemma 6.1]{TZ5} that $[\fa,h_{B}] = 1$, if $h \in K$ and we define $h_{A}$, 
$h_{B}$ as described for $f$. Thus we have a decomposition of $K$ into a central product 
$K_{1} * K_{2}$, with $K_{1} \simeq \Gamma^{+}(A)$ and $K_{2} \simeq \Gamma^{+}(B)$. Notice that 
$g \in K$ by our choice of $A$.

3) Let $C := V_{1}$ and $D := C^{\perp}$. Define $f_{C}$ and $f_{D}$ for any 
$f \in \chi^{-1}(SO(C) \times SO(D)) \cap H$ as we did in 3) for the decomposition $V = A \oplus B$.
Then again $\chi(f_{C}f_{D}) = \chi(f)$ and so $f_{C}f_{D} \in fZ$, and $[f_{C},g_{D}] = 1$.
Choose $\bar{r}_{C}$ to be a generator of $SO(C) = \ZZ_{q+1}$ and take $r := r_{C}$. 
Since $\bar{g}_{C} \in \la \bar{r}_{C} \ra$, we get $\chi(g_{C}) = \chi(r)^{m}$ for some
$m \in \ZZ$ and so $g_{C} \in r^{m}Z$. Thus there is some $t \in Z$ such that 
$g = tr^{m}g_{D}$. Since $Z \leq Z(H)$ and $[r,g_{D}] = 1$, we conclude that $[r,g] = 1$.
Next observe that $\chi(r) \in SO(V) \setminus \Omega(V)$ and $G \lhd H$. By Lemma \ref{easy} 
we may replace $G$ by $\la G,r \ra$ and thereby assume $\chi(G) \geq SO(V)$. Next,  
$Z = \ZZ_{2} \times \ZZ_{(q-1)/2}$ as $q \equiv 3 (\mod 4)$, and $Z(G) \geq O_{2}(Z)$. So we 
may assume that $Z(G) \geq Z$ and $G = H$. It follows that $g \in K = K_{1} * K_{2} \leq G$, with
$K_{1} = \Gamma^{+}(A)$. 

4) Here we show that $o(g)$ (in $G$) and $o(\ga)$ (in $K_{1}$) are the same. Indeed, 
assume we are either in the case 2a), or in the case 2b) and $a_{1} > a_{k}$. Then
$g_{i}^{2^{a_{1}}} = 1_{V_{i}}$, $g_{1}^{2^{a_{1}-1}} = -1_{V_{1}}$, and 
$g_{k}^{2^{a_{1}-1}} = 1_{V_{k}}$. It follows that $\chi(g)^{2^{a_{1}}} = 1_{V}$ but 
$\chi(g)^{2^{a_{1}-1}}$ is not scalar on $V$, whence $g^{2^{a_{1}}} \in Z \leq Z(G)$ but 
$g^{2^{a_{1}-1}} \notin Z(G)$, i.e. $o(g) = 2^{a_{1}}$. The same argument applied to 
$\ga$ in $K_{1}$ yields $o(\ga) = 2^{a_{1}}$. Next assume we are in the case 2b) and 
$a_{1} = a_{k}$; in particular, $n_{s} = 2$ and $\dim(A) = 6$. In this case, 
$g_{1}^{2^{a_{1}-2}}$ is not scalar on $V_{1}$, and $g_{i}^{2^{a_{1}-1}} = -1_{V_{i}}$ for all
$i$. It follows that $\chi(g)^{2^{a_{1}-1}} = -1_{V}$ and so $g^{2^{a_{1}-1}} \in Z(G)$ 
as we mentioned in 1). On the other hand, $\chi(g)^{2^{a_{1}-2}}$ is not scalar on $V$, whence 
$g^{2^{a_{1}-2}} \notin Z(G)$, i.e. $o(g) = 2^{a_{1}-1}$. The same argument applied to 
$\ga$ in $K_{1}$ yields $o(\ga) = 2^{a_{1}-1}$.

5) Next we identify the subgroup $K_{1} = \Gamma^{+}(A)$. If $\dim(A) = 5$, then 
$\Gamma^{+}(A) \simeq CSp_{4}(q)$. Assume $\dim(A) = 6$. Claim 
that $K_{1} \simeq \ZZ_{(q-1)/2} \times L$, with $L \simeq SU_{4}(q) : 2$ and 
$L/Z(L) < PGU_{4}(q)$. Indeed, since $q \equiv 3 (\mod 4)$, we have 
$Z = Z_{1} \times Z_{2}$ with $Z_{1} := O_{2'}(Z) \simeq \ZZ_{(q-1)/2}$ and 
$Z_{2} := O_{2}(Z) \simeq \ZZ_{2}$. Next, $\Gamma^{+}(A) \rhd Spin(A) > Z_{2}$, where
$Spin(A) \simeq SU_{4}(q)$, $Spin(A)/Z_{2} = \Omega(A)$, and $Spin(A) \cap Z_{1} = 1$. If $Q$ 
denotes the quadratic form on $A$, we can find a pair of orthogonal vectors $u,v \in A$ 
such that $Q(u) = -1$ and $Q(v) = 1$. Then inside $C(A)$ we have 
$(uv)^{2} = u(-uv)v = 1$ and $\chi(uv) = \rho_{u}\rho_{v} \in SO(A) \setminus \Omega(A)$, if 
$\rho_{x}$ is the reflection corresponding to $x$. Setting $L := \la Spin(A),uv \ra$, we obtain
that $K_{1} = Z_{1} \times L$. It is well known that $\Out(\bar{S}) = D_{8}$ for 
$\bar{S} := P\Omega(A) \simeq PSU_{4}(q)$, and $D_{8}$ is induced by the action of the 
conformal orthogonal group $CO(A)$ on $\bar{S}$. Notice that 
$GO(A) = \la \Omega(A),\rho_{u},\rho_{v} \ra$ induces the subgroup $\ZZ_{2}^{2}$ of $\Out(\bar{S})$.
Furthermore, $CO(A)$ flips the involutions $\rho_{u}$, $\rho_{v}$ modulo $\Omega(A)$. 
It follows that $\rho_{u}\rho_{v}$ induces the central involution of $\Out(\bar{S})$. Since 
$-1_{A}$ lifts to a central element (of order $4$) of $\Gamma^{+}(A)$, $Z(L) \simeq \ZZ_{4}$. 
Finally, $GU_{4}(q)$ induces the subgroup $\ZZ_{4}$ and so contains the central involution of
$\Out(\bar{S})$. Consequently, $L/Z(L) < PGU_{4}(q)$. 

6) Now we can find an irreducible constituent $\Phi$ of $\Theta|_{K}$ such that 
$\Phi = \Phi_{1} \otimes \Phi_{2}$, $\Phi_{i} \in \IBRL(K_{i})$, and $\dim(\Phi_{1}) > 1$.
Recall that $g = \ga\gb z$ for some $z \in Z < K$, $\ga \in K_{1}$, and $\gb \in K_{2}$. Hence, 
$\deg(\Theta(g)) \geq \deg(\Phi(g)) \geq \deg(\Phi_{1}(\ga))$. Since $o(g) = o(\ga)$, it 
suffices to show that $\deg(\Phi_{1}(\ga)) >  o(\ga)/2$. If $\dim(A) = 5$, then we are done 
by Lemma \ref{sp2}. From now on we assume $\dim(A) = 6$. By Lemma \ref{easy} we may also 
assume that $\ga$ is a $2$-element, whence $\ga \in L$ for the subgroup $L$ defined in 5).  

Recall that $L := \la S,uv \ra$ with $S := SU_{4}(q)$ and $|uv| = 2$. Since $S$ is perfect, all 
irreducible constituents of $\Phi_{1}|_{S}$ are of degree $> 1$. Claim that, if 
$SU_{4}(q) \leq X \leq GU_{4}(q)$, $\Sigma \in \IBRL(X)$, $\dim(\Sigma) > 1$, then  
\begin{equation}\label{sigma-y}
  \deg(\Sigma(y)) \geq \max\{2,2^{d}-1\}
\end{equation}
for any $y \in X$ with $o(y) = 2^{d} \geq 2$. Indeed, if $d = 2$ then the claim follows from 
Corollary \ref{two}(i). So we may assume that $d \geq 2$, and that $y$ is a $2$-element by 
Lemma \ref{easy}. If $y$ belongs to a parabolic subgroup of $S$, then the claim follows from 
\cite{DZ1}. If $y$ does not belong to a parabolic subgroup, then $y$ is reducible (as $y$ is a 
$2$-element) and $y^{q+1} = 1$ by Proposition \ref{su33}(i), whence the claim follows from 
Corollary \ref{su32}.   

Now if $\ga \in S$ then we are done by applying (\ref{sigma-y}) to $y := \ga$. Consider the case 
$\ga \notin S$ and $\Phi_{1}|_{S}$ is reducible. By Lemma \ref{perm} and by (\ref{sigma-y}) applied to
$y := \ga^{2}$, 
$\deg(\Phi_{1}(\ga)) \geq 2 \cdot \deg(\Phi_{1}(\ga^{2})) \geq 2 \cdot \max\{2,2^{a-1}-1\} > 
 2^{a-1}$, and we are done again.

Finally, we consider the case $\ga \notin S$ and $\Phi_{1}|_{S}$ is irreducible. Since 
$L/Z(L) < PGU_{4}(q)$, we can find an element $h \in GU_{4}(q)$ such that the actions of $\ga$ and
of $h$ on $S$ are the same. In particular, $o(\ga)$ (in $L$) and $o(h)$ (in $GU_{4}(q)$) 
are the same, as they are equal to the order of $\ga$ and of $h$ in $\Aut(S)$. We 
use $\Phi_{1}|_{S}$ to define an irreducible representation $\Psi$ of 
$\tilde{L} = \la S,h \ra$ as follows. First, $\Psi(x) = \Phi_{1}(x)$ for any $x \in S$. 
Next, for all $x \in S$, 
$$\Psi^{h}(x) = \Psi(hxh^{-1}) = \Phi_{1}(hxh^{-1}) = \Phi_{1}(\ga x\ga^{-1}) = 
  \Phi_{1}(x) = \Psi(x),$$
i.e. $\Psi$ is $\tilde{L}$-stable. Since $\tilde{L}/S$ is cyclic, we can extend $\Psi$ to an 
irreducible representation of $\tilde{L}$, which we also denote by $\Psi$. Now the actions on 
$\Phi_{1}(S)$ via conjugation by $\Phi_{1}(\ga)$ and by $\Psi(h)$ are the same. But 
$\Phi_{1}|_{S}$ is irreducible, so $\Psi(h) = \beta\Phi_{1}(\ga)$ for some $\beta \in \FF^{\bullet}$.
In particular, $\deg(\Phi_{1}(\ga)) = \deg(\Psi(h))$. Applying (\ref{sigma-y}) to 
$(X,y,\Sigma) = (\tilde{L},h,\Psi)$, we are done.
\end{proof}

\begin{corol}\label{so3}
{\sl Let $V = \FF_{q}^{n}$ be endowed with a nondegenerate quadratic form, where $2|n \geq 8$ and 
$q \equiv 3 (\mod 4)$. Assume that $Spin(V) \lhd G \leq H := \Gamma(V)$, 
$g \in G$ with $o(g) = 2^{a} \geq 2$, and that $\Theta \in \IBRL(G)$ with 
$(\ell,2q) = 1$ and $\dim(\Theta) > 2$. Then $\deg (\Theta(g)) \geq 2^{a-1}+1$.}
\end{corol}

\begin{proof}
1) We will use all the assumptions and notations made in p. 1) of the proof of Lemma \ref{so2}. 
If $\chi(g) \in SO(V)$ and $G \leq \Gamma^{+}(V)$ then we are done by Lemma \ref{so2}. Assume 
$\chi(g) \in SO(V)$ but $G \not\leq \Gamma^{+}(V)$ and set $M := \la Spin(V),g,Z(G) \ra$. Notice
that $\chi(Z(G)) < SO(V)$, so $M \leq \Gamma^{+}(V)$. If the values 
of $o(g)$ in $G$ and in $M$ are the same then we are again done by Lemma \ref{so2}. The remaining
possibility can only occur when $n_{1} = \ldots = n_{s} = 2$, $a_{1} = \ldots = a_{s}$, and 
$-1_{V}$ gives rise to an element $t \in Z(M) \setminus Z(G)$ with $t = g^{2^{a-1}}$ 
(in which case the value of $o(g)$ in $M$ is $2^{a-1}$). By Corollary \ref{two}, $\deg(\Theta(t)) = 2$, 
so we can decompose $\Theta$ into a direct sum of two $\FF M$-representations $\Theta_{i}$, with 
$\Theta_{i}(t) = \al_{i} \cdot \Id$, $i = 1,2$, and $\al_{1} \neq \al_{2}$. Clearly, there is an $i$ 
and an irreducible constituent $\Phi$ of $\Theta_{i}|_{M}$ of degree $> 1$ (otherwise 
$\Ker(\Theta) \geq Spin(V)$). But $M \lhd G$, so by Clifford's theorem there must be a $G$-conjugate
$\Psi$ of $\Phi$ that is contained in $\Theta_{3-i}|_{M}$. Applying Lemma \ref{so2} to 
$M \leq \Gamma^{+}(V)$, we see that $\deg(\Phi(g)) \geq 2^{a-2}+1$ and 
$\deg(\Psi(g)) \geq 2^{a-2}+1$. Furthermore, 
$\Spec(\Phi(g)) \cap \Spec(\Psi(g)) = \emptyset$ as $g^{2^{a-1}} = t$ and $\al_{1} \neq \al_{2}$. 
It follows that $\deg(\Theta(g)) \geq 2(2^{a-2}+1) > 2^{a-1}$ as stated.      

2) From now on we assume $\chi(g) \notin SO(V)$. Let $Z_{1} := O_{2'}(Z) \simeq \ZZ_{(q-1)/2}$, 
$Z_{2} := O_{2}(Z) \simeq \ZZ_{2}$, and let $\bj$ denote the unique central involution of 
$H$. Considering a pair of orthogonal vectors $u,v \in V$ such that $Q(u) = 1$, $Q(v) = -1$ (where $Q$ 
is the quadratic form on $V$), we get the subgroup $L := \la Spin(V),u,v \ra = Spin(V) \cdot D_{8}$, 
and it is easy to see that $H = Z_{1} \times L$. Replacing $G$ by $GZ_{1}$ by Lemma \ref{easy}, we get 
$G = Z_{1} \times (G \cap L)$.   

Let $k$ be the total number of $2$-dimensional $V_{i}$'s,
$l$ the total number of $i$ such that $n_{i} = 1$ and $g_{i} = -1$, and $m$ the total number of 
$i$ such that $n_{i} = 1$ and $g_{i} = 1$. Then $2k + l + m = n$ and $l$ is odd. If $l \geq 3$,
then $g$ fixes a nonzero singular vector in $\oplus^{k+3}_{i=k+1}V_{i}$, a contradiction. Similarly,
$m < 3$. It follows that $l = m = 1$. 

3) Consider the case where $\Theta(\bj) = \Id$. There is no loss to view $\Theta$ as an irreducible 
representation of $K := (G \cap L)/Z_{2} \leq GO(V)$ and replace $g$ by $\chi(g)$ (and assume 
$o(g) > 2$). Define $C := \oplus_{i \neq k+1}V_{i}$, $D := V_{k+1}$. Then 
$g = \diag(g',g_{k+1})$, with $g' \in SO(C)$ and $g_{k+1} = -1$. Since $n_{1} = 2$, one can
find an element $h \in SO(V_{1}) \setminus \Omega(V_{1})$ that centralizes $g$. Similarly,
$g_{k+1} \in GO(V_{k+1}) \setminus SO(V_{k+1})$ centralizes $g$. Therefore, by Lemma \ref{easy}, 
we may replace $K$ by $GO(V)$. It is easy to check that $o(g)$ (in $K$) and $o(g')$ (in $SO(C)$) are 
the same. Restricting $\Theta$ to the subgroup $SO(C) \times GO(D)$ of $K$ and applying Lemma 
\ref{so2} to $g'$, we conclude that $\deg(\Theta(g)) > o(g)/2$.  

4) Now we assume $\Theta(\bj) = -\Id$. Fixing an element $z \in \Gamma(V)$ with $\chi(z) = -1_{V}$, 
we observed in p. 1) of the proof of Lemma \ref{so2} that $z$ centralizes $C^{+}(V)$; in 
particular, $[z,gv^{-1}] = 1$. But $zvz^{-1} = \chi(z)v = \bj v$ (in $G$), hence $zgz^{-1} = \bj g$. 
It follows that $-\gamma \in \Spec(\Theta(g))$ whenever $\gamma \in \Spec(\Theta(g))$. Therefore, 
$\Spec(\Theta(g))$ contains all square roots of elements in $\Spec(\Theta(g^{2}))$. Clearly, 
$g^{2} \in G_{1} := \la x^{2} \mid x \in G \ra \leq \Gamma^{+}(V)$. Moreover, since 
$a_{1} \geq 2$ and $g_{k+2} = 1$, one can check that the values of $o(g^{2})$ in $G_{1}$ and in $G$ 
are the same. Applying Lemma \ref{so2} to an irreducible constituent of degree $> 1$ of 
$\Theta|_{G_{1}}$ and $g^{2}$, we see that $\deg(\Theta(g^{2})) \geq 2^{a-2}+1$.
Consequently, 
$\deg(\Theta(g)) \geq 2 \cdot \deg(\Theta(g^{2})) \geq 2(2^{a-2}+1) > 2^{a-1}$.
\end{proof}

\begin{theor}\label{clas2}
{\sl Under the assumptions of Theorem \ref{main2}, assume that the element $g$ is reducible. 
Then one of the following holds.

{\rm (i)} $p^{a} \geq \deg(\Theta(g)) > p^{a-1}(p-1)$. 

{\rm (ii)} $p > 2$, $\deg(\Theta(g)) = p^{a-1}(p-1)$ and Sylow $p$-subgroups of $G/Z(G)$ are cyclic.
Furthermore, either $a = 1$, or $\ell \neq p$.

{\rm (iii)} $S = PSU_{n}(q)$ with $n \equiv 1 (\mod p^{b})$ for some $b \geq 1$, $q+1 = p$ is a 
Fermat prime, $o(g) = p^{b+1}$, $g^{p^{b}}$ is a pseudoreflection in $GU_{n}(q)$, and 
$\Theta$ is a Weil representation. Furthermore, either $\deg(\Theta(g)) = p^{b}(p-1)$, or 
$(n,p^{b},q) = (4,3,2)$ and $\deg(\Theta(g)) = p^{b}(p-1)-1$. 

{\rm (iv)} $S = PSU_{n}(q)$, $o(g) = p = q+1$ is a Fermat prime, $g$ is contained  
in $GU_{1}(q)^{n}$, $\Theta$ is a Weil representation, and $\deg(\Theta(g)) = o(g)-1$. Furthermore, 
either $n \leq 3$ or $g$ is a pseudoreflection in $GU_{n}(q)$.}
\end{theor}

\begin{proof}
1) Let $V = \FF_{q}^{n}$ denote the natural module for $G$, 
and we write $S = PSU_{n}(q^{1/2})$ in the $SU$-case. We will 
use the notation of p. 1) of the proof of Theorem \ref{clas1}. As we mentioned there, assertion (i) 
holds if $\ell = p = 2$. So we may assume that $p^{a} > 2$, $\ell \neq 2$ 
if $p = 2$, $g$ is a $p$-element, and $d_{\Theta}(g) \leq p^{a-1}(p-1)$. 

Assume $g$ belongs to a parabolic subgroup of $G$. By the main result of \cite{DZ1} and Remark
\ref{csp}, $d_{\Theta}(g) > p^{a-1}(p-1)$, except for the case $(G,g)$ is as in (iii) or (iv), in 
which cases $\Theta|_{L}$ is a Weil representation by \cite[Theorem 3.2]{GMST}. 

Thus we may assume that $g$ is not contained in any parabolic subgroup of $G$, and that 
$S \neq PSL_{n}(q)$. In the remaining cases, $V$ is endowed with a $G$-invariant nondegenerate 
Hermitian, symplectic, or quadratic form, and $g$ cannot fix any nonzero totally singular subspace. 
By Lemma \ref{decom}, we can find an orthogonal decomposition $V = \oplus^{s}_{i=1}V_{i}$, where 
$g$ acts as an irreducible element $g_{i}$ on each (nondegenerate) subspace $V_{i}$, with one 
possible exception. In the exception, $S = P\Omega^{\eps}_{2m}(q)$ with $2|q$ and $m \geq 4$, and 
$g \in H = Sp_{2m-2}(q)$. When this happens, we will
\begin{equation}\label{repl}
   \mbox{ restrict to }H \mbox{ and refer to our result for the symplectic case.}
\end{equation} 
Throughout pp. 2) -- 5) of the proof we will assume that the aforementioned 
decomposition exists. Let $k$ be the smallest positive integer such that $p|(q^{k}-1)$ and let 
$n_{i} = \dim(V_{i})$. Observe that $g^{p^{a}}$ acts scalarly on $V$. We may assume that 
$n_{1} \geq \ldots \geq n_{s}$ and $s \geq 2$. 

\smallskip
2) Here we handle the case $k = 1$. First suppose that $S = PSU_{n}(q^{1/2})$. 
If $n_{1} = 1$, then we are done by Lemma \ref{su2} if $n = 2$, by Proposition \ref{su31} if 
$n = 3$, and by Corollary \ref{su32} if $n \geq 4$. On the other hand, if $n_{1} > 1$ then 
$g^{q^{1/2}+1} \neq 1$ and $p > 2$ by Proposition \ref{icl}(ii), and so we arrive either at 
(i) or (iii) by Proposition \ref{su33}.

Next we consider the case $k = 1$ for the remaining classical groups. By Proposition 
\ref{icl}, this implies that $p = 2$ and $q$ is odd. If $q \equiv 3 (\mod 4)$, then we are done
by Lemmas \ref{sp2}, \ref{so2} and Corollary \ref{so3}. Assume $q \equiv 1 (\mod 4)$. 
If $V = \FF_{q}^{n}$ is endowed with a quadratic form, then it is easy to see that 
$g_{i}^{2} = \Id$ for each $i$, whence $o(g) = 2$. Assume $V$ is endowed with a symplectic form 
$(\cdot,\cdot)$. Then for every $i$, $n_{i} = 2$, $|\lam_{i}| = 2(q-1)$ if 
$\lam_{i} \in \Spec(g_{i})$, and $g_{i}$ is conjugate to $\diag(\lam_{i},\lam_{i}^{q})$ in 
$GL_{2}(\overline{\FF}_{q})$. If $g$ changes the form $(\cdot,\cdot)$ by $\tau$, then 
$\tau = \det(g_{i}) = -\lam_{i}^{2}$. Hence $g^{2} = -\tau \cdot \Id$ is scalar on $V$ and 
$o(g) = 2$. Thus in all cases $o(g) = 2$ and we are done by Corollary \ref{two}(i).  
 
\smallskip
3) Henceforth we assume that $k > 1$; in particular, $p > 2$. The assumptions on $G,p$ now imply 
that each $g_{i}$ belongs to $SU(V_{i})$ in the unitary case, $Sp(V_{i})$ in the symplectic case, and 
$\Omega(V_{i})$ in the orthogonal case; moreover, $g^{p^{a}} = 1$. Choose 
$X_{i} := SU(V_{i})$, $Sp(V_{i})$, or $Spin(V_{i})$, respectively, and let $X := X_{1}$. Without
loss we may identify $g_{i}$ with an inverse image of $p$-power order of it in $X_{i}$. Since any 
central $p$-element of $I(V_{i})$ is trivial and since $n_{1} \geq \ldots \geq n_{s}$, we may assume 
that the order of $g_{1}$ modulo $Z(X)$ is $p^{a}$. 

4) Suppose Theorem \ref{clas1} is applicable to the element $g_{1}$ in the group $X$, and a 
nontrivial irreducible constituent $\Phi_{1}$ of $\Theta|_{X}$ but $X \neq SL_{2}(q)$ for even $q$; 
in particular, either $n_{1} \geq 3$, or $n_{1} = 2$ and $G$ is of symplectic type. Then
$\deg(\Phi_{1}(g_{1})) \leq \deg(\Theta(g_{1})) \leq p^{a-1}(p-1)$. By Theorem \ref{clas1} 
and Corollary \ref{p-bound}, $k = n_{1}$, $d_{\Phi}(g_{1}) = p^{a-1}(p-1)$, and either $a = 1$ or
$\ell \neq p$ (notice that $X/Z(X) \not\simeq PSL_{r}(q)$ for any odd prime $r$). 

We aim to show that {\it the Sylow $p$-subgroups of $G$ are cyclic} and so (ii) holds. Assume the 
contrary. This implies that $n \geq 2n_{1}$. If $n_{2} < n_{1}$, then $g$ acts trivially on the 
nondegenerate subspace $V_{1}^{\perp}$, which contains a nonzero singular vector as 
$\dim(V_{1}^{\perp}) \geq n_{1}$, contrary to our assumption on $g$. Thus $n_{2} = n_{1}$. We 
consider two subcases: $|g_{2}| = p^{a}$ and $|g_{2}| < p^{a}$.

\smallskip
4a) Subcase I: $|g_{2}| = p^{a}$. Assume in addition that 
\begin{equation}\label{cond}
  \mbox{Either \ref{prod4} or \ref{prod5} or \ref{prod6} applies to the subgroup }
  X_{1} * X_{2} \mbox{ of }G.
\end{equation} 
Then $\Phi_{1}$ can be chosen such that $\Phi_{1} \otimes \Phi_{2}$ is an irreducible 
constituent of $\Theta|_{X_{1} * X_{2}}$, $\Phi_{2} \in \IBRL(X_{2})$, and $\dim(\Phi_{2}) > 1$. 
In particular, $\deg(\Phi_{2}(g_{2})) \geq p^{a-1}(p-1)$ by Theorem \ref{clas1}. If $\ell \neq p$, 
then $\deg(\Theta(g)) \geq \deg(\Phi_{1}(g_{1}) \otimes \Phi_{2}(g_{2})) = p^{a}$ by 
Lemma \ref{prod1}(i). Assume $\ell = p$. In this case $a = 1$, whence 
$\deg(\Theta(g)) \geq \deg(\Phi_{1}(g_{1}) \otimes \Phi_{2}(g_{2})) = p = o(g)$ by 
Lemma \ref{prod1}(ii). 

\smallskip
4b) Subcase II: $|g_{2}| < p^{a}$; in particular, $a \geq 2$. By Theorem \ref{clas1}(i), we have
two possibilities: $X = SU_{n_{1}}(q^{1/2})$ with $n_{1} \geq 3$ odd, or 
$X = Spin^{-}_{n_{1}}(q)$ with $2 < n_{1} \equiv 2 (\mod 4)$. Suppose the first possibility occurs.
Embed $g' := \diag(g_{1},g_{2})$ in a subgroup $Y := GU_{2}(q^{n_{1}/2})$ of $GU_{n}(q^{1/2})$. The 
assumption $|g_{2}| < p^{a}$ implies that $o(g') (\mbox{in } Y) = |g'| = p^{a}$. Applying Lemmas 
\ref{sl21} and \ref{su2} to an irreducible constituent $\Psi$ of degree $> 1$ of 
$\Theta|_{Y \cap G}$, we see that $d_{\Theta}(g) \geq d_{\Psi}(g') \geq p^{a}-2 > p^{a-1}(p-1)$,
a contradiction. 

Consider the second possibility. Then by Proposition \ref{icl}(iv) we can embed $g_{1}$, resp. 
$g_{2}$ in a subgroup $SU_{n_{1}/2}(q)$, and then embed the above 
element $g'$ in a subgroup $Y := SU_{n_{1}}(q)$ of $G$. Again we have  
$o(g') (\mbox{in } Y) = |g'| = p^{a}$. By our assumptions, there is an irreducible 
constituent $\Psi$ of degree $> 1$ of $\Theta|_{Y}$ such that $d_{\Psi}(g') \leq p^{a-1}(p-1)$. 
Notice that $n_{1} \geq 6$ here. Now if $k = 2$, that is $p|(q+1)$, then we get a contradiction
by the results of 2). So $k > 2$; in particular, $(n_{1},q) \neq (6,2)$. Hence Theorem \ref{clas1}
applies to $X$ and (\ref{cond}) holds as well. In this case, we get a contradiction by the 
$SU$-part of 4b).

\smallskip   
5) In what follows, we will examine the cases where either Theorem \ref{clas1} does not apply or
(\ref{cond}) fails.
  
5a) In the $SU$-case, $n_{1} \geq 3$ (as $g_{1}$ is a nontrivial irreducible $p$-element)
and $(n_{1},q) \neq (3,4)$ (as $SU_{3}(2)$ does not contain any irreducible $3'$-elements of 
$GL_{3}(4)$). Thus Theorem \ref{clas1} applies and (\ref{cond}) holds, and so we are done.

\smallskip
5b) Assume $S = PSp_{n}(q)$ with $n \geq 4$. Notice that $(n_{1},q) \neq (2,3)$ as $k > 1$; 
furthermore, if $(n_{1},q) \neq (2,2)$ then (\ref{cond}) holds. Applying Theorem \ref{clas1} to 
$X = Sp(V_{1})$, we are done if $n_{1} > 2$, or if $n_{1} = 2$ but $q$ is odd. Consider the case 
where $n_{1} = 2$ and $2|q \geq 4$. Then $n_{2} = 2$ (as $n \geq 4$) and $|g_{2}| = p$ (as 
otherwise $g$ fixes a nonzero singular vector). By Lemma \ref{sl21} we are done unless 
$p = |g_{1}| = q+1$. In the exceptional case, by Corollary \ref{prod4} we can find an 
irreducible constituent $\Phi_{1} \otimes \Phi_{2}$ of $\Theta|_{X_{1} \times X_{2}}$ such 
that $\Phi_{i} \in \IBRL(Sp(V_{i}))$ and $\dim(\Phi_{i})) > 1$. Again by Lemma \ref{sl21}, 
$\deg(\Phi_{i}(g_{i}) \geq p-2$, whence 
$\deg(\Theta(g)) \geq \deg(\Phi_{1}(g_{1}) \otimes \Phi_{2}(g_{2})) = p = o(g)$ by 
Lemma \ref{prod1}. Finally, assume $(n_{1},q) = (2,2)$. Then $p = 3$, and $g_{1}$ and $g_{2}$ 
are conjugate in $Sp_{2}(q)$. It follows by Lemma \ref{equal} that $g$ fixes a nonzero totally 
singular subspace of $V_{1} \oplus V_{2}$, a contradiction. 

\smallskip    
5c) Assume $S = P\Omega^{\eps}_{n}(q)$ with $n \geq 7$. By Proposition \ref{icl}(iv), either 
$\eps_{j} = -$ and $n_{j}$ is even, or $n_{j} = 1$. First we consider the case $n_{1} > 4$ 
and $(n_{1},q) \neq (6,2)$. Then (\ref{cond}) holds. Observe that Theorem \ref{clas1} can be applied 
to the element $g_{1}$ of $X = Spin(V_{1})$. (This is clear if $n_{1} \geq 8$. Suppose 
$n_{1} = 6$. Then $Spin^{-}_{6}(q) \simeq SU_{4}(q)$, and $o(g)$ divides $q^{3}+1$ but not
$q+1$ (by irreducibility of $g_{1}$) by Proposition \ref{icl}(iv). Now we can apply the already 
proven results for $X \simeq SU_{4}(q)$ and obtain that $\deg(\Psi(g_{1})) \geq p^{a-1}(p-1)$ for any 
$\Psi \in \IBRL(X)$ of degree $> 1$ and that the $p$-Sylow subgroups of $X$ are cyclic. Moreover, 
if the equality attains, then either $a = 1$ or $\ell = p$.) Thus we are done in this case.

Now we assume that $(n_{1},q) = (6,2)$; in particular, $o(g) = o(g_{1}) = 9$. Claim that 
$\deg(\Theta(g)) \geq 7$. For, $g^{3}$ belongs to a parabolic subgroup of $G$ by Proposition 
\ref{icl}(iv). Hence by the main result of \cite{DZ1} $\deg(\Theta(g^{3})) = o(g^{3}) = 3$. 
Therefore if $\ell = 3$ then we are done by Lemma \ref{root1}(ii). It remains to consider 
the case $\ell \neq 3$. First we suppose that $n_{2} = 6$. Then $|g_{2}| = 9$ (as otherwise 
$g_{2}$ is reducible on $V_{2}$). It suffices to prove the claim for $G = \Omega^{+}_{12}(2)$.
Notice that $G > (\Omega(V_{1}) \times \Omega(V_{2})) \cdot \la t \ra$ for some involution
$t$ and $\la \Omega(V_{i}),t \ra \simeq GO(V_{i})$ for $i = 1,2$. Furthermore, $g_{i}$ is 
rational in $GO(V_{i})$, cf. \cite{Atlas}. Thus $g = \diag(g_{1},g_{2})$ is rational in $G$. 
Since $|g| = 9$, the rationality of $g$ combined with $\deg(\Theta(g^{3})) = 3$ implies that 
$\deg(\Theta(g) \geq 7$. Finally, we consider the case $n_{2} < 6$. Since $q = 2$ and 
$n > 6$, we must have $n_{2} = 2$ or $4$ (and $\eps_{2} = -$). But $p = 3$ and the $3$-element 
$g_{2}$ is irreducible on $V_{2}$, so $n_{2} = 2$ and $|g_{2}| = 3$. It now suffices to prove
the claim for $G = \Omega^{+}_{8}(2)$. Notice that $\Omega^{+}_{8}(2)$ has three classes of 
elements of order $9$, which are permuted by the triality automorphism $\tau$, and at least one of 
them intersects a subgroup $Y \simeq Sp_{6}(2)$ of $G$, cf. \cite{Atlas}. So we may assume that 
$g' := g^{\tau^{l}} \in Y$ for some $\ell = 0,1,2$. Clearly, $g'$ is irreducible in $Y$. 
Applying Proposition \ref{sp1} to $g'$ and an irreducible constituent of $\Theta^{\tau^{-l}}|_{Y}$, 
we see that $\deg(\Theta(g)) = \deg(\Theta^{\tau^{-l}}(g^{\tau^{l}})) \geq 7$, and so we are done.
      
\smallskip
5d) We continue the case $S = P\Omega^{\eps}_{n}(q)$ with $n \geq 7$ and $n_{1} \leq 4$.
Consider the case $n_{1} = 4$ and $2|q$. First suppose that $q \geq 4$. Then (\ref{cond}) holds. 
Notice that $\Omega^{-}_{4}(q) \simeq PSL_{2}(q^{2})$. Now we can apply Lemma \ref{sl21} and argue 
as in 5b). Assume $q = 2$, whence $o(g) = 5$. Since $g$ is not contained in any parabolic 
subgroup of $G$, $n_{2} = 4$. By restriction, it suffices to prove 
$\deg(\Theta(g)) = 5$ for $n = 8$. Notice that $\Omega^{+}_{8}(2)$ has three classes of elements of 
order $5$, which are permuted by the triality automorphism $\tau$, and one of them intersects a 
parabolic subgroup of $G$, cf. \cite{Atlas}. So we may assume that $g^{\tau^{l}}$ is contained in a
parabolic subgroup for some $\ell = 0,1,2$. By \cite{DZ1} applied to $\Theta^{\tau^{-l}}$, 
$\deg(\Theta(g)) = \deg(\Theta^{\tau^{-l}}(g^{\tau^{l}})) = 5$, and so we are done.      

Consider the case $n_{1} = 4$ and $q$ is odd. First suppose that $q \geq 5$. Then (\ref{cond}) holds. 
Notice that $Spin^{-}_{4}(q) \simeq SL_{2}(q^{2})$. Now we can apply Lemma \ref{sl21} and argue 
as in 5b). Assume $q = 3$, whence $o(g) = 5$. Since $g$ is not contained in any parabolic 
subgroup of $G$, $n_{2} = 4$. By restriction, it suffices to prove 
$\deg(\Theta(g)) = 5$ for $n = 8$. Notice that $Spin^{+}_{8}(3)$ has three classes of elements of 
order $5$, which are permuted by the triality automorphism $\tau$, and one of them intersects a 
parabolic subgroup of $G$, cf. \cite{Atlas}. So we may assume that $g^{\tau^{l}}$ is contained in a
parabolic subgroup for some $\ell = 0,1,2$. By \cite{DZ1} applied to $\Theta^{\tau^{-l}}$, 
$\deg(\Theta(g)) = \deg(\Theta^{\tau^{-l}}(g^{\tau^{l}})) = 5$, and so we are done.
  
Now we may assume that $n_{1} = 2$. If $m$ is the largest index such that $n_{1} = n_{m}$, then 
$n \leq 2m+2$ (otherwise $g$ acts trivially on the nondegenerate subspace 
$V_{m+1} \oplus \ldots \oplus V_{m+s}$ of dimension $\geq 3$ and so $g$ fixes a nonzero singular 
vector). In particular, $m \geq 3$. Recall that $g_{i}^{q+1} = 1$ for all $i$, and 
$o(g_{1}) = p^{a}$. Now we can use the isomorphism $P\Omega^{-}_{6}(q) \simeq PSU_{4}(q)$ and apply 
Corollary \ref{su32} to the element $h := \diag(g_{1},g_{2},g_{3})$ (inside an inverse image of 
$P\Omega^{-}_{6}(q)$ in $G$) to get that $q+1 = p^{a}$ and $h$ is a pseudoreflection in 
$GU_{4}(q)$. Notice that the isomorphism $P\Omega^{-}_{6}(q) \simeq PSU_{4}(q)$ is realized by 
letting $SU_{4}(q)$ act on the alternating square of its $4$-dimensional module, cf. \cite[p. 45]{KL}. 
It follows that the spectrum of $h$ on the $6$-dimensional module for $\Omega^{-}_{6}(q)$ is 
of the form $\{\beta,\beta,\beta,\beta^{-1},\beta^{-1},\beta^{-1}\}$ for some $\beta \in \FF_{q^{2}}$. 
In particular, the $g_{i}$-modules $V_{i}$ are isomorphic for $i = 1,2,3$.
By Lemma \ref{equal}, $h$ belongs to a parabolic subgroup of $GO^{-}_{6}(q)$, contrary to 
our assumption. 

\smallskip
6) Finally, we come back to the exception specified before (\ref{repl}). Recall that, in this case
$n = 2m$, $S = \Omega^{\eps}_{n}(q)$ and $g \in H = Sp_{n-2}(q)$. By the result proved for $H$, we 
see that, under the assumption $\deg(\Theta(g)) \leq p^{a-1}(p-1)$, $p$ is odd and Sylow 
$p$-subgroups of $H$ are cyclic (in particular, $k \geq 2 \lfloor (m-1)/2 \rfloor + 1$ by 
Lemma \ref{sylow}); furthermore, either $a = 1$ or $\ell \neq p$. If Sylow $p$-subgroups of $G$ are 
cyclic then we are done. Assume the contrary. Then $\eps = +$, $k = m$ is even, and we can 
embed $g$ in a subgroup $R \simeq Sp_{m}(q)$ of $H$. Applying Proposition \ref{sp1} to an 
irreducible constituent of degree $> 1$ of $\Theta|_{R}$, we see that $(m,q,|g|) = (4,2,5)$. 
Using the argument with the triality automorphism of $\Omega^{+}_{8}(2)$ as in 5d), we conclude
that $\deg(\Theta(g)) = 5 = o(g)$.         
\end{proof}

We can say more about case (ii) of Theorems \ref{main1}, \ref{main2}, and \ref{clas2}. 

\begin{propo}\label{p-cyclic}
{\sl Under the assumptions of Theorem \ref{main2}, suppose that 
$\deg(\Theta(g)) = p^{a-1}(p-1)$ and Sylow $p$-subgroups of $G/Z(G)$ are cyclic. 
Let $m$ be the smallest positive integer such that $p|(q^{m}-1)$. Then one of the following 
holds.

{\rm (i)} $S = PSL_{n}(q)$, $n \geq 3$, and $m = n$. Moreover, if $2|n$ then 
$p^{a} = p = (q^{n/2}+1)/(2,q+1)$.

{\rm (ii)} $S = PSU_{n}(q)$, $n \geq 3$, and $m = 4\lfloor (n-1)/2 \rfloor + 2$.

{\rm (iii)} $S = PSp_{2n}(q)$, $n \geq 1$, $m = 2n$, and  
$p^{a} = p = (q^{n}+1)/(2,q+1)$.

{\rm (iv)} $S = P\Omega^{+}_{2n}(q)$, $n \geq 4$, and $m = 2n-2$. Moreover, if $n$ is odd then 
$p^{a} = p = (q^{n-1}+1)/(2,q+1)$.

{\rm (v)} $S = P\Omega^{-}_{2n}(q)$, $n \geq 4$, and $m = 2n$. Moreover, if $2|n$ then 
$p^{a} = p = (q^{n}+1)/(2,q+1)$.}
\end{propo}

\begin{proof}
1) Our choice of $m$ is equivalent to the condition that $p$ divides $\Phi_{m}(q)$ but not  
$\Phi_{i}(q)$ for any $i < m$. By the main result of \cite{DZ1}, the conditions that
$\deg(\Theta(g)) = p^{a-1}(p-1)$ and Sylow $p$-subgroups of $G/Z(G)$ are cyclic imply that
$g$ cannot belong to any parabolic subgroup of $G$. 

2) Assume $S = PSL_{n}(q)$ with $n \geq 3$. By assumption, Sylow $p$-subgroups of $S/Z(S)$ are 
cyclic, whence $(p,q-1) = 1$. In this case, we see that Sylow $p$-subgroups of $GL_{n}(q)$ 
are cyclic, so without loss we may replace $G$ by $GL_{n}(q)$. Claim that $p$ does not divide
the order of any parabolic subgroup $P$ of $G$. Assume the contrary. Write 
$|P|$ and $(G:P)$ (in a unique way, as cyclotomic polynomials are irreducible
over $\QQ$) as $q^{r}\prod_{i \in I}\Phi^{k_{i}}_{i}(q)$ and 
$q^{s}\prod_{j \in J}\Phi^{l_{j}}_{j}(q)$. By Lemma \ref{sylow} and since $p$ divides $|P|$, 
$m \in I$ and $p \not{|} \Phi_{j}(q)$ for all $j \in J$. Thus $(G:P)$ is coprime to $p$, and
so a conjugate of the $p$-element $g$ belongs to $P$, a contradiction. Now our claim 
implies that $m = n$. In particular, $g$ is irreducible, and we can apply Corollary 
\ref{sp-gl} if $2|n$.

3) Assume $S = PSU_{n}(q)$ with $n \geq 3$. Arguing as in 2), we see that $(p,q+1) = 1$, and
that $p$ does not divide the order of any parabolic subgroup $P$ of $GU_{n}(q)$. Since 
we can embed $X := GU_{n-2}(q)$ in such a $P$, we see that $(p,|X|) = 1$. Furthermore, since 
$Y := GL_{\lfloor n/2 \rfloor}(q^{2})$ embeds in another such a $P$, we see that 
$(p,|Y|) = 1$. The last two conditions readily imply that $m = 2n-2$ for even $n$ and 
$m = 2n$ for odd $n$, and so we are done.

4) Assume $S = PSp_{2n}(q)$ with $n \geq 2$. Arguing as in 3), we see that $p > 2$, and
that $p$ is coprime to $|Sp_{2n-2}(q)|$ and $|GL_{n}(q)|$. It follows that $m = 2n$. In 
particular, $g$ is irreducible, and we can apply Proposition \ref{sp1}.

Assume $S = P\Omega^{+}_{2n}(q)$ with $n \geq 4$. Arguing as in 3), we see that $p > 2$, and
that $p$ is coprime to $|SO^{+}_{2n-2}(q)|$ and $|GL_{n}(q)|$. It follows that $m = 2n-2$. 
If $n$ is odd (so $n \geq 5$), then we can embed $g$ in $Spin^{-}_{2n-2}(q)$ as an 
irreducible element and apply Lemma \ref{so1}.

Assume $S = P\Omega^{-}_{2n}(q)$ with $n \geq 4$. Arguing as in 3), we see that $p > 2$, and
that $p$ is coprime to $|SO^{-}_{2n-2}(q)|$ and $|GL_{n-1}(q)|$. It follows that $m = 2n$
and so $g$ is irreducible. If $n$ is even, then we can apply Lemma \ref{so1}.   

Assume $S = PSL_{2}(q)$. Arguing as in 2), we see that $p > 2$, and that $(p,q-1) = 1$. 
It follows that $m = 2$, and we can apply Lemma \ref{sl21}.
\end{proof}

Now Theorems \ref{main1} and \ref{main2} follow immediately from Theorems \ref{clas1} and \ref{clas2}.

\end{document}